\documentclass[11pt,twoside]{article}
\newcommand{\ifims}[2]{#1} % 1 - Latex, 2 - ims version
\newcommand{\ifAoS}[2]{#1}   % 1 - AoS, 2 - full
\newcommand{\ifAMS}[2]{#1}   % 1 - AMS, 2 - JEL
\newcommand{\ifau}[3]{#2}  % #1 - 1 author, #2 - 2 authors, #3 - 3 authors

\newcommand{\ifbook}[2]{#1}   % 1 - paper; 2 - book

\def\thetitle{Bootstrap tuning in ordered model selection}
\def\thanksa%{}
{The author is supported by RFBR, research project No. 13-01-12007 ofi\_m.
Financial support by the German Research Foundation (DFG) through the Collaborative 
Research Center 649 ``Economic Risk'' is gratefully acknowledged
} 
\def\theruntitle{Bootstrap-tuned model selection}
\def\thanksb{}
\def\theabstract{
In the problem of model selection for a given family of linear estimators,
ordered by their variance,
we offer a new ``smallest accepted'' approach motivated by Lepski's method 
and multiple testing theory.
The procedure selects the smallest model which satisfies an acceptance rule based 
on comparison with all larger models.
The method is completely data-driven and does not use any prior information about 
the variance structure of the noise: 
its parameters are adjusted to the underlying possibly heterogeneous noise
by the so-called ``propagation condition'' using a wild bootstrap 
method.
The validity of the bootstrap calibration is proved for finite samples with an explicit error bound.
We provide a comprehensive theoretical study of the method and 
describe in detail the set of possible values of the selector \( \hat{\mm} \).
We also establish some precise 
oracle error bounds for the corresponding estimator 
\( \hat{\thetav} = \tilde{\thetav}_{\hat{\mm}} \)
which equally applies to estimation of the whole parameter vectors, 
some subvector or linear mapping, as well as the estimation of a linear functional. 
}

\def\kwdp{62G05}
\def\kwds{62G09, 62J15}

\def\thekeywords{smallest accepted, oracle, propagation condition}

\def\authora{Vladimir Spokoiny}
\def\runauthora{spokoiny, v.}
\def\addressa{
    Weierstrass-Institute and \\
    Humboldt University Berlin, \\ 
    HSE and IITP RAS, Moscow, 
    \\
    Mohrenstr. 39, \\
    10117 Berlin, Germany,   
    }
\def\emaila{spokoiny@wias-berlin.de}
\def\affiliationa{Weierstrass-Institute and Humboldt University Berlin}
\def\authorb{Niklas Willrich}
\def\runauthorb{willrich, n.}
\def\addressb{
    \\
    Weierstrass-Institute ,     
    \\
    \\
    Mohrenstr. 39, \\
    10117 Berlin, Germany,    
    }
\def\emailb{willrich@wias-berlin.de}
\def\affiliationb{Weierstrass-Institute Berlin}

\renewcommand{\(}{$\,}
\renewcommand{\)}{\,$}

\def\nquad{\hspace{-1cm}}
\def\eqdef{\stackrel{\operatorname{def}}{=}}

\def\ND{\mathcal{N}}

\newcommand{\cc}[1]{\mathscr{#1}}
\newcommand{\bb}[1]{\boldsymbol{#1}}

\renewcommand{\bar}[1]{\overline{#1}}
\renewcommand{\hat}[1]{\widehat{#1}}
\renewcommand{\tilde}[1]{\widetilde{#1}}

\renewcommand{\Gamma}{\varGamma}
\renewcommand{\Pi}{\varPi}
\renewcommand{\Sigma}{\varSigma}
\renewcommand{\Delta}{\varDelta}
\renewcommand{\Lambda}{\varLambda}
\renewcommand{\Psi}{\varPsi}
\renewcommand{\Phi}{\varPhi}
\renewcommand{\Theta}{\varTheta}
\renewcommand{\Omega}{\varOmega}
\renewcommand{\Xi}{\varXi}
\renewcommand{\Upsilon}{\varUpsilon}

\def\Var{\operatorname{Var}}

\def\argmin{\operatornamewithlimits{argmin}}

\def\tr{\operatorname{tr}}

\def\R{I\!\!R}
\def\E{I\!\!E}
\def\P{I\!\!P}

\def\kappa{\varkappa}

\def\T{\top}
\def\diag{\operatorname{diag}}

\def\rank{\operatorname{rank}}

\def\av{\bb{a}}
\def\bv{\bb{b}}

\def\fv{\bb{f}}

\def\uv{\bb{u}}

\def\wv{\bb{w}}
\def\xv{\bb{x}}

\def\zv{\bb{z}}

\def\Cv{\bb{C}}

\def\Mv{\bb{M}}
\def\Sv{\bb{S}}
\def\Uv{\bb{U}}
\def\Xv{\bb{X}}
\def\Yv{\bb{Y}}
\def\Zv{\bb{Z}}

\def\betav{\bb{\beta}}
\def\gammav{\bb{\gamma}}

\def\epsv{\bb{\varepsilon}}
\def\etav{\bb{\eta}}
\def\gammav{\bb{\gamma}}

\def\xiv{\bb{\xi}}

\def\CONST{\mathtt{C}}

\newcommand{\tobedone}[1]{\par\textbf{\color{red}To be done:} {\color{magenta}#1}}

\usepackage{color}

\definecolor{blue(pigment)}{rgb}{0.2, 0.2, 0.6}
\definecolor{ultramarine}{rgb}{0.07, 0.04, 0.56}
\definecolor{darkspringgreen}{rgb}{0.09, 0.45, 0.27}
\definecolor{hookersgreen}{rgb}{0.0, 0.44, 0.0}
\definecolor{plum(traditional)}{rgb}{0.56, 0.27, 0.52}
\definecolor{purple(html/css)}{rgb}{0.5, 0.0, 0.5}
\definecolor{magenta(dye)}{rgb}{0.79, 0.08, 0.48}

\def\Frob{\operatorname{Fr}}
\def\TV{\operatorname{TV}}
\def\oper{\operatorname{op}}
\def\Fr{\operatorname{Fr}}

\def\CONST{\mathtt{C} \hspace{0.1em}}
\def\cond{\, \big| \,}

\def\nsize{{n}}

\def\ex{\mathrm{e}}

\def\ND{\cc{N}}
\def\gaussv{\bb{\gauss}}
\def\gauss{\gamma}

\def\Id{I\!\!\!I}
\def\Ind{\operatorname{1}\hspace{-4.3pt}\operatorname{I}}

\def\alp{\alpha}

\def\alpb{\alp_{1}}		% conflict
\def\alpb{\alp_{+}}
\def\alpb{\beta}		% conflict

\def\alpd{\bar{\alp}}

\def\AA{A}
\def\AAb{\AA^{\sbt}}

\def\bias{b}
\def\biasv{\bb{b}}

\def\Bias{B}		% conflict
\def\Bias{\bb{B}}

\def\BB{I\!\!B}

\def\BB{B}

\def\BBB{\cc{B}}

\def\bxiv{\xiv^{\sbt}}

\def\bvb{\bv^{\sbt}}
\def\bvs{\bv}

\def\Covm{\Sigma}			% conflict
\def\Covm{\mathbb{V}}
\def\covm{\sigma}
\def\dimp{p}

\def\dimA{\mathtt{p}}

\def\dimn{\dimp_{\nsize}}
\def\dimq{q}

\def\dimAb{\dimA^{\sbt}}

\def\dimAbt{\dimA^{\sbt}}

\def\dimw{\mathtt{q}}
\def\dimwf{\dimw_{2}}

\def\dPsi{\td_{\Psi}}

\def\dPsis{\td^{*}}
\def\dPsip{\dPsi}

\def\Eb{\E^{\sbt}}

\def\errSi{\errXY}

\def\eps{\epsilon}			% conflict
\def\eps{\varepsilon}
\def\epsvb{\epsv^{\sbt}}
\def\epsvr{\breve{\epsv}}
\def\epsr{\breve{\eps}}

\def\errSi{\Delta}
\def\errV{\errSi_{0}}		% conflict
\def\fvs{\fv}
\def\fs{f}

\def\fvs{\fv^{*}}

\def\Gam{\Xi}
\def\Gam{\mathcal{S}}

\def\GQF{Q}

\def\IF{\Bbb{F}}

\def\IFr{\breve{\IF}}

\def\kullb{\cc{K}} %{\wp}

\def\LL{\cc{L}}

\def\losst{\varrho}

\def\lambdav{\bb{\lambda}}

\def\mm{m}
\def\ms{\mm^{*}}

\def\mc{\mm'}

\def\md{\mm^{\circ}}

\def\mmmax{\mm_{\max}}	% conflict
\def\mmmax{M}
\def\mmmin{\mm_{0}}

\def\mmprec{\mm_{(-1)}}
\def\mres{\mm^{\dag}}
\def\MM{\cc{M}}
\def\MMu{\MM^{+}}
\def\MMd{\MM^{-}}
\def\MMi{\MM^{\circ}}
\def\MMc{\MM^{c}}

\def\Meta{\mathbb{M}}

\def\Pb{\P^{\sbt}}

\def\PPsi{\mathbb{Q}}
\def\PPsib{\mathbb{Q}^{\sbt}}

\def\QL{W}

\def\qq{q}

\def\qqbt{\qq^{\sbt}}

\def\rdl{\epsilon}
\def\rd{\bb{\rdl}}

\def\riskt{\cc{R}}

\def\riskr{\riskt^{+}}

\def\rr{\mathtt{r}}

\def\rei{\delta}

\def\SPiS{\Upsilon}

\def\Sigmar{\breve{\Sigma}}
\def\Sam{\cc{T}}

\def\score{\nabla}

\def\scoreb{\score^{\sbt}}

\def\sbt{\hspace{1pt} \flat}

\def\supA{\lambda}
\def\supAB{\supA^{*}}

\def\supepsi{\delta_{\eps}}
\def\supeps{\delta_{1}}

\def\thetav{\bb{\theta}}
\def\thetavs{\thetav^{*}}

\def\thetas{\theta^{*}}

\def\test{\tau}

\def\td{\delta}
\def\tdn{\td_{n}}

\def\tdPi{\tdn}

\def\tar{\phi}
\def\tars{\tar^{*}}
\def\tarv{\bb{\tar}}
\def\tarvs{\tarv^{*}}

\def\Tam{\cc{K}}

\def\Tmd{\mathbb{T}}
\def\Tmdb{\Tmd^{\sbt}}

\def\Tb{\Tmd^{\sbt}}

\def\upsv{\bb{\varkappa}}

\def\UV{\mathcal{U}}
\def\UVcol{\bb{\omega}}

\def\vA{\mathtt{v}}

\def\Varb{\Var^{\sbt}}
\def\Varxi{S}
\def\Varxib{S^{\sbt}}

\def\VP{V}
\def\VPb{\cc{V}}	% conflict

\def\VPb{{\VP^{\sbt}}}

\def\VQ{\mathbb{V}}

\def\vp{\mathbf{v}}	% conflict

\def\vp{\mathrm{v}}

\def\vpi{\mathtt{b}}

\def\wv{\bb{w}}

\def\wb{w^{\sbt}}		% conflict
\def\Wb{\uv^{\sbt}}		% conflict
\def\WB{\Uv^{\sbt}}		% conflict
\def\Wb{\wv^{\sbt}}
\def\WB{\cc{W}^{\sbt}}

\def\xivb{\xiv^{\sbt}}

\def\xx{\mathtt{x}}

\def\xxn{\xx_{\nsize}}
\def\xxp{\xx_{\dimp}}

\def\xxt{t}
\def\xxb{\xx^{\sbt}}
\def\xxu{\xx_{\mathtt{s}}}

\def\Yr{\breve{Y}}
\def\Yvr{\breve{\Yv}}

\def\zq{z}

\def\zqb{\bar{\zq}}

\def\zqb{\zq^{\sbt}}

\def\ZZ{\mathtt{z}}		% conflict
\def\ZZ{\mathfrak{z}}

\def\ZZs{\bar{\ZZ}}

\def\ZZbt{\ZZ^{\sbt}}

\def\mm{m}
\usepackage{amsmath,amssymb,amsthm}
\usepackage{natbib}
\usepackage{epsfig,graphicx}
\usepackage{comment}
\usepackage{color}
\usepackage{srcltx}
\usepackage[mathscr]{eucal}
\usepackage[math]{easyeqn}
\usepackage{etoolbox}
\usepackage{hyperref}
\hypersetup{%linktocpage,
            colorlinks,
            linkcolor=hookersgreen,
            linktoc=blue,
            citecolor=ultramarine,
            urlcolor=black,
            filecolor=black
            }
\usepackage{nameref}

\ifims{
\textheight=23cm
\textwidth=14.8cm
\topmargin=0pt
\oddsidemargin=1.0cm
\evensidemargin=1.0cm
\linespread{1.3}
\renewenvironment{abstract}
    {\centerline{\textbf{Abstract}}\bigskip
      \begin{center}
       \begin{minipage}{11cm}
        \begin{small}
    }
    {   \end{small}
       \end{minipage}
      \end{center}
     \bigskip
    }

}{ %ims style
}

\numberwithin{equation}{section}
\numberwithin{figure}{section}
\newcounter{example}[section]
\numberwithin{example}{section}
\newcounter{remark}[section]
\numberwithin{remark}{section}
\newtheorem{theorem}{Theorem}[section]
\newtheorem{proposition}[theorem]{Proposition}
\newtheorem{lemma}[theorem]{Lemma}
\newtheorem{corollary}[theorem]{Corollary}

\newtheorem{exmp}[example]{Example}
\newtheorem{rmrk}[remark]{Remark}
\newenvironment{example}{\begin{exmp}\rm}{\end{exmp}}
\newenvironment{remark}{\begin{rmrk}\rm}{\end{rmrk}}

\begin{document}
\thispagestyle{empty}
\ifims{
\title{\thetitle}
\ifau{ % 1 author
  \author{
    \authora
    \ifdef{\thanksa}{\thanks{\thanksa}}{}
    \\[5.pt]
    \addressa \\
    \texttt{ \emaila}
  }
}
{  % 2 authors
  \author{
    \authora
    \ifdef{\thanksa}{\thanks{\thanksa}}{}
    \\[5.pt]
    \addressa \\
    \texttt{ \emaila}
    \and
    \authorb
    \ifdef{\thanksb}{\thanks{\thanksb}}{}
    \\[5.pt]
    \addressb \\
    \texttt{ \emailb}
  }
}
{   % 3 authors
  \author{
    \authora
    \ifdef{\thanksa}{\thanks{\thanksa}}{}
    \\[5.pt]
    \addressa \\
    \texttt{ \emaila}
    \and
    \authorb
    \ifdef{\thanksb}{\thanks{\thanksb}}{}
    \\[5.pt]
    \addressb \\
    \texttt{ \emailb}
    \and
    \authorc
    \ifdef{\thanksc}{\thanks{\thanksc}}{}
    \\[5.pt]
    \addressc \\
    \texttt{ \emailc}
  }
}

\maketitle
\pagestyle{myheadings}
\markboth
 {\hfill \textsc{ \small \theruntitle} \hfill}
 {\hfill
 \textsc{ \small
 \ifau{\runauthora}
      {\runauthora and \runauthorb}
      {\runauthora, \runauthorb, and \runauthorc}
 }
 \hfill}
\begin{abstract}
\theabstract
\end{abstract}

\ifAMS
    {\par\noindent\emph{AMS 2000 Subject Classification:} Primary \kwdp. Secondary \kwds}
    {\par\noindent\emph{JEL codes}: \kwdp}

} % end front latex
{ % front ims
\begin{frontmatter}
\title{\thetitle}
%\thankstext{T1}{\thankstitle}

% "Title of the paper"

\runtitle{\theruntitle}

\ifau{ % 1 author
\begin{aug}
    \author{\authora\ead[label=e1]{\emaila}}
    \address{\addressa \\
     \printead{e1}}
\end{aug}

 \runauthor{\runauthora}
\affiliation{\affiliationa} }
{ % 2 authors
\begin{aug}
    \author{\authora\ead[label=e1]{\emaila}\thanksref{t21}}
    \and
    \author{\authorb\ead[label=e2]{\emailb}\thanksref{t22}}
    
    \address{\addressa \\
     \printead{e1}}
    \address{\addressb \\
     \printead{e2}}
    \thankstext{t21}{\thanksa}
    \thankstext{t22}{\thanksb}
    \affiliation{\affiliationa, \affiliationb} 
    \runauthor{\runauthora and \runauthorb}
\end{aug}
} 
{ % 3 authors
\begin{aug}
    \author{\authora\ead[label=e1]{\emaila}\thanksref{t21}}
    \and
    \author{\authorb\ead[label=e2]{\emailb}\thanksref{t22}}
    \and
    \author{\authorc\ead[label=e3]{\emailc}\thanksref{t23}}
    
    \address{\addressa \\
     \printead{e1}}
    \address{\addressb \\
     \printead{e2}}
    \address{\addressc \\
     \printead{e3}}
    \thankstext{t21}{\thanksa}
    \thankstext{t22}{\thanksb}
    \thankstext{t23}{\thanksc}
    \affiliation{\affiliationa, \affiliationb, \affiliationc} 
    \runauthor{\runauthora, \runauthorb, and \runauthorc}
\end{aug}}

\begin{abstract}
\theabstract
\end{abstract}



\end{frontmatter}
} % end front ims

\ifims{ %full
\tableofcontents
}{ %annals
}
% !TEX root = modsel_2015_7.tex

\def\suppSW{the supplement [SW2015]}

\section{Introduction}
Model selection is one of the key topics in mathematical statistics. 
A choice between models of differing complexity can often be viewed as a trade-off between
overfitting the data by choosing a model which has too many degrees of freedom
and smoothing out the underlying structure in the data by choosing a model which
has too few degrees of freedom. 
This trade-off which shows up in most methods
as the classical bias-variance trade-off is at the heart of every model
selection method (as for example in 
unbiased risk estimation, \cite{Kneip1994} or in
penalized model selection,
\cite{barronmassartbirge1999},
\cite{massartconcentration2007}).
This is also the case in Lepski's method, \cite{lepski1990}, \cite{lepski1991},
\cite{lepski1992}, \cite{lepskispokoiny1997}, \cite{lepski1997},
\cite{birgelepski} and risk hull minimization, \cite{cavaliergolubev2006}. 
Many of these methods allow their strongest theoretical results only for highly
idealized situations (for example sequence space models), are very specific to
the type of problem under consideration (for instance, signal or functional estimation),
require to know the noise behavior (like homogeneity) and the exact noise level.
Moreover, they typically involve an unwieldy number of calibration constants whose
choice is crucial to the applicability of the method and is not addressed by the theoretical 
considerations.
For instance, any Lepski-type method requires to fix a numerical constant in the definition 
of the threshold,
the theoretical results only apply if this constant is sufficiently large while 
the numerical results benefit from the choice of a rather small constant.   
\cite{spokoinyvial} offered a propagation approach to calibration of 
Lepski's method in the case of the estimation of a one-dimensional quantity of interest.
However, the proposal still requires the exact knowledge of the noise level 
and only applies to linear functional estimation.
A similar approach has been applied to local constant density estimation with sup-norm risk in
\cite{GaNiSp2012} and to local quantile estimation in \cite{SpWe2012}. 

In the case of unknown but homogeneous noise, generalized cross validation can be used instead of 
unbiased risk estimation method. 
For the penalized model selection, recently a number of proposals appeared
to apply one or another resampling method. 
\cite{arlot2009} suggested the use of resampling methods for the choice of an optimal penalization, following the framework of  penalized model selection,
\cite{barronmassartbirge1999},
%\cite{massartconcentration2007}. 
%Another approach of data-driven calibration in the face of an unknown error structure is proposed in \cite{arlotminimal2009} using the concept of minimal penalties, 
\cite{massartminimalpenalties}.
The validity of a bootstrapping procedure for Lepski's method has also been studied in \cite{chernozhukov2014} 
with new innovative technical tools with applications to
honest adaptive confidence bands.

An alternative approach to adaptive estimation is based on 
aggregation of different estimates; see \cite{goldenshluger2009} and \cite{Dala2012} for an overview 
of the existing results.
However, the proposed aggregation procedures either require  
two independent copies of the data or involves a data splitting for estimating the noise variance. 
Each of these requirements is very restrictive for practical applications.

Another point to mention is that the majority of the obtained results on adaptive estimation 
focus on the quality of estimating the unknown response, that is, the loss is measured 
by the difference between the true response and its estimate.
At the same time, inference questions like confidence estimation would require 
to know some additional information about the right model parameter. 
Only few results address the issue of estimating the true (oracle) model.
Moreover, there are some negative results showing that 
a construction of adaptive honest confidence sets is impossible without 
special conditions like self-similarity;
see, e.g. \cite{GiNi2010}.

This paper aims at developing a unified approach to the problem of ordered model selection 
with the focus on the quality of model selection rather than on accuracy of adaptive 
estimation under realistic assumptions on the model. 
Our setup covers linear regression and linear inverse problems, 
and equally applies to estimation of the whole parameter vectors, 
a subvector or linear mapping, as well as estimation of a linear functional. 
The proposed procedure and the theoretical study are also unified and do not distinguish 
between models and problems. 
In the case of a linear inverse problem, it is applicable to mild and severely ill-posed 
problems without prior knowledge of the type and degree of ill-posedness; cf. \cite{Tsybakov2000835},  \cite{cavalier2002}.
Another important issue is that the procedure does not use any prior information about 
the variance structure of the noise under assumption of minimal H\"older smoothness 1/4 on the underlying function. 
The method automatically adjusts the parameters to the underlying possibly heterogeneous noise:
the resampling technique allows to achieve the same quality of estimation as if the noise 
structure were precisely known.
Also we allow for a model misspecification: the linear structure of the response 
can be violated, in this case the procedure adaptively recovers the best linear projection.

\medskip

Consider a linear model \( \Yv = \Psi^{\T} \thetavs + \epsv \) in \( \R^{n} \) 
for an unknown parameter vector \( \thetavs \in \R^{\dimp} \) and 
a given \( \dimp \times n \) design matrix \( \Psi \).
Suppose that a family of linear smoothers \( \tilde{\thetav}_{\mm} = \Gam_{\mm} \Yv \) 
is given, where \( \Gam_{\mm} \) is for each \( \mm \in \MM \) a given \( \dimp \times n \) matrix.
We also assume that this family is ordered by the complexity of the method.
The task is to develop a data-based model selector \( \hat{\mm} \) which performs 
nearly as good as the optimal choice, which depends on the model and is not available.
The proposed procedure called the ``smallest accepted'' (SmA) rule can be viewed as 
a calibrated Lepski-type method.
The idea how the parameters of the method can be tuned,
originates from \cite{spokoinyvial} and is related to a multiple testing problem.
The whole procedure is based on family of pairwise tests,
each model is tested against all larger ones.
Finally the smallest accepted model is selected.
The critical values for this multiple testing procedure are fixed
using the so-called \emph{propagation condition}.
Theorem~\ref{ToracleSmA} presents finite sample results on the behavior of the proposed selector \( \hat{\mm} \) and 
the corresponding estimator \( \hat{\thetav} = \tilde{\thetav}_{\hat{\mm}} \). 
In particular, it describes a concentration set \( \MMi \) for the selected index 
\( \hat{\mm} \) and states an oracle bound for the resulting estimator 
\( \hat{\thetav} = \tilde{\thetav}_{\hat{\mm}} \).
Usual rate results can be easily derived from these statements.
Further results address the important issue called ``the payment for adaptation''
which can be defined as the gap between oracle and adaptive bounds.
Theorem~\ref{Tpayment} gives a general description of this quantity.
Then we specify the results to important special cases like projection estimation and 
estimation of a linear functional.
It appears, that in some cases the obtained results yield sharp asymptotic bounds.
In some other cases they lead to the usual log-price for adaptation. 
However, all these results require a known noise distribution.
Section~\ref{Sbootcalibr} explains how the proposed procedure can be tuned 
in the case of unknown noise using a bootstrap procedure. 
We establish explicit error bounds on the accuracy of the bootstrap approximation and show
that the procedure with 
bootstrap tuning does essentially the same job as the ideal procedure designed for the known noise. 
The study is quite involved because the procedure uses the same data twice 
for parameter tuning and for model selection.

The paper is structured as follows. 
The next section presents the procedure and the results for an idealistic situation
when the noise distribution is precisely known. 
We introduce the SmA selector \( \hat{\mm} \) and explain how it can be calibrated.
% to the problem at hand using the so called \emph{propagation condition}.
Then we describe the set of possible \( \hat{\mm} \)-values and establish probabilistic oracle 
bounds.
\ifAoS{% full 
Section~\ref{Spolynomloss} explains how the method and the results can be extended to the case 
of a polynomial loss function.}
{ % annals
An extension to the case of a polynomial loss function is given in Section A of \suppSW.}
The results are also specified to the particular problems of projection and linear 
functional estimation. 
Section~\ref{Sbootcalibr} extends the method and the study to the realistic case with 
unknown heteroscedastic noise by using a resampling technique.
The proofs and a detailed study of the bootstrap procedure in the linear Gaussian setup
are given in the appendix.
\ifAoS{% full
We also collect there some useful technical statements for Gaussian quadratic forms
and sums of random matrices.}
{%annals
The proofs of some technical results as well as some useful statements for Gaussian quadratic forms
and sums of random matrices are collected in \suppSW.}

\section{Model and problem. Known noise variance}
\label{SKnownnoice}
This section presents the model selector for the idealistic case when the noise distribution
is precisely known.
In the next section we explain how the unknown noise structure can be recovered from the data
using a resampling technique.
First we specify our setup. 
%Below we consider two models: the ``ideal'' linear Gaussian model will be used for 
%motivation and inspiration. 
%The real model will be used in our analysis.
%
We consider the following linear Gaussian model:
\begin{EQA}
	Y_{i}
	&=&
	\Psi_{i}^{\T} \thetavs + \eps_{i} \, ,
	\qquad
	\eps_{i} \sim \ND(0,\sigma^{2}) \, \, \text{ i.i.d. } ,
	\qquad
	i=1,\ldots,n,
\label{TiPiTtei}
\end{EQA}
with given design \( \Psi_{1},\ldots,\Psi_{n} \) in \( \R^{\dimp} \).
We also write this equation in vector form 
\( \Yv = \Psi^{\T} \thetavs + \epsv \in \R^{n} \), where \( \Psi \) is \( \dimp \times n \)
design matrix and \( \epsv \sim \ND(0,\sigma^{2} \Id_{n}) \).
Below we assume a deterministic design, otherwise one can understand the results 
conditioned on the design realization. 

In what follows, we allow the model \eqref{TiPiTtei} to be completely misspecified.
We mainly assume that the observations \( Y_{i} \) are independent and define
the response vector \( \fvs = \E \Yv \) with entries \( \fs_{i} \).
Such a model can be written as  
\begin{EQA}
	Y_{i} 
	&=&
	\fs_{i} + \eps_{i} \, .
\label{Yifieiiid}
\end{EQA}
Our study allows that
%Both kinds of misspecification are allowed.
the linear parametric assumption \( \fvs = \Psi^{\T} \thetavs \) is
violated, and the underlying noise \( \epsv = (\eps_{i}) \) can be 
heterogeneous and non-Gaussian.
However, in this section we assume the noise distribution to be known.
The main oracle results of Theorem~\ref{ToracleSmA} below do not require any 
further conditions on the noise. 
Some upper bounds on the quantities \( \ZZs_{\ms} \) entering in the oracle bounds are established 
under i.i.d. Gaussian noise, but can be easily extended to non-Gaussian heterogeneous 
noise under moment conditions.
%These results will be used in the next section as benchmark 
%to judge whether the problem of adaptive estimation can be solved 
%under lacking information about the noise.  
%
For the linear model \eqref{Yifieiiid},
define \( \thetavs \in \R^{\dimp} \) as the vector providing the best linear fit:
\begin{EQA}
	\thetavs
	& \eqdef &
	\argmin_{\thetav} \E \| \Yv - \Psi^{\T} \thetav \|^{2}
	=
	\bigl( \Psi \Psi^{\T} \bigr)^{-1} \Psi \fvs .
\label{tsYvPTts2}
\end{EQA}
As usual, a pseudo-inversion is assumed if the matrix \( \Psi \Psi^{\T} \) is degenerated.

Below we assume a family \( \bigl\{ \tilde{\thetav}_{\mm} \bigr\} \) 
of linear estimators \( \tilde{\thetav}_{\mm} = \Gam_{\mm} \Yv \)
of \( \thetavs \) to be given.
Typical examples include projection estimation on an \( \mm \)-dimensional subspace 
or regularized estimation with a regularization parameter \( \alp_{\mm} \),
penalized estimators with a quadratic penalty function, etc.
To include specific problems like subvector/functional estimation, 
we also introduce a weighting \( \dimq \times \dimp \)-matrix \( \QL \) for some fixed 
\( \dimq \geq 1 \) and 
%some linear map \( \QL \thetavs \) for a given weighting matrix \( \QL \). 
%
%\subsection{Loss and risk}
%First we have to define the loss and risk.
define quadratic loss and risk with this weighting matrix \( \QL \):
\begin{EQ}[rcl]
	\losst_{\mm}
	\eqdef 
	\| \QL (\tilde{\thetav}_{\mm} - \thetavs) \|^{2} ,
	&\qquad &
	\riskt_{\mm}
	\eqdef 
	\E \| \QL (\tilde{\thetav}_{\mm} - \thetavs) \|^{2} .
\label{lossQL}
\end{EQ}
Of course, the loss and the risk depend on the choice of \( \QL \).
We do not indicate this dependence explicitly but it is important to keep in mind 
the role of \( \QL \) in the definition of \( \losst_{\mm} \).
Typical examples of \( \QL \) are as follows.

\paragraph{Estimation of the whole vector \( \thetavs \)}
Let \( \QL \) be the identity matrix \( \QL = \Id_{\dimp} \) with \( \dimq = \dimp \).
This means that the estimation loss is measured by the usual squared Euclidean norm 
\( \| \tilde{\thetav}_{\mm} - \thetavs \|^{2} \).

\paragraph{Prediction}
Let \( \QL \) be the square root of the total Fisher information matrix 
\( \IF = \sigma^{-2} \Psi \Psi^{\T} \), that is, \( \QL^{2} = \IF \).
%Then minimization of the loss \( \losst_{\mm} \) of \( \tilde{\thetav}_{\mm} \) 
%is equivalent to maximizing 
%the log-likelihood \( - (2 \sigma^{2})^{-1} \| \Psi^{\T} (\thetav - \thetavs) \|^{2} \)
%over the corresponding \( \mm \)-dimensional subspace.
Such a type of loss is usually referred to as \emph{prediction loss} because 
it measures the fit and the prediction ability of the true model by the model with the parameter \( \thetav \).

\paragraph{Semiparametric estimation}
Let the target of estimation not be the whole vector \( \thetavs \) but some subvector 
\( \thetavs_{0} \) of dimension \( \dimq \).
The estimate \( \Pi \tilde{\thetav}_{\mm} \) is called the 
\emph{profile maximum likelihood estimate}.
The matrix \( \QL \) can be defined as the projector \( \Pi_{0} \) 
on the \( \thetavs_{0} \) subspace.
The corresponding loss is equal to the squared Euclidean norm in this subspace:
\begin{EQA}
	\losst_{\mm}
	&=&
	\| \Pi_{0} \bigl( \tilde{\thetav}_{\mm} - \thetavs \bigr) \|^{2} .
\label{lossmmsemi}
\end{EQA} 
Alternatively, one can select \( \QL^{2} \) as the efficient Fisher information matrix 
defined by the relation
\begin{EQA}
	\QL^{2}
	& \eqdef &
	\IFr
	=
	\bigl( \Pi_{0} \IF^{-1} \Pi_{0}^{\T} \bigr)^{-1} .
\label{QLsemiIFmm}
\end{EQA}

\paragraph{Linear functional estimation}
The choice of the weighting matrix \( \QL \) can be adjusted to address the problem of estimating 
some functionals of the whole parameter \( \thetavs \).
\ifbook{}{
For instance, in the regression problem \( \E Y_{i} = \fs(X_{i}) \) with 
the Fourier expansion the target function \( \fs \) can be represented as
\begin{EQA}
	\fs(x) 
	&=& 
	\sum_{j \geq 0} \thetas_{j} \psi_{j}(x) 
	=
	\sum_{j \geq 0} 
	\bigl\{ \thetas_{2j} \cos(2\pi j x ) + \thetas_{2j+1} \sin( 2 \pi j x) \bigr\} .
\label{fsxcosjx}
\end{EQA} 
The value of this function at zero coincides with the functional 
\(
	\fs(0)
	=
	\sum_{j} \thetas_{2j}
\).
The first derivative of this function leads to the functional
\( \fs'(0) = 2 \pi \sum_{j \geq 0} j \thetas_{2j+1} \).
}

\medskip
In all cases, the most important feature of the estimators \( \tilde{\thetav}_{\mm} \) is \emph{linearity}.
It greatly simplifies the study of their properties including the prominent
bias-variance decomposition of the risk of \( \tilde{\thetav}_{\mm} \).
%
%\begin{lemma}
%\label{Tlinsmoothr}
Namely, for the model \eqref{Yifieiiid} with \( \E \epsv = 0\), it holds
%\( \Var(\epsv) = \sigma^{2} \Id_{n} \).
%Then 
\begin{EQA}
	\E \tilde{\thetav}_{\mm}
	&=&
	\thetavs_{\mm}
	=
	\Gam_{\mm} \fvs ,
%	\\
%	\Var\bigl( \tilde{\thetav}_{\mm} \bigr)
%	&=&
%	\VP_{\mm}^{2}
%	=
%%	\sigma^{2} 
%	\Gam_{\mm} \Var(\epsv) \Gam_{\mm}^{\T},
	\\
	\riskt_{\mm}
	&=&
	\| \QL \bigl( \thetavs_{\mm} - \thetavs \bigr) \|^{2} 
%	+ \sigma^{2} 
	+ \tr \bigl\{ \QL \Gam_{\mm} \, \Var(\epsv) \, \Gam_{\mm}^{\T} \QL^{\T} \bigr\}
	\\
	&=&
	\| \QL (\Gam_{\mm} - \Gam) \fvs \|^{2} 
%	+ \sigma^{2} 
	+ \tr \bigl\{ \QL \Gam_{\mm} \, \Var(\epsv) \, \Gam_{\mm}^{\T} \QL^{\T} \bigr\}.
\label{Erisktmm}
\end{EQA}
%\end{lemma}
%The expression for the risk \( \riskt_{\mm} \) in \eqref{Erisktmm} can be called the 
%``bias-variance'' decomposition.
The optimal choice of the parameter \( \mm \) can be defined by risk minimization:
\begin{EQA}
	\ms
	& \eqdef &
	\argmin_{\mm \in \MM} \riskt_{\mm} .
\label{msQldef}
\end{EQA}
The \emph{model selection} problem can be described as the choice of \( \mm \) by data 
which \emph{mimics the oracle}, that is, 
we aim at constructing a selector \( \hat{\mm} \) leading to the adaptive estimate 
\( \hat{\thetav} = \tilde{\thetav}_{\hat{\mm}} \) with properties similar to the oracle
estimate \( \tilde{\thetav}_{\ms} \).

\medskip
Below we discuss the \emph{ordered case}.
The parameter \( \mm \in \MM \) is treated as complexity of the method 
\( \tilde{\thetav}_{\mm} \).
%and this parameter may grows until the maximal value \( \mmmax \).
In some cases the set \( \MM \) of possible \( \mm \) choices can be countable and/or 
continuous and even unbounded.
For simplicity of presentation, we assume 
that \( \MM \) is a finite set of positive numbers,
\( |\MM| \) stands for its cardinality.
Typical examples are given by the number of terms in the Fourier expansion,
or by the bandwidth in the kernel smoothing. 
%We also assume without loss of generality that \( \mm \) is an integer non-negative number.
In general, complexity can be naturally expressed via 
the variance of the stochastic term of the 
estimator \( \tilde{\thetav}_{\mm} \): the larger \( \mm \), the larger is the variance
\( \Var (\QL \tilde{\thetav}_{\mm}) \).
In the case of projection estimation with \( \mm \)-dimensional 
projectors \( \Gam_{\mm} \), this variance 
is linear in \( \mm \), \( \Var(\tilde{\thetav}_{\mm}) = \sigma^{2} \mm \).
In general, dependence of the variance term on \( \mm \) may be more complicated but 
the monotonicity constraint \eqref{Varmono} has to be preserved.

Further, it is implicitly assumed that 
the bias term \( \| \QL (\thetavs - \thetavs_{\mm}) \|^{2} \)
becomes small when \( \mm \) increases.
The smallest index \( \mm = \mmmin \) corresponds to the simplest (zero) model with probably a large bias,
while \( \mm \) large ensures a good approximation quality 
\( \thetavs_{\mm} \approx \thetavs \) 
and a small bias at cost of a big complexity measured by the variance term.
In the case of projection estimation, the bias term in \eqref{Erisktmm} describes 
the accuracy of approximating the response \( \fvs \) by an \( \mm \)-dimensional 
linear subspace and this approximation improves as \( \mm \) grows.
However, in general, in constrast to the case of projection estimation,
one cannot require that the bias term \( \| \QL (\thetavs - \thetavs_{\mm}) \| \)
monotonously decreases with \( \mm \).
One example is given by an estimation-at-a-point problem.

\begin{example}
Suppose that a signal \( \thetavs \) is observed with noise:
\( Y_{i} = \thetas_{j} + \eps_{j} \).
Consider the set of projection estimates 
\( \tilde{\thetav}_{\mm} \) 
on the first \( \mm \) coordinates and the target is 
\( \tars \eqdef \QL \thetav = \sum_{j} \theta_{j} \). 
If \( \thetavs \) is composed of alternating blocks of \( 1 \)'s and \( -1 \)'s with equal length, then the bias \( |\tars - \tars_{\mm}| \) for \( \tars_{\mm} = \sum_{j \leq \mm} \thetas_{j} \) 
is not monotonous in \( \mm \). 
\end{example}

\subsection{Smallest accepted (SmA) method in ordered model selection}
\label{SSmA}
First we recall our setup.
Due to the linear structure of the estimators \( \tilde{\thetav}_{\mm} = \Gam_{\mm} \Yv \) 
and of the loss function \( \QL \), one can consider 
\( \tilde{\tarv}_{\mm} = \Tam_{\mm} \Yv \) with 
\( \Tam_{\mm} = \QL \Gam_{\mm} \colon \R^{n} \to \R^{\dimq} \),
\( \mm \in \MM \),
as a family of linear estimators 
of the \( \dimq \)-dimensional target of estimation 
\( \tarvs = \QL \thetavs = \QL \Gam \fvs = \Tam \fvs \) for \( \Tam = \QL \Gam \).

Now we discuss a general approach to model selection problems based on multiple testing.
Suppose that the given family \( \bigl\{ \tilde{\tarv}_{\mm} \, , \, \mm \in \MM \bigr\} \) 
of estimators is naturally ordered by their complexity (variance).
Due to \eqref{Erisktmm}, this condition can be written as 
\begin{EQA}
	\Tam_{\mm} \Var(\epsv) \, \Tam_{\mm}^{\T}
	& \leq &
	\Tam_{\mc} \Var(\epsv) \, \Tam_{\mc}^{\T} ,
	\qquad
	\mc > \mm .
\label{Varmono}
\end{EQA}
One would like to pick up a smallest possible index \( \mm \in \MM \) which still 
provides a reasonable fit. 
The latter means that the bias component 
\begin{EQA}
	\| \bias_{\mm} \|^{2}
	&=&
	\| \tarvs_{\mm} - \tarvs \|^{2} 
	= 
	\| (\Tam_{\mm} - \Tam) \fvs \|^{2}
\label{biasmmTTs2}
\end{EQA}
in the risk decomposition \eqref{Erisktmm} % with \( \Tam = \QL \Gam \) 
is not significantly larger than the variance 
\begin{EQA}
	\tr \bigl\{ \Var\bigl( \tilde{\tarv}_{\mm} \bigr) \bigr\}
	&=&
	\tr \bigl\{ \Tam_{\mm} \Var(\epsv) \, \Tam_{\mm}^{\T} \bigr\} .
\label{trVarttarm}
\end{EQA}
If \( \md \in \MM \) is such a ``good'' choice, then our ordering assumption yields that 
a further increase of the index \( \mm \) over \( \md \) only increases 
the complexity (variance) of the method without real gain in the quality of approximation.
This latter fact can be interpreted in term of pairwise comparison:
whatever \( \mm \in \MM \) with \( \mm > \md \) we take, there is no significant bias reduction in 
using a larger model \( \mm \) instead of \( \md \). 
This leads to a multiple testing procedure: for each pair \( \mm > \md \) from \( \MM \), we consider 
%Let \( H_{\md} \) vs. \( H_{\mm} \) mean 
a hypothesis of no significant bias 
between the models \( \md \) and \( \mm \), and let \( \test_{\mm,\md} \) be 
the corresponding test.
The model \( \md \) is accepted if \( \test_{\mm,\md} = 0 \) for all \( \mm > \md \).
Finally, the selected model is the ``smallest accepted'':
\begin{EQA}
	\hat{\mm}
	& \eqdef &
	\argmin \bigl\{ \md \in \MM \colon \test_{\mm,\md} = 0, \, \forall \mm > \md \bigr\}.
\label{SmAphimmd}
\end{EQA}
Usually the test \( \test_{\mm,\md} \) can be written in the form 
\begin{EQA}
	\test_{\mm,\md} 
	&=& 
	\Ind\bigl\{ \Tmd_{\mm,\md} > \ZZ_{\mm,\md} \bigr\}
\label{phimmmdIZZ}
\end{EQA} 
for some \emph{test statistics} \( \Tmd_{\mm,\md} \) and for \emph{critical values}
\( \ZZ_{\mm,\md} \).
The information-based criteria like AIC or BIC use the likelihood ratio test statistics
\( \Tmd_{\mm,\md} = \sigma^{-2} \bigl\| \Psi^{\T} \bigl( \tilde{\thetav}_{\mm} - \tilde{\thetav}_{\md} \bigr) \bigr\|^{2} \).
A great advantage of such tests is that the test statistic \( \Tmd_{\mm,\md} \) 
is pivotal (\( \chi^{2} \) with \( \mm - \md \) degrees of freedom) 
under the correct null hypothesis, this makes it simple to compute
the corresponding critical values. 
Below we apply another choice corresponding to Lepski-type procedure and based on the 
norm of differences \( \tilde{\tarv}_{\mm} - \tilde{\tarv}_{\md} \):
\begin{EQA}
	\Tmd_{\mm,\md}
	&=&
	\| \tilde{\tarv}_{\mm} - \tilde{\tarv}_{\md} \|
	=
	\| \Tam_{\mm,\md} \Yv \| ,
	\qquad
\label{Tmdmmmd}
%\end{EQA} 
%where 
%\begin{EQA}
	\Tam_{\mm,\md}
	\eqdef 
	\Tam_{\mm} - \Tam_{\md} \, .
\label{Tammmmddeff}
\end{EQA}
The main issue for such a method is a proper choice of the critical values 
\( \ZZ_{\mm,\md} \).
One can say that the procedure is specified by a way of selecting these critical values.
Below we offer a novel way of carrying out this choice in a general situation by using a so-called
\emph{propagation condition}:
if a model \( \md \) is ``good'' it has to be accepted with a high probability.
This rule can be seen as an analog of the family-wise level condition in a multiple testing problem.
Rejecting a ``good'' model is the family-wise error of first kind, 
and this error has to be controlled.

\subsection{Oracle choice}
To specify precisely the meaning of a good model, we use below for each pair
\( \mm > \md \) from \( \MM \) the decomposition 
\begin{EQA}
	\Tmd_{\mm,\md}
	&=&
	\| \tilde{\tarv}_{\mm} - \tilde{\tarv}_{\md} \|
	=
	\| \Tam_{\mm,\md} \Yv \|
	=
	\| \Tam_{\mm,\md} (\fvs + \epsv) \|
	=
	\| \bias_{\mm,\md} + \xiv_{\mm,\md} \| ,
\label{Tmdmmmddec}
\end{EQA}
where with \( \Tam_{\mm,\md} = \Tam_{\mm} - \Tam_{\md} \)
\begin{EQA}
	\bias_{\mm,\md} 
	\eqdef  
	\Tam_{\mm,\md} \fvs ,
	&\quad &
	\xiv_{\mm,\md} 
	\eqdef 
	\Tam_{\mm,\md} \epsv .
\label{biasxivmmmd}
\end{EQA}
We also define 
 \begin{EQA} 
	\biasv_{\mm} 
	 \eqdef  
	\Tam_{\mm} \fvs, \quad &  
	\xiv_{\mm} 
	 \eqdef 
	\Tam_{\mm} \epsv . & \label{eq:biasvecdefsingle} 
\end{EQA}
It obviously holds \( \E \xiv_{\mm,\md} = 0 \).
Introduce the \( \dimq \times \dimq \)-matrix \( \VQ_{\mm,\md} \) as the variance 
of \( \tilde{\tarv}_{\mm} - \tilde{\tarv}_{\md} \):
\begin{EQA}
	\VQ_{\mm,\md}
	& \eqdef &
	\Var\bigl( \tilde{\tarv}_{\mm} - \tilde{\tarv}_{\md} \bigr)
	=
	\Var\bigl( \Tam_{\mm,\md} \Yv \bigr)
	=
	\Tam_{\mm,\md} \Var(\epsv) \, \Tam_{\mm,\md}^{\T} .
\label{VQmmmddef}
\end{EQA}
If the noise \( \epsv \) is homogeneous
with \( \Var(\epsv) = \sigma^{2} \Id_{\nsize} \),
it holds
\begin{EQA}
	\VQ_{\mm,\md}
	& = &
	\sigma^{2} \, \Tam_{\mm,\md} \, \Tam_{\mm,\md}^{\T} .
\label{VQmmmdhomo}
\end{EQA}
Further,
\begin{EQA}
	\E \, \Tmd_{\mm,\md}^{2}
	&=&
	\| \bias_{\mm,\md} \|^{2} + \E \| \xiv_{\mm,\md} \|^{2}
	=
	\| \bias_{\mm,\md} \|^{2} + \dimA_{\mm,\md} ,
\label{ETmmmd}
	\\
	\dimA_{\mm,\md} 
	&\eqdef& 
	\tr (\VQ_{\mm,\md}) = \E \| \xiv_{\mm,\md} \|^{2} .
\label{dimAmmmdd}
\end{EQA}
The bias term \( \bias_{\mm,\md} \eqdef \Tam_{\mm,\md} \fvs \) 
is significant if 
its squared norm is competitive with the variance term 
\( \dimA_{\mm,\md} = \tr (\VQ_{\mm,\md}) \).
We say that \( \md \) is a ``good'' choice if there is no significant bias 
\( \bias_{\mm,\md} \) for any \( \mm > \md \).
This condition can be quantified in the following ``bias-variance trade-off'':
\begin{EQA}
	\| \bias_{\mm,\md} \|^{2}
	& \leq &
	\alpb^{2} \, \dimA_{\mm,\md} \, ,
	\qquad
	\mm > \md
\label{biasmmmdalpzqu}
\end{EQA}
for a given parameter \( \alpb \) which controls the bias component in the risk due to 
decomposition \eqref{ETmmmd}.
Now define the \emph{oracle} \( \ms \) as the minimal \( \md \) 
with the property \eqref{biasmmmdalpzqu}:
\begin{EQA}
	\ms
	& \eqdef &
	\min \Bigl\{ \md \colon 
		\max_{\mm > \md} 
			\bigl\{ \| \bias_{\mm,\md} \|^{2} - \alpb^{2} \, \dimA_{\mm,\md} \bigr\} 
		\leq 0
	\Bigr\} .
\label{msdeflin}
\end{EQA}

\subsection{Tail function, multiplicity correction, critical values \( \ZZ_{\mm,\md} \)}
Now we explain a possible choice of critical values \( \ZZ_{\mm,\md} \) in the situation
when the noise distribution is known. 
A particular example is the case of Gaussian errors 
\( \epsv \sim \ND(0,\sigma^{2} \Id_{n}) \).
Then the distribution of the stochastic component \( \xiv_{\mm,\md} \) is known as well.
In the Gaussian case, 
it is \( \ND(0,\VQ_{\mm,\md}) \) with the covariance matrix 
\( \VQ_{\mm,\md} \).
Introduce for each pair \( \mm > \md \) from \( \MM \) a \emph{tail function} 
\( \zq_{\mm,\md}(\xxt) \) 
of the argument \( \xxt \) such that %for any \( \xx \)
\begin{EQA}
	\P\Bigl( 
%		\zqm(\VQ_{\mm,\md},\xx) 
%		\leq 
		\| \xiv_{\mm,\md} \| 
		> \zq_{\mm,\md}(\xxt)
	\Bigr)
	& = &
	\ex^{-\xxt} .
\label{Pximdmmub}
\end{EQA}
Here we assume that the distribution of \( \| \xiv_{\mm,\md} \| \) is continuous 
and the value \( \zq_{\mm,\md}(\xxt) \) is well defined.
Otherwise one has to define \( \zq_{\mm,\md}(\xxt) \) as the smallest value for which the error probability is smaller than \( \ex^{-\xxt} \).
 
For checking the propagation condition, we need a uniform in \( \mm > \md \) version of the probability bound \eqref{Pximdmmub}.
Let 
\begin{EQA}
	\MM^{+}(\md)
	& \eqdef &
	\bigl\{ \mm \in \MM \colon \mm > \md \bigr\} .
\label{MMpmddef}
\end{EQA}
Given \( \xx \),
by \( \qq_{\md} = \qq_{\md}(\xx) \) denote the corresponding multiplicity correction: 
\begin{EQA}
	\P\Bigl( 
		\bigcup_{\mm \in \MM^{+}(\md)} 
		\bigl\{ 
%			\zzm(\dimq,\xx + \qq_{\md}) \leq 
			\| \xiv_{\mm,\md} \|
			\geq \zq_{\mm,\md}(\xx + \qq_{\md}) 
		\bigr\}
	\Bigr)
	& = &
	\ex^{-\xx} .
\label{Pximdmmubu}
\end{EQA}
%We write \( \xx_{\md} \eqdef \xx + \qq_{\md}(\xx) \).
A simple way of computing the multiplicity correction \( \qq_{\md} \) is based 
on the Bonferroni bound: \( \qq_{\md} = \log(\#\MM^{+}(\md)) \).
However, it is well known that the Bonferroni bound is very conservative and leads 
to a large correction \( \qq_{\md} \), especially if the random vectors 
\( \xiv_{\mm,\md} \) are strongly correlated.
This is exactly the case under consideration.
Note that the joint distribution of the \( \xiv_{\mm,\md} \)'s is 
precisely known.
This allows to define the correction \( \qq_{\md} = \qq_{\md}(\xx) \) just by 
condition \eqref{Pximdmmubu}.
Finally we define the critical values \( \ZZ_{\mm,\md} \) by one more correction 
for the bias:
\begin{EQA}
	\ZZ_{\mm,\md} %
	& \eqdef &
	\zq_{\mm,\md}(\xx + \qq_{\md}) + \alpb \sqrt{\dimA_{\mm,\md}} 
\label{ZZmmmddef}
\end{EQA}
for \( \dimA_{\mm,\md} = \tr( \VQ_{\mm,\md}) \).
This definition still involves two numerical tuning constants \( \xx \) and \( \alpb \).
The first value \( \xx \) controls the nominal rejection probability under the null,
a usual choice \( \xx = 3 \) does a good job in most of cases. 
The value \( \alpb \) controls the amount of admissible bias in the definition 
of a good choice; cf. \eqref{biasmmmdalpzqu} and \eqref{msdeflin}.
%Any choice in the range between 0 and 1 works well, 
This value is mainly for theoretical 
study, in practice one can always take \( \alpb = 0 \).

\subsection{SmA choice and the oracle inequality}
Define the selector \( \hat{\mm} \) by the ``smallest accepted'' (SmA) rule.
Namely, with \( \ZZ_{\mm,\md} \) from \eqref{ZZmmmddef},
the acceptance rule reads as follows:
\begin{EQA}
	\bigl\{ \md \text{ is accepted} \bigr\}
	& \Leftrightarrow &
	\Bigl\{ \max_{\mm \in \MMu(\md)} 
		\bigl\{ 
			\Tmd_{\mm,\md} - \ZZ_{\mm,\md} 
		\bigr\} \leq 0 
	\Bigr\} .
\end{EQA}
The SmA rule is
\begin{EQA}
	\hat{\mm}
	& \eqdef &
	\text{``smallest accepted''} 
	\\
	&=&
	\min \Bigl\{ \md \colon
	\max_{\mm \in \MMu(\md)} 
		\bigl\{ 
			\Tmd_{\mm,\md} - \ZZ_{\mm,\md} 
		\bigr\} \leq 0 
	\Bigr\}.
\label{hammlindef}
\end{EQA} 
Our study mainly focuses on the behavior of the selector \( \hat{\mm} \).
The performance of the resulting estimator 
\( \hat{\tarv} = \tilde{\tarv}_{\hat{\mm}} \) is a kind of corollary from 
statements about the selected model \( \hat{\mm} \).
The ideal solution would be \( \hat{\mm} \equiv \ms \), then the adaptive estimator 
\( \hat{\tarv} \) coincides with the oracle estimate \( \tilde{\tarv}_{\ms} \).

The bound \eqref{Pximdmmub} automatically ensures the desired \emph{propagation property}:
any good model \( \md \) in the sense \eqref{biasmmmdalpzqu} will be accepted with probability 
at least \( 1 - \ex^{-\xx} \).
In some sense, this property is built-in by the construction of the procedure.
By definition, the oracle \( \ms \) is also a ``good'' choice, this yields 
\begin{EQA}
	\P\bigl( \ms \text{ is rejected} \bigr)
	& \leq &
	\ex^{-\xx} .
\label{Ppropmslin}
\end{EQA}
Therefore, the selector \( \hat{\mm} \) typically takes its value in \( \MMd(\ms) \), where
\begin{EQA}
	\MMd(\ms) 
	&=& 
	\bigl\{ \mm \in \MM \colon \mm < \ms \bigr\} 
\label{MMdmsdef}
\end{EQA}
is the set of all models in \( \MM \) smaller than \( \ms \).
It remains to check the performance of the method in this region.
%
%We already know that the construction ensures the propagation property:
%\( \hat{\mm} \leq \ms \) with a high probability.
%In other words, with high probability \( \hat{\mm} \in \MMd(\ms) \).
%
The next step is to specify a subset \( \MMi \) of \( \MMd(\ms) \) of highly probable 
\( \hat{\mm} \)-values. 
We will refer to this subset as the \emph{zone of insensitivity}.
The definition of \( \ms \) implies that there is a significant bias for each 
\( \mm \in \MMd(\ms) \).
If this bias is really large, then, again, the probability of selecting \( \mm \) 
can be bounded from above by a small value. 
Therefore, the zone of insensitivity is composed of \( \mm \)-values for which 
the bias is significant but not very large.

Now we present a formal description which specifies 
a subset \( \MMc \) of \( \MMd(\ms) \) 
for which the bias \( \| \bias_{\ms,\mm} \| \) is sufficiently large and hence, 
the probability of the event \( \bigl\{ \hat{\mm} \in \MMc \bigr\} \) is
negligible.

\begin{theorem}
\label{ToracleSmA}
Let \( \zq_{\mm,\md}(\cdot) \) be the tail function from \eqref{Pximdmmub} 
for each pair \( \mm > \md \in \MM \).
Given \( \xx \) and \( \alpb \),
let \( \ZZ_{\mm,\md} \) be given by \eqref{Pximdmmubu} and \eqref{ZZmmmddef}.
Then the propagation property \eqref{Ppropmslin} is fulfilled for the SmA selector 
\( \hat{\mm} \).
Moreover, for any subset \( \MMc \subseteq \MMd(\ms) \) s.t. 
\begin{EQA}
	\| \bias_{\ms,\mm} \|
	& > &
	\ZZ_{\ms,\mm} + \zq_{\ms,\mm}(\xxu) ,
	\qquad 
	\mm \in \MMc,
\label{MMxxdeflin}
\end{EQA}
for \( \xxu \eqdef \xx + \log(|\MMc|) \) with \( |\MMc| \) being the cardinality of 
\( \MMc \), it holds
\begin{EQA}
	\P\bigl( \hat{\mm} \in \MMc \bigr)
	& \leq &
	\ex^{-\xx} .
\label{PMMcSmA}
\end{EQA}
The SmA estimator
\( \hat{\tarv} = \tilde{\tarv}_{\hat{\mm}} \) satisfies the following bound:
\begin{EQA}
	\P\Bigl( \bigl\| \hat{\tarv} - \tilde{\tarv}_{\ms} \bigr\|
		> \ZZs_{\ms} 
	\Bigr)
	& \leq &
	2 \ex^{-\xx} \, ,
\label{QLtmmtmshh}
\end{EQA}
where \( \ZZs_{\ms} \) is defined 
with \( \MMi \eqdef \MMd(\ms) \setminus \MMc \) 
as 
\begin{EQA}
	\ZZs_{\ms}
	& \eqdef &
	\max_{\mm \in \MMi} \ZZ_{\ms,\mm} \, .
\label{ZZbmsdef}
\end{EQA}
This implies the probabilistic oracle bound:
with probability at least \( 1 - 2 \ex^{-\xx} \)
\begin{EQA}
	\bigl\| 
		\hat{\tarv} - \tarvs 
	\bigr\|
	& \leq &
	\bigl\| \tilde{\tarv}_{\ms} - \tarvs \bigr\| + \ZZs_{\ms} .
\label{Proboralin}
\end{EQA}
\end{theorem}

\begin{remark}
Note that the choice \( \xxu = \xx + \log(|\MMc|) \) relies on crude Bonferroni arguments and 
the definition of \( \MMc \) can be refined by choosing \( \xxu \) more
carefully.
However, this value only enters in the theoretical bound and is not used 
in the procedure, a fine tuning for this value is not required.
Obviously \( \xxu \leq \xx + \log (|\MMd(\ms)|) \).
\end{remark}

\begin{remark}
The result \eqref{Proboralin} is called the \emph{oracle bound} 
because it compares the loss of the data-driven 
selector \( \hat{\mm} \) and of the optimal choice \( \ms \).
The value \( \ZZs_{\ms} \) in \eqref{ZZbmsdef} can be viewed as a ``payment for adaptation''.
An interesting feature of the presented result is that not only the oracle quality but also the payment for adaptation depend upon the unknown response \( \fvs \) and 
the corresponding oracle choice \( \ms \).
In the worst case of a model with a flat risk profile \( \riskt_{\mm} \), the set 
\( \MMi \) can coincide with the whole range \( \MMd(\ms) \).
Even in this case the bounds \eqref{QLtmmtmshh} and \eqref{Proboralin} are meaningful.
However, the payment for adaptation \( \ZZs_{\ms} \) in this case 
can be larger than the oracle risk.
In the contrary, if the risk function \( \riskt_{\mm} \) grows rapidly as 
\( \mm \) decreases below \( \ms \), then the set \( \MMi \) is small 
and the value \( \ZZs_{\ms} \) is much smaller than the oracle risk \( \riskt_{\ms} \).
\end{remark}

\subsection{Analysis of the payment for adaptation \( \ZZs_{\ms} \)}
Here we present an upper bound on \( \ZZs_{\ms} \) for a special case of 
Gaussian independent errors \( \eps_{i} \).
%The results can be easily extended to heterogeneous sub-Gaussian errors \( \eps_{i} \). 
%
The benefit of considering the Gaussian case is that each 
vector \( \xiv_{\mc,\mm} \) is 
Gaussian as well, which simplifies the analysis of the tail function 
\( \zq_{\mc,\mm}(\cdot) \).
%
%These arguments yield the following result.
However, the results can be extended to non-Gaussian errors \( \eps_{i} \) under 
exponential moment conditions.
Below \( \mmmin \) denotes the smallest model in \( \MM \). Writing \( \VQ_{\mm} \eqdef \sigma^{2} \, \Tam_{\mm} \, \Tam_{\mm}^{\T}\), we define
\begin{EQA}
 \dimA_{\mm} & \eqdef &  \tr (\VQ_{\mm}) \\
 \supA_{\mm} & \eqdef &  \| \VQ_{\mm} \|_{\oper}.
\end{EQA}

\begin{theorem}
\label{Tpayment}
Assume the conditions of Theorem~\ref{ToracleSmA}.
Let also \( \dimA_{\mm,\md} = \tr (\VQ_{\mm,\md}) \) and 
\( \supA_{\mm,\md} = \| \VQ_{\mm,\md} \|_{\oper} \) with \( \VQ_{\mm,\md} = \Var(\xiv_{\mm,\md}) \) satisfy 
\( \dimA_{\ms,\mm} \leq \dimA_{\ms,\mmmin} \leq \dimA_{\ms} \) and 
\( \supA_{\ms,\mm} \leq \supA_{\ms,\mmmin} \leq \supA_{\ms} \)
 for all \( \mmmin \leq \mm < \ms \).
If the errors \( \eps_{i} \) are normal zero mean then
the critical values \( \ZZ_{\mm,\md} \) given by \eqref{ZZmmmddef} satisfy 
\begin{EQA}
	\ZZ_{\mm,\md}
	& \leq &
	(1 + \alpb) \sqrt{\dimA_{\mm,\md}} 
		+ \sqrt{2 \supA_{\mm,\md} \, \{ \xx + \log(|\MM|) \} } \, ,
\label{ZZmmmdub}
\end{EQA}
while the payment for 
adaptation \( \ZZs_{\ms} \) follows the bound
\begin{EQA}
	\ZZs_{\ms}
	&\leq &
	(1 + \alpb) \sqrt{\dimA_{\ms,\mmmin}} 
		+ \sqrt{2 \supA_{\ms,\mmmin} \, \{ \xx + \log(|\MMd(\ms)|) \} } 
%	\\
%	&\leq &
%	(1 + \alpb) \sqrt{\dimA_{\ms}} + \sqrt{\xx + \log(|\MMd(\ms)|) }
	\\
	&\leq &
	(1 + \alpb) \sqrt{\dimA_{\ms}} + \sqrt{2 \supA_{\ms} \{ \xx + \log(|\MM|) \} } \, .
\label{ZZbmsupxx}
\end{EQA}
\end {theorem}

Some special cases of this result for projection and linear functional estimation 
will be discussed in Sections~\ref{Sprojesta} and \ref{Slinf} below.

\ifAoS{
\subsection{Power loss function}
\label{Spolynomloss}

The probabilistic oracle bound of Theorem~\ref{ToracleSmA} provides some statement about
typical behavior of the adaptive SmA estimate 
\( \hat{\tarv} = \tilde{\tarv}_{\hat{\mm}} \).
Unfortunately, this bound does not yield a risk bound for quadratic or polynomial losses:
even if big losses occur with a small probability, the related risk can still be large.
It happens that the SmA procedure can be easily tuned to secure an oracle risk bound.

For simplicity of notation, we only consider the quadratic risk
\begin{EQA}
	\riskt(\hat{\tarv})
	& \eqdef &
	\E \| \hat{\tarv} - \tarvs \|^{2} .
\label{risk2SmA}
\end{EQA}
We aim at comparing the risk of the SmA procedure with the risk \( \riskt_{\ms} \)
of the oracle estimate \( \tilde{\tarv}_{\ms} \).
Recall the representation
%\begin{EQA}
%	\QL (\tilde{\tarv}_{\mm} - \tarvs)
%	&=&
%	\xiv_{\mm} + \bias_{\mm}
%\label{QLtttsmPL}
%\end{EQA}
%with \( \xiv_{\mm} = \Tam_{\mm} \epsv \) and \( \bias_{\mm} = (\Tam_{\mm} - \Tam) \fvs \)
%yielding the risk decomposition
\begin{EQA}
	\riskt_{\mm}
	& \eqdef &
	\E \bigl\| \tilde{\tarv}_{\mm} - \tarvs \bigr\|^{2}
	=
	\E \| \xiv_{\mm} \|^{2} + \| \bias_{\mm} \|^{2}
	=
	\dimA_{\mm} + \| \bias_{\mm} \|^{2} 
\label{risktmPL}
\end{EQA}
with \( \dimA_{\mm} = \tr \bigl( \VQ_{\mm} \bigr) \) and 
\( \VQ_{\mm} = \Var (\xiv_{\mm}) \).
%A similar decomposition holds for the quadratic risk of the difference 
%\( \tilde{\tarv}_{\mc} - \tilde{\tarv}_{\mm} \) for any pair \( \mc > \mm \):
%\begin{EQA}
%	\E \bigl\| \QL (\tilde{\tarv}_{\mc} - \tilde{\tarv}_{\mm}) \bigr\|^{2}
%	&=&
%	\E \| \xiv_{\mc,\mm} \|^{2} + \| \bias_{\mc,\mm} \|^{2}
%	=
%	\dimA_{\mc,\mm} + \| \bias_{\mc,\mm} \|^{2} .
%\label{risktmcmmPL}
%\end{EQA}
%%
%Usually the oracle choice is defined by minimizing the risk \( \riskt_{\mm} \).
For our analysis, we have to slightly modify the definition
of the oracle \eqref{msdeflin}.
Namely, to ensure an oracle risk bound, we require that not only the model \( \ms \) is ``good'' but also 
all the larger models \( \mm > \ms \) are ``good'' as well:
\begin{EQA}
	\ms
	& \eqdef &
	\min \Bigl\{ \md \colon 
		\max_{\mc , \mm \in \MMu(\md) \colon \mc > \mm} 
			\bigl\{ \| \bias_{\mc,\mm} \|^{2} - \alpb^{2} \, \dimA_{\mc,\mm} \bigr\} 
		\leq 0
	\Bigr\} .
	\qquad
\label{msdeflinPL}
\end{EQA}
Below we also suppose that the bias component \( \| \bias_{\mm} \|^{2} \) fulfills
\begin{EQA}
	\| \bias_{\mm} \|
	& \leq &
	\| \bias_{\ms} \|,
	\qquad
	\mm > \ms .
\label{biasmmmsPL}
\end{EQA}
%:
Otherwise, one can define 
\( \| \bias_{\ms} \| \eqdef \max_{\mm \in \MMu(\ms)} \| \bias_{\mm} \| \).

The choice of the critical values 
\( \ZZ_{\mc,\mm} \) for the SmA procedure 
has to be slightly changed to ensure a risk bound
for quadratic loss.
For this, we need a bit more detailed analysis of the SmA procedure in the 
propagation zone \( \mm > \ms \).
In this zone the variance dominates the bias,
%If such \( \mm \) is selected then \( \hat{\tarv} \) coincides with 
%\( \tilde{\tarv}_{\mm} \) whose loss and risk can be much larger than ones 
%for the oracle estimate \( \tilde{\tarv}_{\ms} \) because of its larger complexity.
%There is no essential bias component in this zone.
therefore, the SmA procedure can be tuned in the situation when there is no signal 
and hence no bias at all:
\begin{EQA}
	\Tmd_{\mc,\mm}
	=
	\bigl\| \tilde{\tarv}_{\mc} - \tilde{\tarv}_{\mm} \bigr\|
	&=&
	\| \xiv_{\mc,\mm} \| \, .
\label{QLtmstmPL}
\end{EQA}
The analysis is based on a simple but important observation that if 
\( \hat{\mm} = \mm > \ms \), then 
the good model \( \md = \mm_{\lbrack -1 \rbrack} \) is rejected, where \( \mm_{\lbrack -1 \rbrack}\) denotes the next smaller model with respect to \(\mm\). 
The latter means that 
at least one check based on \( \Tmd_{\mc,\mmprec} \) fails.
The same can be expressed as follows: 
the maximum of the r.v.'s 
\( \Tmd_{\mc,\mmprec} \Ind( \Tmd_{\mc,\mmprec} > \ZZ_{\mc,\mmprec}) \)
is positive.
For a formal description, introduce for each \( \mm \) and \( \xx \) a random event
\begin{EQA}
	A_{\mm}(\xx)
	& \eqdef &
	\Ind\Bigl( 
			\max_{\mc \in \MMu(\mm)} 
			\bigl\{ 
				\| \xiv_{\mc,\mm} \| - \zq_{\mc,\mm}(\xx)
			\bigr\} > 0 
		\Bigr) 
\label{Ammxxmcmm122def}
\end{EQA}
on which at least one of the test statistics 
\( \Tmd_{\mc,\mm} = \| \xiv_{\mc,\mm} \| \) exceeds the critical value 
\( \zq_{\mc,\mm}(\xx) \).
%
%The procedure is tuned in the ideal situation when there is no bias 
%and \( \Tmd_{\mc,\mmprec} = \| \xiv_{\mc,\mmprec} \| \). 
%Let \( \zq_{\mc,\mmprec}(\xx) \) be the tail function of the 
%\( \| \xiv_{\mc,\mmprec} \| \) 
%for  \( \xiv_{\mc,\mmprec} \sim \ND(0,\VQ_{\mc,\mmprec}) \).
The case of probabilistic loss focuses on the probability of this event, 
the value \( \xx \) is selected to make it small enough.
Now, under the polynomial loss function, we need a bound for the moment of 
the corresponding loss.
Namely, for each \( \mm \),
 consider the expectation of \( \dimA_{\mm}^{-1} \| \xiv_{\mm} \|^{2} \) 
on the random set \( A_{\mmprec}(\xx) \) :
\begin{EQA}
	\riskr_{\mm}(\xx)
	& \eqdef &
	\E \Bigl[ (\dimA_{\mm}^{-1} \| \xiv_{\mm} \|^{2} \vee 1)
		\Ind\Bigl( 
			\max_{\mc \in \MMu(\mmprec)} 
			\bigl\{ \| \xiv_{\mc,\mmprec} \| - \zq_{\mc,\mmprec}(\xx) \bigr\} > 0 
		\Bigr) 
	\Bigr] .
	\qquad
\label{riskrmmxxdef}
\end{EQA}
Similarly one can consider any other power loss function by replacing 
\( (\dimA_{\mm}^{-1/2} \| \xiv_{\mm} \|)^{2} \) with 
\( (\dimA_{\mm}^{-1/2} \| \xiv_{\mm} \|)^{q} \).
In particular, \( q = 0 \) yields the probability loss considered before.
%Here \( \dimA_{\mm} = \E \| \xiv_{\mm} \|^{2} \).

Now we define the value \( \xx_{\mmprec} \) in such a way that the related 
deviation risk \( \riskr_{\mm}(\xx) \) can be controlled from above.
Let \( \alp_{\mm} \) be a given decreasing sequence.
Its choice will be discussed below.
We fix for each \( \mm \) the value \( \xx_{\mmprec} \) such that
\begin{EQA}
	\riskr_{\mm}(\xx_{\mmprec})
	&=&
	\alp_{\mm} \, .
\label{riskrmmPL}
\end{EQA}
It implies 
\begin{EQ}[rcl]
	\E \Bigl[ \| \xiv_{\mm} \|^{2} 
		\Ind\Bigl( A_{\mmprec}(\xx)
%			\max_{\mc \in \MMu(\mmprec)} 
%			\bigl\{ \| \xiv_{\mc,\mmprec} \| - \zq_{\mc,\mmprec}(\xx_{\mmprec}) \bigr\} > 0 
		\Bigr) 
	\Bigr]
	& \leq &
	\alp_{\mm} \dimA_{\mm} \, ,
	\\
	\P \Bigl( A_{\mmprec}(\xx)
%			\max_{\mc \in \MMu(\mmprec)} 
%			\bigl\{ \| \xiv_{\mc,\mmprec} \| - \zq_{\mc,\mmprec}(\xx_{\mmprec}) \bigr\} > 0 
		\Bigr) 
	& \leq &
	\alp_{\mm} \, .
\label{riskrmmPL2}
\end{EQ}
%for a given sequence \( \alp_{\mm} > 0 \).
Now define the critical values \( \ZZ_{\mm,\md} \) of the SmA procedure as 
%Suppose that the SmA procedure is applied with 
\begin{EQA}
	\ZZ_{\mm,\md}
	&=&
	\zq_{\mm,\md}(\xx_{\md}) + \alpb \dimA_{\mm,\md}^{1/2} .
\label{zqmmmdPL}
\end{EQA}
The resulting procedure reads exactly as in the case of probabilistic loss:
\begin{EQA}
	\hat{\mm}
	&=&
	\min \Bigl\{ \md \colon
	\max_{\mm \in \MMu(\md)} 
		\bigl\{ 
			\Tmd_{\mm,\md} - \ZZ_{\mm,\md}
		\bigr\} \leq 0 
	\Bigr\}.
\label{SmAPL}
\end{EQA}
It is worth mentioning that the procedure is the same, 
and even the critical values \( \ZZ_{\mm,\md} \) are given by the same formula,
as in the case of probabilistic loss.
The only difference is in the propagation condition \eqref{riskrmmPL} which 
is a bit stronger than a similar condition for indicator loss.
This implies that the values \( \xx_{\md} \) and \( \ZZ_{\mm,\md} \) are a bit larger
in the case of a power loss function.
%As a corollary, the same procedure can be applied 

\begin{theorem}
\label{TSmApoly}
Let the SmA procedure \eqref{SmAPL} be applied with the critical values 
\( \ZZ_{\mm,\md} \) from \eqref{zqmmmdPL}, where the values \( \xx_{\mm} \) are
defined by \eqref{riskrmmPL} with the coefficients \( \alp_{\mm} \) satisfying 
\begin{EQA}
	\sum_{\mm \in \MMu(\ms)} \alp_{\mm} \dimA_{\mm}
	& \leq &
	\alpd_{\ms} \dimA_{\ms} 
\label{alpmmalpb}
\end{EQA}
for some \( \alpd_{\ms} \).
If the errors \( \eps_{i} \) are normal zero mean, then
\begin{EQA}
	\E \bigl\| \hat{\tarv} - \tarvs \bigr\|^{2}
	& \leq &
	2 \alpd_{\ms} \riskt_{\ms} 
	+ \bigl( \riskt_{\ms}^{1/2} + \ZZs_{\ms} \bigr)^{2} ,
\label{EQLttts2ms}
\end{EQA}
where
\begin{EQA}
	\ZZs_{\ms}
	& \eqdef &
	\max_{\mm \in \MMd(\ms)} \ZZ_{\ms,\mm} \, .
\label{ZZbmspoldef}
\end{EQA}
\end{theorem}

Similarly to the probabilistic loss function, the result can be refined by 
considering the zone of insensitivity in the region \( \mm < \ms \). 
%
%\tobedone{Define the zone of insensitivity. An oracle bound refined}

Now we briefly discuss the choice of constants \( \alp_{\mm} \) entering into \eqref{alpmmalpb}.
Suppose that the \( \dimA_{\mm} \)'s satisfy
\begin{EQA}
	\sum_{\mm \in \MMu(\ms)} (\dimA_{\ms}/\dimA_{\mm})^{a}
	& \leq &
	\CONST
\label{summmdimAmm}
\end{EQA}
for some \( a > 0 \) and a fixed constant \( \CONST \).
Then one can take 
\begin{EQA}
	\alp_{\mm}
	&=&
	(\dimA_{\mm}/\dimA_{\mmmin})^{-1-a} .
\label{alpmmdimAa}
\end{EQA}
Below we focus on a situation when the effective dimension \( \dimA_{\mm} \) 
grows exponentially with \( \mm \).
Note that this situation is typical in model selection and 
often one can reduce the general case to this one by a proper discretization.
Then \eqref{summmdimAmm} is fulfilled for any \( a > 0 \) with \( \CONST = \CONST(a) \).

The further step is an upper bound on the values \( \xx_{\mm} \), 
\( \zq_{\mm,\md}(\xx_{\mm}) \), and \( \ZZ_{\mm,\md} \),
as well as on the payment for adaptation \( \ZZs_{\ms} \).
These bounds require some exponential moment conditions on the errors \( \eps_{i} \).
To reduce the computational burden, we again focus on the case of Gaussian errors. 

\begin{proposition}
\label{Talpmxxm}
Suppose \eqref{summmdimAmm} for \( a > 0 \).
If the errors \( \eps_{i} \) are normal zero mean, then the choice 
\begin{EQA}
	\alp_{\mm} 
	&=& 
	\sqrt{3} (\dimA_{\mm} / \dimA_{\mmmin})^{-1-a} ,
	\qquad
	\xx_{\mmprec}
	=
	2 (1+a) \log(\dimA_{\mm}/\dimA_{\mmmin}) ,
\label{xxmm121amm}
\end{EQA}
ensures conditions \eqref{alpmmalpb}, \eqref{riskrmmPL}, 
and therefore, the oracle bound \eqref{EQLttts2ms} with 
\( \alpd_{\ms} = \sqrt{3} \CONST  (\dimA_{\mmmin}/\dimA_{\ms})^{1+a} \).
Furthermore,
\begin{EQA}
	\ZZs_{\ms}
%	& \preccurlyeq &
	& \leq &
	\alpb \sqrt{\dimA_{\ms}} 
	+  \sqrt{2 \supA_{\ms} \{ 2 (1+a) \log (\dimA_{\ms}/\dimA_{\mmmin}) + \log(|\MM|) \}} . 
\label{Zpaymentpower}
\end{EQA}
\end{proposition}
}{ % annals version
The procedure and the results can be extended to the case of polynomial loss,
see Section~A in the \suppSW.}

\subsection{Application to projection estimation}
\label{Sprojesta}
An important feature of the obtained oracle statements is their universality:
they equally apply to various setups and problems and provide some  
quantitative explicit error bounds even for finite samples.
Below we briefly comment on two popular cases of projection estimation and 
estimation of a linear functional.
In some sense, these are two extreme cases of relation between 
\( \dimA_{\ms} \) and \( \supA_{\ms} \).

This section discusses the case of projection estimation in the linear model 
\( \Yv = \Psi^{\T} \thetavs + \epsv \) with homogeneous errors \( \eps_{i} \):
\( \Var(\eps_{i}) = \sigma^{2} \).
All the conclusions can be easily extended to heterogeneous errors 
whose variances are contained in some fixed interval.
We also focus on probabilistic loss, the case of polynomial loss can be considered 
in the same way.

Let us assume an ordering on the features of \( \Psi_{\mm} \) and let for each \( \mm \in \mathbb{N}\) denote
\( \Psi_{\mm} \) as the submatrix \( \Psi \) corresponding to first \(\mm\) features, i. e. the projector onto the first \(\mm\) features. We use \(\mm\) to denote the model and the number of features.
The related estimator \( \tilde{\thetav}_{\mm} \) is the standard LSE with 
\( \Gam_{\mm} = \bigl( \Psi_{\mm} \Psi_{\mm}^{\T} \bigr)^{-1} \Psi_{\mm} \)
and the prediction problem with \( \QL = \Psi^{\T} \) yields 
\( \Tam_{\mm} \Yv = \Pi_{\mm} \Yv \) where 
\( \Pi_{\mm} = \Psi_{\mm}^{\T} \bigl( \Psi_{\mm} \Psi_{\mm}^{\T} \bigr)^{-1} \Psi_{\mm} \) 
is the projector onto the corresponding feature subspace.
For homogeneous errors \( \eps_{i} \) with \( \Var(\eps_{i}) = \sigma^{2} \), 
the variance \( \VQ_{\mm} = \Var\bigl( \Pi_{m} \Yv \bigr) \) satisfies 
\begin{EQA}
	\dimA_{\mm}
	&=&
	\tr \bigl\{ \Var\bigl( \Pi_{m} \Yv \bigr) \bigr\}
	=
	\sigma^{2} \tr \bigl( \Pi_{\mm} \bigr)
	=
	\sigma^{2} \mm .
\label{dimAmmPr}
\end{EQA}
Moreover, for each pair \( \mm > \md \), 
it holds 
\begin{EQA}
	\Psi^{\T} \bigl( \tilde{\thetav}_{\mm} - \tilde{\thetav}_{\md} \bigr)
	&=&
	\bigl( \Pi_{\mm} - \Pi_{\md} \bigr) \Yv
	=
	\Pi_{\mm,\md} \Yv ,
\label{PsiTmmmdPr}
\end{EQA}
where \( \Pi_{\mm,\md} \) projects on the 
subspace of features entering in \( \mm \) but not in \( \md \).
%Obviously
%\begin{EQA}
%	\dimA_{\mm,\md}
%	&=&
%	\sigma^{2} (\mm - \md),
%	\qquad
%	\supA_{\mm,\md}
%	=
%	\sigma^{2} .
%\label{dAmmmdsAmd}
%\end{EQA}
%
%and \( \supA_{\ms} \) are bounded by a fixed constant (equal to one for homogeneous errors).

\begin{corollary}
Consider the problem of projection estimation with homogeneous Gaussian errors \( \eps_{i} \)
and probabilistic loss.
Then \( \dimA_{\mm,\md} = \sigma^{2} (\mm - \md) \), 
\( \supA_{\mm,\md} = \sigma^{2}  \), and 
\begin{EQA}
	\ZZ_{\mm,\md}
	&\leq &
	\sigma (1 + \alpb) \sqrt{\mm - \md} + \sigma \sqrt{2 \xx + 2 \log(|\MM|) } ,
	\\
	\ZZs_{\ms}
	&\leq &
	\sigma (1 + \alpb) \sqrt{\ms} + \sigma \sqrt{2 \xx + 2 \log(|\MM|) } .
\label{ZZmmmdubpro}
\end{EQA}
\end{corollary}

The first term in the expression for \( \ZZs_{\ms} \) 
is of order \( \sqrt{\ms} \) and it is a leading one provided that the effective dimension 
\( \ms \) is essentially larger than \( \log(|\MM|) \).
Usually the cardinality of the set \( \MM \) is only logarithmic in the sample size \( n \);
cf. \cite{lepski1991, lepski1997}. 
Then \( \log(|\MM|) \approx \log \log n \) and 
\( \ZZs_{\ms} \approx \sigma \sqrt{\ms} \) for \( \ms \gg \log \log n \).
%Note that \( \dimA_{\ms} \leq \riskt_{\ms} = \dimA_{\ms} + \| \bias_{\ms} \|^{2} \).
%
For the oracle risk \( \riskt_{\ms} \), it holds 
\( \riskt_{\ms} = \dimA_{\ms} + \| \bias_{\ms} \|^{2} \geq \sigma^{2} \ms \).
Therefore, the payment for adaptation \( \ZZs_{\ms} \) is of the same order as the square root 
of the oracle risk, and the result of Proposition~\ref{Tpayment} has a surprising corollary:
rate adaptive estimation is possible if the oracle dimension \( {\ms} \) 
is significantly larger than \( \log \log n \).

\begin{remark}
The payment for adaptation can be drastically reduced in the situations
with a narrow zone of insensitivity.
If the bias grows rapidly when \( \mm \) decreases from \( \ms \) to \( \mmmin \), 
more precisely, 
if \( \| \bias_{\ms,\mm} \|^{2} 
\geq \CONST \sigma^{2} \bigl( \ms - \mm + 2\xx + 2 \log(|\MM|) \bigr) \) 
for some fixed constant \( \CONST \) and all \( \mm \leq \md \) with \( \md < \ms \), then  
\begin{EQA}
	\ZZs_{\ms}
	&\leq &
	\sigma (1 + \alpb) \sqrt{\ms - \md} + \sigma \sqrt{2 \xx + 2 \log(|\MM|) } .
\label{Zbmsmsmdin}
\end{EQA}
So, if \( (\ms - \md)/\ms \) is small, the payment for adaptation is smaller in order 
than the oracle risk, and the procedure is sharp adaptive.
In particular, one can easily see that the self-similarity condition of \cite{GiNi2010} 
ensures a rapid growth of the bias when the index \( \mm \) becomes smaller than \( \ms \).
This in turn yields a narrow zone of insensitivity and hence, a sharp adaptive estimation.  
\end{remark}

\begin{remark}
It is worth mentioning the relation of the proposed procedure to the popular 
Akaike (AIC) criterion.
AIC defines \( \hat{\mm} \) by minimizing
\begin{EQA}
	\hat{\mm}
	&=&
	\argmin_{\mm} \bigl\{ \| \Yv - \Pi_{\mm} \Yv \|^{2} + 2 \sigma^{2} \mm \bigr\}.
\label{hammAIC}
\end{EQA}
One can easily see that this rule is equivalent to the SmA rule \eqref{hammlindef}
with \( \ZZ_{\mm,\md}^{2} = 2 \sigma^{2} (\mm - \md) \).
However, this choice does not guarantee the propagation condition \eqref{Ppropmslin}.
\end{remark}

\subsection{Linear functional estimation}
\label{Slinf}
In this section, we discuss the problem of linear functional estimation.
As previously, we assume a family of estimators 
\( \tilde{\tar}_{\mm} = \Tam_{\mm} \Yv \), \( \mm \in \MM \),
to be given,
where the rank of each \( \Tam_{\mm} \) is equal to 1.
The ordering condition means that these estimators are ordered by their variance:
\begin{EQA}
	\vp_{\mm}^{2}
	& \eqdef &
	\Var\bigl( \Tam_{\mm} \Yv \bigr) 
	=
	\Tam_{\mm} \, \Var(\epsv) \, \Tam_{\mm}^{\T} \, 
\label{vpmm2linf}
\end{EQA}
grows with \( \mm \).
Further, each stochastic component \( \xi_{\mm,\md} = \Tam_{\mm,\md} \epsv \)
is one-dimensional, and it holds 
\begin{EQA}
	\supA_{\mm,\md}
	&=& 
	\dimA_{\mm,\md} 
	=
	\vp_{\mm,\md}^{2} 
	= 
	\Tam_{\mm,\md} \, \Var(\epsv) \, \Tam_{\mm,\md}^{\T} .
\label{supAdimAvpmmmd}
\end{EQA}
Note that in the case of Gaussian errors,
\( \xi_{\mm,\md} \) is also Gaussian: \( \xi_{\mm,\md} \sim \ND(0,\vp_{\mm,\md}^{2}) \).
The tail function \( \zq_{\mm,\md}(\xx) \) of \( \xi_{\mm,\md} \) 
can be upper-bounded by \( \vp_{\mm,\md} \sqrt{2 \xx} \).
In the case of probabilistic loss,
a Bonferroni correction and a bias adjustment lead to the upper bound for 
the critical values \( \ZZ_{\mm,\md} \):
\begin{EQA}
	\ZZ_{\mm,\md}
	& \leq &
	\vp_{\mm,\md} \Bigl( \alpb + \sqrt{2 \xx + 2 \log(|\MM|)} \Bigr) ,
\label{ZZmmmdlinf}
\end{EQA}
where \( |\MM| \) is the number of elements in \( \MM \).
This implies 
\begin{EQA}
	\ZZs_{\ms}
	& \leq &
	\vp_{\ms} \Bigl( \alpb + \sqrt{2 \xx + 2 \log(|\MM|)} \Bigr) .
\label{Zbmslinf}
\end{EQA}

\begin{theorem}
Let the errors \( \eps_{i} \) be Gaussian zero mean. 
Consider a problem of linear functional estimation of \( \tars = \Tam \fvs \) 
by a given family \( \tilde{\tar}_{\mm} = \Tam_{\mm} \Yv \)
with \( \rank(\Tam_{\mm}) = \rank(\Tam) = 1 \), \( \mm \in \MM \).
Then the critical values \( \ZZ_{\mm,\md} \) from 
\eqref{ZZmmmddef} fulfill \eqref{ZZmmmdlinf} and the oracle inequality 
\eqref{Proboralin} holds with the payment for adaptation \( \ZZs_{\ms} \) obeying \eqref{Zbmslinf}.
\end{theorem}

\begin{remark}
One can conclude that for the problem of functional estimation with
 probabilistic loss, 
the squared payment for adaptation \( \ZZs_{\ms}^{2} \) is by factor \( \log(|\MM|) \) 
larger than the oracle variance \( \vp_{\ms}^{2} \).
If \( |\MM| \) itself is logarithmic in the sample size \( n \),
we end up with the extra (log log n)-factor in the accuracy of adaptive estimation.
\end{remark}

\ifAoS{ %full version
In the case of \emph{polynomial loss},
similar arguments yield due to \eqref{zqmmmdPL} and \eqref{xxmm121amm}
\begin{EQA}
	\ZZ_{\mm,\md}
	& \leq &
	\vp_{\mm,\md} \bigl( \alpb + \sqrt{2 \xx_{\md} + 2 \log(|\MM|)} \bigr)
	\\
	& \leq &
	\vp_{\mm,\md} \bigl( 
		\alpb 
		+ \sqrt{2 (1+a) \log(\dimA_{\md}/\dimA_{\mmmin}) + 2 \log(|\MM|)} 
	\bigr)
\label{ZZmmmdlinfpol}
\end{EQA}
\cite{spokoinyvial} showed that the bound 
\( \ZZ_{\mm,\md}^{2} \geq \CONST \vp_{\mm,\md}^{2} (\md - \mmmin) \)
is necessary to ensure a propagation condition for geometrically growing variance 
\( \dimA_{\mm} = \vp_{\mm}^{2} \).
The bound \eqref{ZZmmmdlinf} yields 
\begin{EQA}
	\ZZs_{\ms}
	& \leq &
	\vp_{\ms} \Bigl(  
		\alpb 
		+ \sqrt{2 (1+a) \log(\vp_{\ms}^{2}/\vp_{\mmmin}^{2}) + 2 \log(|\MM|)} 
	\Bigr) .
\label{ZZbmslinf}
\end{EQA}

\begin{theorem}
Suppose that the errors \( \eps_{i} \) are Gaussian zero mean.
Let the family of functional estimators \( \Tam_{\mm} \Yv \) be such that 
the variances \( \dimA_{\mm} = \vp_{\mm}^{2} \) from \eqref{vpmm2linf}
fulfill the condition \eqref{summmdimAmm} with \( a > 0 \).
Then the critical values \( \ZZ_{\mm,\md} \) from \eqref{zqmmmdPL} for the SmA procedure 
fulfill \eqref{ZZmmmdlinf}. 
For the resulting selector \( \hat{\mm} \), the oracle inequality \eqref{EQLttts2ms}
holds and the payment for adaptation \( \ZZs_{\ms} \) follows \eqref{ZZbmslinf}.
\end{theorem}

\begin{remark}
It appears that polynomial loss yields a larger price for adaptation:
\( \ZZs_{\ms}^{2} \asymp \vp_{\ms}^{2} \log (\vp_{\ms}^{2}/\vp_{\mmmin}^{2}) \).
%Then the leading term in \( \ZZs_{\ms} \) is of order
%\( \sqrt{\dimA_{\ms} \log \dimA_{\ms}} \) and thus, the price for adaptation is larger 
%than the oracle risk \( \riskt_{\ms} \) by a log factor \( \log(\dimA_{\ms}) \).
This conclusion is consistent with the results by \cite{lepski1992} and \cite{cailow2003,cailow2005} which show that the log-price for adaptation cannot be avoided
if a polynomial loss function is considered.
Our result seems to be even more informative because it delivers a non-asymptotic 
error bound which adapts to the underlying unknown model.  
\end{remark}
}{ %annals version
}

\section{Bootstrap tuning}
\label{Sbootcalibr}
This section explains how the proposed SmA procedure can be applied if no information about 
the noise \( \epsv = \Yv - \E \Yv \) is available.

\subsection{Presmoothing and wild boostrap}
Let the observed data \( \Yv \) follow the model \( \Yv = \fvs + \epsv \) with \(\epsv \sim \ND(0,\Sigma)\),
where \( \Sigma = \diag\bigl( \sigma_{1}^{2},\ldots, \sigma_{n}^{2} \bigr) \)
is an unknown diagonal covariance matrix. 
We assume that the response vector \( \fvs \) can be well approximated by a linear 
expansion for a given basis \( \Psi \) in the form \( \fvs \approx \Psi^{\T} \thetavs \).
The vector \( \thetavs \) can be naturally treated as target of estimation.
Assume we are given the ordered family of the estimators 
\( (\tilde{\thetav}_{\mm}) \) of \( \thetavs \):
\begin{EQA}
\label{eq:estimatorsmodelbs}
  	\tilde{\thetav}_{\mm} 
	& = & 
	\Gam_{\mm} \Yv = (\Psi_{\mm} \Psi_{\mm}^{\T})^{-1} \Psi_{\mm} \Yv,
	\quad
	\mm \in \MM.
\end{EQA}
For each pair \( \mm > \md \) from \( \MM \), we consider the test statistic 
\( \Tmd_{\mm,\md} \) and its decomposition from \eqref{Tmdmmmddec}:
with \( \Tam_{\mm,\md} = \QL (\Gam_{\mm} - \Gam_{\md}) \)
\begin{EQA}
	\Tmd_{\mm,\md}
	&=&
	\| \Tam_{\mm,\md} \Yv \|
	=
	\| \Tam_{\mm,\md} (\fvs + \epsv) \|
	=
	\| \bias_{\mm,\md} + \xiv_{\mm,\md} \| ,
\label{Tmdmmmddeca}
\end{EQA}
Calibration of the SmA model selection procedure requires to know the joint distribution 
of all corresponding stochastic terms \( \| \xiv_{\mm,\md} \| \) for \( \mm > \md \) 
which is uniquely determined by the noise covariance matrix \( \Sigma \).
In the case when this matrix is unknown, we are going to use a bootstrapping procedure to approximate this distribution.

The proposed procedure relates to the concept of the \emph{wild} bootstrap, \cite{wu1986}, \cite{beran1986}. 
In the framework of a regression problem, it 
%with normal errors 
suggests to model the unknown heteroscedastic noise using randomly weighted residuals from pilot estimation. 
We apply normal weights. For other possible bootstrap weights see for example \cite{mammen1993}.

Suppose we are given a pilot estimator (presmoothing) \( \tilde{\fv} \) of the response vector 
\( \fvs \in \R^{n} \).
Define the residuals:
\begin{EQA}
 	\Yvr 
	& \eqdef & \
	\Yv - \tilde{\fv}.
\end{EQA}
This pilot is supposed to undersmooth, that is, the bias is negligible and the variance 
of \( \Yvr \) is close to \( \Sigma \). 
This pre-smoothing requires some minimal smoothness of the regression function, 
and this condition seems to be unavoidable if no information about the noise is given:
otherwise one cannot distinguish between signal and noise. 
Below we suppose that \( \tilde{\fv} \) is a linear predictor, \( \tilde{\fv} = \Pi \Yv \),
where \( \Pi \) is a sub-projector in the space \( \R^{n} \). 
For example, one can take 
\( \Pi =  \Psi_{\mres}^{\T} \bigl( \Psi_{\mres} \Psi_{\mres}^{\T} \bigr)^{-1} \Psi_{\mres} \)
where \( \mres \) is a large model, e.g. the largest model \( \mmmax \) in our collection.

The wild bootstrap proposes to resample from the heteroscedastic Gaussian noise 
\( \Pb = \ND(0,\Sigmar) \) with 
\begin{EQA}
	\Sigmar 
	&=&
	\diag(\Yvr \cdot \Yvr),
\label{Sigmardef}
\end{EQA}
where \(\Yvr \cdot \Yvr\) denotes the coordinate-wise product of the vector \(\Yvr\) 
with itself and \(\diag(\Yvr \cdot \Yvr) \) denotes the diagonal matrix with entries 
from \(\Yvr \cdot \Yvr\).
These entries depend on \( \Yv \) and thus are random. 
Therefore, the bootstrap distribution \( \Pb \) is a random measure on \( \R^{n} \) 
and the aim of our study is to show 
that this random measure mimics well the underlying data distribution for typical 
realizations of \( \Yv \). 
Clearly \( \diag(\Yvr \cdot \Yvr) \) is a very bad estimator of the covariance matrix 
\( \Sigma \). 
However, below we show that under realistic conditions on the pilot \( \tilde{\fv} \) and 
on the model, it does a good job and allows to obtain essentially the same results
as in the case of known \( \Sigma \).

Let \( \Wb \) denote the \( n \)-vector of bootstrap weights 
\( \Wb \sim \ND(0, \Id_{n}) \).
Clearly the product \( \epsvb = \diag(\Yvr) \Wb \) is conditionally on \( \Yv \) normal,
\begin{EQA}
	\epsvb = \diag(\Yvr) \Wb \cond_{\Yv} 
	& \sim & 
	\ND(0,\Sigmar) .
\label{diagYvrepsvbcYv}
\end{EQA}
Bootstrap analog of \( \xiv_{\mm,\md} = \Tam_{\mm,\md} \epsv \) reads 
\( \xivb_{\mm,\md} = \Tam_{\mm,\md} \epsvb = \Tam_{\mm,\md} \diag(\Yvr) \Wb \) and
\begin{EQA} 
\label{xivbmmmddef}
	\| \xivb_{\mm,\md} \|  
	& \eqdef & 
	\| \Tam_{\mm,\md} \diag(\Yvr) \Wb \| \, .
\end{EQA}
The idea is to calibrate the SmA procedure under the bootstrap measure \( \Pb \) using 
\( \| \xivb_{\mm,\md} \| \) in place of \( \| \xiv_{\mm,\md} \| \).
The bootstrap quantiles \( \zqb_{\mm,\md}(\xxt) \) are given by analog of \eqref{Pximdmmub}:
\begin{EQA}
	\Pb\Bigl( 
%		\zqm(\VQ_{\mm,\md},\xx) 
%		\leq 
		\| \xivb_{\mm,\md} \| 
		> \zqb_{\mm,\md}(\xxt)
	\Bigr)
	& = &
	\ex^{-\xxt} .
\label{Pximdmmubst}
\end{EQA}
%For a given \( \xx \),
The multiplicity correction \( \qqbt_{\md} = \qqbt_{\md}(\xx) \) is specified by the condition 
\begin{EQA}
	\Pb\Biggl( 
		\bigcup_{\mm \in \MM^{+}(\md)} 
		\bigl\{ 
%			\zzm(\dimq,\xx + \qq_{\md}) \leq 
			\| \xivb_{\mm,\md} \|
			\geq \zqb_{\mm,\md}(\xx + \qqbt_{\md}) 
		\bigr\}
	\Biggr)
	& = &
	\ex^{-\xx} .
\label{Pximdmmububst}
\end{EQA}
Finally, the bootstrap critical values are fixed by the analog of \eqref{ZZmmmddef}:
\begin{EQA}
	\ZZbt_{\mm,\md} %
	& \eqdef &
	\zqb_{\mm,\md}(\xx + \qqbt_{\md}) + \alpb \sqrt{\dimAbt_{\mm,\md}} 
\label{ZZmmmddefbst}
\end{EQA}
for \( \dimAbt_{\mm,\md} = \Eb \| \xivb_{\mm,\md} \|^{2} \) given by
\begin{EQA}
 	\dimAbt_{\md,\mm} 
	& \eqdef & 
	\tr\bigl( \Tam_{\md,\mm}^{\T} \, \diag(\Yvr \cdot \Yvr) \, \Tam_{\md,\mm} \bigr).
\end{EQA}
Recall that all these quantities are data-driven and depend upon the original data.
Now we apply the SmA procedure with the critical values \( \ZZbt_{\mm,\md} \) defined in such a way.
Our main result claims that this choice still ensures the propagation condition \eqref{Pximdmmubu}
and therefore, all the obtained results including the oracle bounds, apply for this choice 
as well; see Theorem~\ref{TGaussbootB2}.
Moreover, we evaluate the distance between the unknown underlying 
data distribution \( \P \) and the bootstrap distribution \( \Pb \).
The latter is random, however, we show that with high probability, 
it is close to its deterministic counterpart \( \P \).
To make the results transparent and concise we assume a heterogeneous Gaussian noise
\( \epsv \).
All the statements can be extended to a non-Gaussian noise under some exponential moment conditions
at the cost of many technical details. 

Let \( \PPsi \) denote the joint distribution of all stochastic vectors \( \xiv_{\mm,\md} \)
entering in the decomposition of the test statistics \( \Tmd_{\mm,\md} \)
for \( \mm > \md \).
Let also \( \PPsib \) be the similar distribution of the 
bootstrapized stochastic vectors \( \xivb_{\mm,\md} \) entering in the test statistics 
\( \Tmdb_{\mm,\md} \).
The next result allows to upper bound the total variation distance between 
\( \PPsi \) and \( \PPsib \) in terms of the following quantities:

\begin{description}
	\item [Design Regularity] is measured by the value \( \dPsi \) 
\begin{EQA}
	\dPsi
	\eqdef
	\max_ {i=1,\ldots,n} \| \Varxi^{-1/2} \Psi_{i} \| \sigma_{i} \, ,
	\quad 
	& \text{where} &
	\quad
	\Varxi
	\eqdef
	\sum_{i=1}^{n} \Psi_{i} \Psi_{i}^{\T} \sigma_{i}^{2} \, ;
\label{dPsim1nsiSmA}
\end{EQA}
Obviously
\begin{EQA}
	\sum_{i=1}^{n} \| \Varxi^{-1/2} \Psi_{i} \|^{2} \sigma_{i}^{2}
	&=&
	\tr\Bigl( \sum_{i=1}^{n} \Varxi^{-2} \Psi_{i} \Psi_{i}^{\T} \sigma_{i}^{2} \Bigr)
	=
	\tr \Id_{\dimp}
	=
	\dimp ,
\label{sumSm1Psii}
\end{EQA}
and therefore in typical situations the value \( \dPsi \) is of order \( \sqrt{\dimp/n} \).

	\item [Presmoothing bias] for a projector \( \Pi \) 
	is described by the vector 
\begin{EQA}
	\Bias 
	&=& 
	\Sigma^{-1/2} (\fvs - \Pi \fvs) .
\label{BiasSm12fPi}
\end{EQA} 
We will use the sup-norm \( \| \Bias \|_{\infty} = \max_{i} |\bias_{i}| \)
and the squared \( \ell_{2} \)-norm \( \| \Bias \|^{2} = \sum_{i} \bias_{i}^{2} \)
to measure the bias after presmoothing. 

	\item[Stochastic noise after presmoothing] is described via the 
	covariance matrix \( \Var(\epsvr) \) of the smoothed noise 
	\( \epsvr = \Sigma^{-1/2} (\epsv - \Pi \epsv) \).
	Namely, this matrix is assumed to be sufficiently close to the unit matrix \( \Id_{n} \),
	in particular, its diagonal elements should be close to one.
	This is measured by the operator norm of \( \Var(\epsvr) - \Id_{n} \) and 
	by deviations of the individual variances \( \E \epsr_{i}^{2} \) from one:
\begin{EQ}[rcl]
	\supeps
	& \eqdef &
	\| \Var(\epsvr) - \Id_{n} \|_{\oper},
	\\
	\supepsi
	& \eqdef & 
	\max_{i} |\E \epsr_{i}^{2} - 1|.
\label{supepsdef}
\end{EQ}

In particular, in the case of homogeneous errors \( \Sigma = \sigma^{2} \Id_{n} \)
and the smoothing operator \( \Pi \) as a \( \dimp \)-dimensional projector, it holds
\begin{EQA}[rcl]
	\Var(\epsvr)
	& = &
	(\Id_{n} - \Pi)^{2}
	=
	\Id_{n} - \Pi
	\leq 
	\Id_{n} \, ,
\label{Covepsvrle}
%\end{EQA}
%and 
%\begin{EQA}[rcccl]
	\\
	\supeps
	&=&
	\| \Var(\epsvr) - \Id_{n} \|_{\oper}
	=
	\| \Pi \|_{\oper}
	=
	1 ,
\label{VarepsvrPiop1}
	\\
	\supepsi
	&=&
	\max_{i} |\E \epsr_{i}^{2} - 1|
	= 
	\max_{i} |\Pi_{ii}| .
\label{supepsPiii}
\end{EQA}
One can check that \( \Pi_{ii} \asymp \sqrt{\dimp/n} \) for typical smoothing 
operators like local average or kernel smoothing.
Similar bounds with an additional constant can be established for general regular noise 
\( \epsv \) and a general smoothing operator \( \Pi \).

	\item [Regularity of the smoothing operator \( \Pi \)] is required 
in Theorem~\ref{TGaussbootB2}.
This condition will be expressed via the norm of the rows \( \SPiS_{i}^{\T} \) of the matrix 
\( \SPiS \eqdef \Sigma^{-1/2} \Pi \Sigma^{1/2} \)
fulfill 
\begin{EQA}
	\| \SPiS_{i}^{\T} \|
	& \leq &
	\dPsip ,
	\qquad
	i=1,\ldots, n.
\label{SPiSidPidef}
\end{EQA}
This condition is in fact very close to the design regularity condition \eqref{dPsim1nsiSmA}.
To see this, consider the case of a homogeneous noise with \( \Sigma = \sigma^{2} \Id_{n} \)
and \( \Pi = \Psi^{\T} \bigl( \Psi \Psi^{\T} \bigr)^{-1} \Psi \).
Then \( \SPiS = \Pi \) and \eqref{dPsim1nsiSmA} implies
\begin{EQA}
	\| \SPiS_{i}^{\T} \|
	&=&
	\| \Psi^{\T} \bigl( \Psi \Psi^{\T} \bigr)^{-1} \Psi_{i} \| 
	= 
	\| \bigl( \Psi \Psi^{\T} \bigr)^{-1/2} \Psi_{i} \|
	\leq 
	\dPsi \, .
\label{SPiSidPsis}
\end{EQA}
In general one can expect that \eqref{SPiSidPidef} is fulfilled with some other constant
which however, is of the same magnitude as \( \dPsi \).
For simplicity, we use the same symbol.
\end{description}

\subsection{Bootstrap validation. Range of applicability}
This section states the main results justifying the proposed bootstrap procedure.
They claim that the joint distribution \( \PPsib \) of the bootstrap stochastic components 
\( \xivb_{\mm,\md} \) for \( \mm > \md \) nicely reproduces 
the underlying distribution \( \PPsi \) of the \( \xiv_{\mm,\md} \)'s,
and hence, all the probabilistic results obtained in Section~\ref{SKnownnoice} for  
known noise continue to apply after bootstrap parameter tuning.
In the next result, we give a bound on the total variation distance 
\( \| \PPsi - \PPsib \|_{\TV} \)
between \( \PPsi \) and \( \PPsib \).

\begin{theorem}
\label{TGaussbootB}
Let \( \Yv = \fvs + \epsv \) be a Gaussian vector in \( \R^{n} \) with independent components, 
\( \Yv \sim \ND(\fvs,\Sigma) \) for 
\( \Sigma = \diag(\sigma_{1}^{2}, \ldots, \sigma_{n}^{2}) \).
Let also \( \Psi \) be a \( \dimp \times n \) feature matrix  
such that the 
\( \dimp \times \dimp \)-matrix \( \Varxi = \Psi \, \Sigma \, \Psi^{\T} \) is non-degenerated.
For a given presmoothing operator \( \Pi \colon \R^{n} \to R^{n} \),
assume that \( \supeps \) from \eqref{supepsdef} satisfies \( \supeps \leq 1 \).
Let \( \PPsi = \LL\bigl( \xiv_{\mm,\md} , \mm,\md \in \MM \bigr) \) and 
let \( \PPsib \) be the joint conditional distribution of the bootstrap stochastic terms 
\( \xivb_{\mm,\md} \) for \( \mm,\md \in \MM \) given the data \( \Yv \).
%\( \PPsib = \LL\bigl( \xivb_{\mm,\md} , \mm,\md \in \MM \bigr) \).
%
Then it holds on a random set \( \Omega_{2}(\xx) \) with 
\( \P\bigl( \Omega_{2}(\xx) \bigr) \geq 1 - 3 \ex^{-\xx} \):
\begin{EQA}
	\| \PPsi - \PPsib \|_{\TV}
	& \leq &
	\frac{1}{2} \errSi_{2}(\xx) ,
\label{PPsiPPbTV}
	\\
	\errSi_{2}(\xx)
	& \eqdef &
	2 \sqrt{\dPsi^{2} \, \dimp \, \xxn}
	+ \sqrt{\supepsi^{2} \, \dimp} 
	+ \sqrt{\| \Bias \|_{\infty}^{4} \, \dimp}
%	+ \sqrt{\dPsi^{4} \| \Bias \|^{4} \dimp}
	+ 4 \, \dPsi^{2} \, \| \Bias \| \, \bigl( 1 + \sqrt{\xx} \bigr) .
\label{errSiuBcu}
\end{EQA}
where \( \xxn = \xx + \log(n) \), the bias \( \Bias \) is given by \eqref{BiasSm12fPi} and 
 \( \supeps \), \( \supepsi \) by \eqref{supepsdef}. 
%assume that \( \supeps \leq 1 \), and 
%\( \supepsi^{2} \, \dimp \) and \( \dimp \, \| \Bias \|_{\infty}^{4} \) are small values. 
\end{theorem}

The result \eqref{PPsiPPbTV} gives us a way to control differences \( \PPsi(A) - \PPsib(A) \)
for fixed sets \( A \).
To justify the propagation property for the bootstrap-based set of
critical values \( \zqb_{\mm,\md}(\xx + \qqbt_{\md}) \),
%the values \( \zqb_{\mm,\md}(\xx) \) and \( \qqbt_{\md} \) for all \( \mm > \md \) be defined 
given according to \eqref{xivbmmmddef}, \eqref{Pximdmmubst}, and \eqref{Pximdmmububst}
with \( \Yvr = \Yv - \Pi \Yv \), we also need to take into account the $\Yv$-dependence of \( \zqb_{\mm,\md}(\xx + \qqbt_{\md}) \). 
This is done by the following theorem.

\begin{theorem}
\label{TGaussbootB2}
Assume the conditions of Theorem~\ref{TGaussbootB},
and let the rows \( \SPiS_{i}^{\T} \) of the matrix \( \SPiS \eqdef \Sigma^{-1/2} \Pi \Sigma^{1/2} \)
satisfy \eqref{SPiSidPidef}.
Then for each \( \md \in \MM \) 
%on a set 
%\( \Omega_{0}(\xx) \) with \( \P\bigl( \Omega_{0}(\xx) \bigr) \geq 1 - 6 \ex^{-\xx} \)
\begin{EQA}
\label{multcorrBmix}
	\P\biggl( 
		\max_{\mm > \md} 
			\Bigl\{ \| \xiv_{\mm,\md} \| - \zqb_{\mm,\md}(\xx + \qqbt_{\md}) \Bigr\} 
		\geq 
		0
	\biggr)
	& \leq &
	6 \ex^{-\xx} + \sqrt{\dimp} \, \errSi_{0}(\xx) ,
	\qquad
\label{errS0def}
\end{EQA}
where with \( \xxn = \xx + \log(n) \) and \( \xxp = \xx + \log(2 \dimp) \)
\begin{EQA}
	\errSi_{0}(\xx)
	& \eqdef &
	\| \Bias \|_{\infty}^{2} + \dPsi^{2} \| \Bias \| \sqrt{2 \xx} 
	+ 2 \dPsip \xxn + \dPsip^{2} \xxn
	+ 2 \dPsi \sqrt{ \xxp} + 2 \dPsi^{2} \xxp .
\label{errS0x222}
\end{EQA}
\end{theorem}

The SmA procedure also involves the values \( \dimA_{\mm,\md} \), which are unknown 
and depend on the noise \( \epsv \).
The next result shows the bootstrap counterparts \( \dimAbt_{\mm,\md} \) can be well used 
in place of \( \dimA_{\mm,\md} \).

\begin{theorem}
\label{TdimAbtdimA}
Assume the conditions of Theorem~\ref{TGaussbootB}.
Then it holds on a set 
\( \Omega_{1}(\xx) \) with \( \P\bigl( \Omega_{1}(\xx) \bigr) \geq 1 - 3 \ex^{-\xx} \)
for all pairs \( \mm < \md \in \MM \)  
\begin{EQA}
	\biggl| \frac{\dimAbt_{\mm,\md}}{\dimA_{\mm,\md}} - 1 \biggr|
	& \leq &
	\errSi_{\dimA} \, ,
	\\
	\errSi_{\dimA}
	& \eqdef &
	\| \Bias \|_{\infty}^{2} + 4 \, \xx_{\MM}^{1/2} \, \tdn^{2} \, \| \Bias \|
	+ 4 \xx_{\MM}^{1/2} \, \tdn + 4 \, \xx_{\MM} \, \tdn^{2} 
	+ \supepsi  ,
\label{dimAbootb}
\end{EQA}
where \( \dimAbt_{\mm,\md} = \Eb \| \xivb_{\mm,\md} \|^{2} \), 
\( \dimA_{\mm,\md} = \E \| \xiv_{\mm,\md} \|^{2} \), and
\( \xx_{\MM} = \xx + 2 \log(|\MM|) \).
\end{theorem}

The above results immediately imply all the oracle bounds for probabilistic loss of 
Section~\ref{SKnownnoice} with the obvious correction of the error terms.

\medskip

Now we discuss the sense of the required conditions for bootstrap validity. 
Our results are only meaningful and the bootstrap approximation is accurate 
if the values  \( \errSi_{2}(\xx) \) and
\( \sqrt{\dimp} \, \errSi_{0}(\xx) \) are small.
One easily gets
\begin{EQA}[c]
	\errSi_{2}(\xx) 
	\asymp 
	\sqrt{\dimp} \, \errSi_{0}(\xx)
	\leq 
	\CONST \dimp^{1/2} \bigl( \| \Bias \|_{\infty}^{2} + \dPsi ^{2} \, \| \Bias \|
	+ \dPsi 
	+ \supepsi \bigr) ,
\label{errSi2xxp12e0}
\end{EQA}
where \( \CONST \) is a generic notation for absolute constants and log-terms like 
\( \xxn, \xxp \) etc.
So, keeping the errors of bootstrap approximation small requires that the values \( \dPsi^{2} \, \dimp \),
\( \supepsi^{2} \, \dimp \), \( \| \Bias \|_{\infty}^{4} \, \dimp \), and
\( \dPsi^{2} \, \| \Bias \| \) are sufficiently small.
Now we spell this condition in the typical situation with \( \dPsi \asymp \sqrt{\dimp/n} \) and
\( \supepsi \asymp \sqrt{\dimp/n} \).
Then we need that
\( \dimp^{2} \log(n)/n \) is small.
Further, the bias component does not 
destroy the bootstrap validity result if the values \( \| \Bias \|_{\infty}^{4} \, \dimp \) 
and 
\( \dimp \, n^{-1} \, \| \Bias \| \leq \dimp \, n^{-1/2} \| \Bias \|_{\infty} \) are small.
If \( \fvs \) is H\"older-smooth with the parameter \( s \): 
\begin{EQA}
	\| \Bias \|_{\infty}
	& \leq &
	\CONST \dimp^{-s}
\label{f1nB2Cma}
\end{EQA}
%yielding \( \| \Bias \| \leq \CONST \dimp^{-s} n^{1/2} \),
then the bootstrap procedure is justified for \( s > 1/4 \) if 
\( \dimp = \dimn \to \infty \) but \( \dimn^{2}/n \to 0 \) as \( n \to \infty \).
We state one asymptotic result of this sort. 

\begin{corollary}
Assume the conditions of Theorem~\ref{TGaussbootB2} and
let \( \dimp = \dimn \) fulfill \( \dimn^{2} \log(n) / n \to 0 \) as \( n \to \infty \), 
and \eqref{f1nB2Cma} hold for \( s > 1/4 \).
Then the results of Theorem~\ref{TGaussbootB} and \ref{TGaussbootB2} apply with a small value 
\( \errSi_{n} = \bigl( \sqrt{\dimn} \, \errSi_{0}(\xxn) \bigr) \vee \errSi_{2}(\xxn) \to 0 \) as \( n \to 0 \).
\end{corollary}

\section{Simulations}
This section illustrates the performance of the proposed procedure by means of 
simulated examples. 
We consider a regression problem for an unknown univariate %periodic 
function on \( [0,1] \) with
unknown inhomogeneous noise.
The aim is to compare the bootstrap-calibrated procedure with the SmA procedure for the known noise
and with the oracle estimator.
We also check the sensitivity of the method to the choice of the presmoothing parameter \( \mres \). 

\begin{figure}
%\mygraphics{0.49}{0.2}{./plots/smallvariance.pdf}
%\mygraphics{0.49}{0.2}{./plots/mediumvariance.pdf}
%\mygraphics{0.49}{0.2}{./plots/largevariance.pdf}
 \begin{center}
 \includegraphics[width = 0.49\textwidth,height = 0.13\textheight]{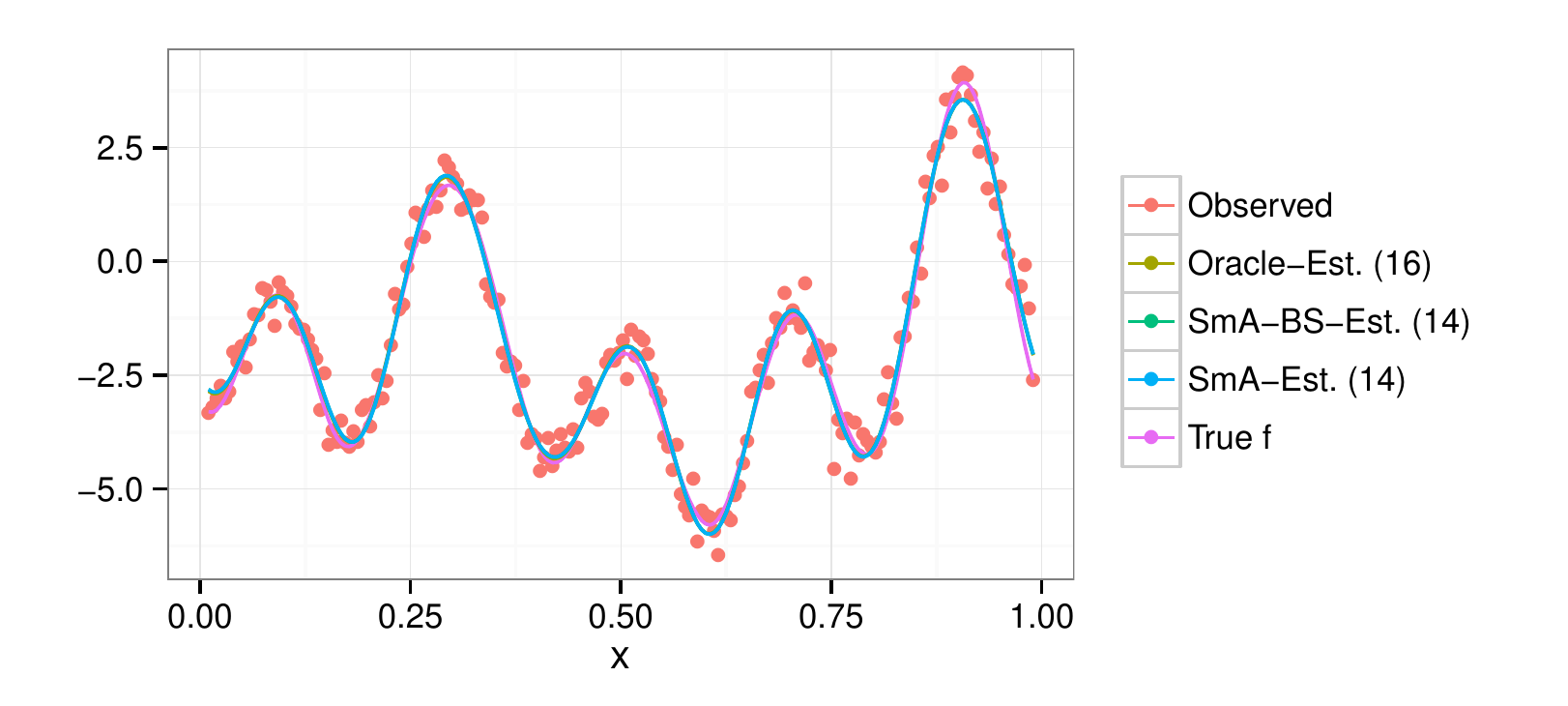}
 \includegraphics[width = 0.49\textwidth,height = 0.13\textheight]{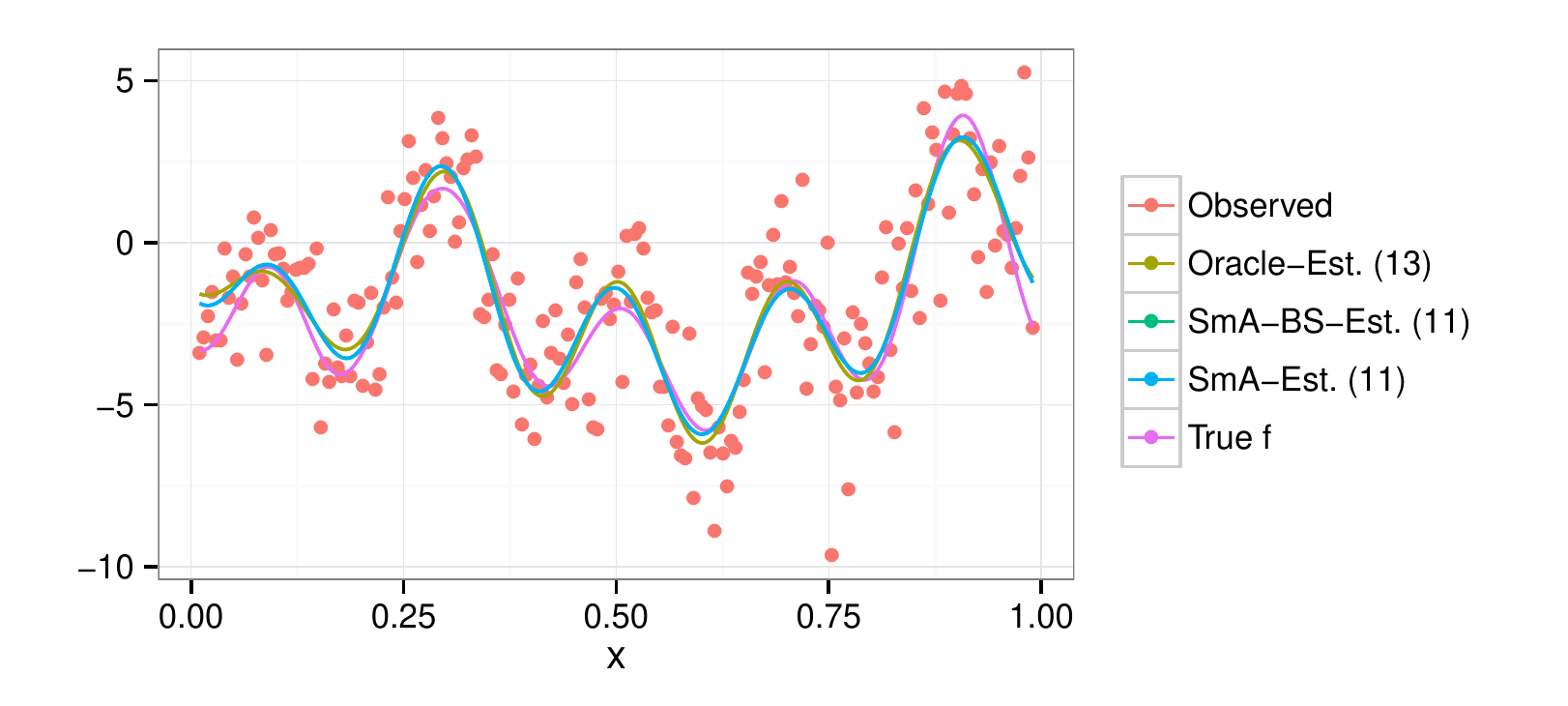}
 \includegraphics[width = 0.49\textwidth,height = 0.13\textheight]{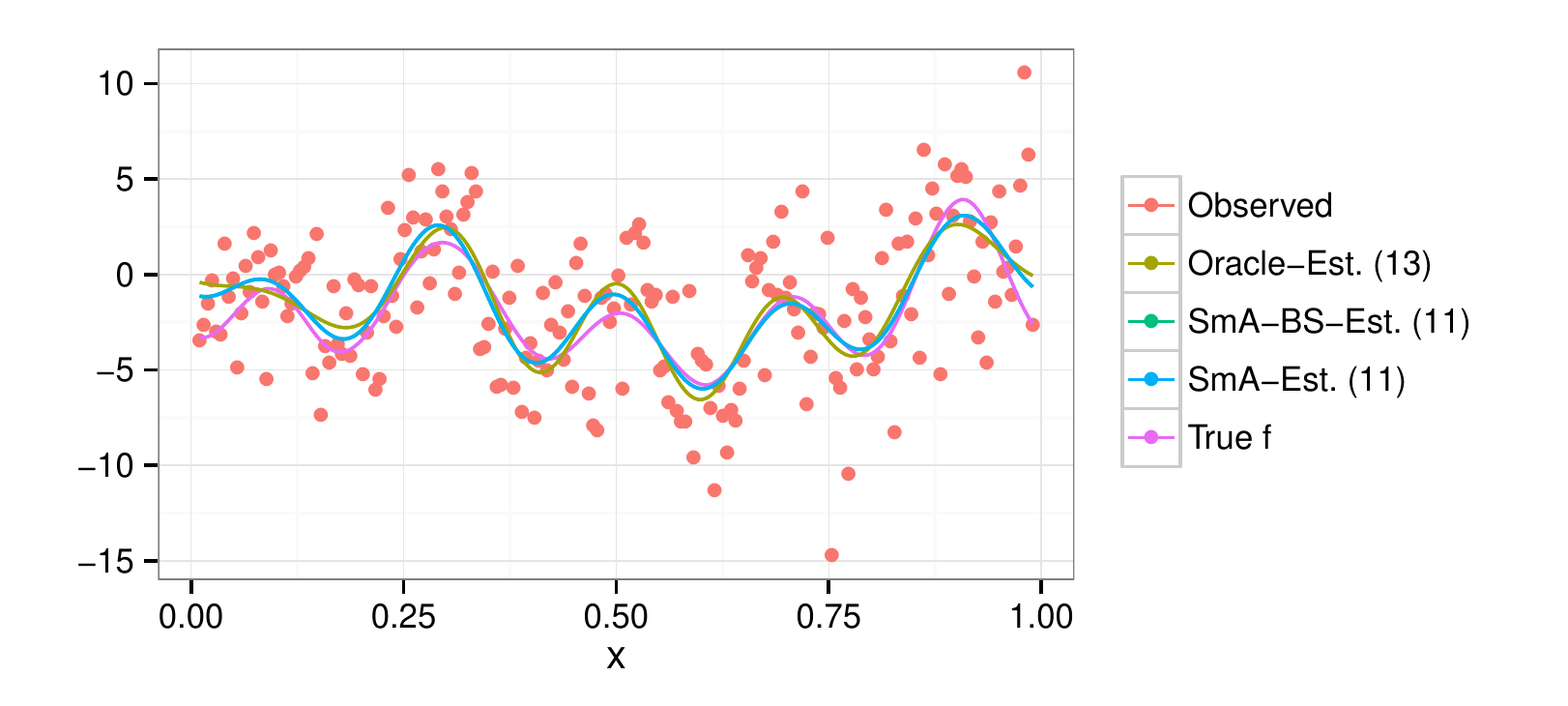}
 % lowvarheterofull.pdf: 423x216 pixel, 72dpi, 14.92x7.62 cm, bb=0 0 423 216
\end{center}
\caption{True functions and observed values plotted with oracle estimator, the known-variance SmA-Estimator (SmA-Est.) and the Bootstrap-SmA-Estimator (SmA-BS-Est.) for 3 different functions with different noise structure going from low noise to high noise. The numbers in parentheses indicate the chosen model dimension.}
\label{fig:fullvectorexample}
\end{figure}

We use a uniform design on \( [0,1] \) and the Fourier basis 
\( \{ \psi_{j}(x) \}_{j=1}^{\infty} \) 
for approximation of the regression function \( \fs \)
which is modelled in the form
\begin{EQA}[c]
 	\fs(x) = c_{1} \psi_{1}(x) + \ldots + c_{\dimp} \psi_{\dimp}(x),
\end{EQA}
where the \( (c_{j})_{ 1\le j \le \dimp} \) are chosen randomly: with
\( \gauss_{j} \) i.i.d. standard normal
\begin{EQA}[c]
 	c_{j} 
	= 
	\begin{cases}
        \gauss_{j}, & 1 \le j \le 10, \\
        \gauss_{j}/(j-10)^{2} , & 11 \le j \le 200.
	\end{cases}
\end{EQA} 
The noise intensity grows from low to high as \( x \) increases to one.
We use \( n_\text{sim-bs} = n_\text{sim-theo} = n_\text{sim-calib} = 1000 \) samples 
for computing the bootstrap marginal quantiles and the theoretical quantiles and for checking 
the calibration condition. 
The maximal model dimension is \( \mmmax = 37 \) and we also choose \( \mres = 20 \). 
The calibration is run with \( \xx = 2 \) and \( \alpb = 1 \). 
%We use a trigonometric basis for the construction of our projection estimators.

We start by considering examples for  \( W = \Psi_{n}^{\top} \), i.e. the estimation of the whole function vector with prediction loss. 
One can see in Figure \ref{fig:fullvectorexample} three examples with different intensity of the noise term comparing the Bootstrap-method to the oracle estimator and the known-variance SmA-Method.
Figure \ref{fig:residualvarplot} illustrates the dependence of the choice of the estimated dimension on our calibration dimension \( \mres \) and the sample size \( n \). 
We see that in the specific example we are considering, the sensitivity of the chosen dimension 
\( \tilde{\mm} \) on \( \mres \) decreases very fast. 
In the case \( n = 200 \), we have no variation in the choice of \( \tilde{\mm} \) with respect to 
\( \mres \). 
The oracles are respectively \( \ms = 12 \) for \( n = 100,200 \) and \( \ms = 10 \) for \( n = 50 \).
\begin{figure}[h]
%  	\mygraphics{0.49}{0.2}{./plots/residualex50.pdf}
%   	\mygraphics{0.49}{0.2}{./plots/residualex100.pdf}
%  	\mygraphics{0.49}{0.2}{./plots/residualex200.pdf}
% 	\mygraphics{0.49}{0.2}{./plots/residualvarplot.pdf}

 	\begin{center}
  	\includegraphics[width = 0.49\textwidth,height = 0.12\textheight]{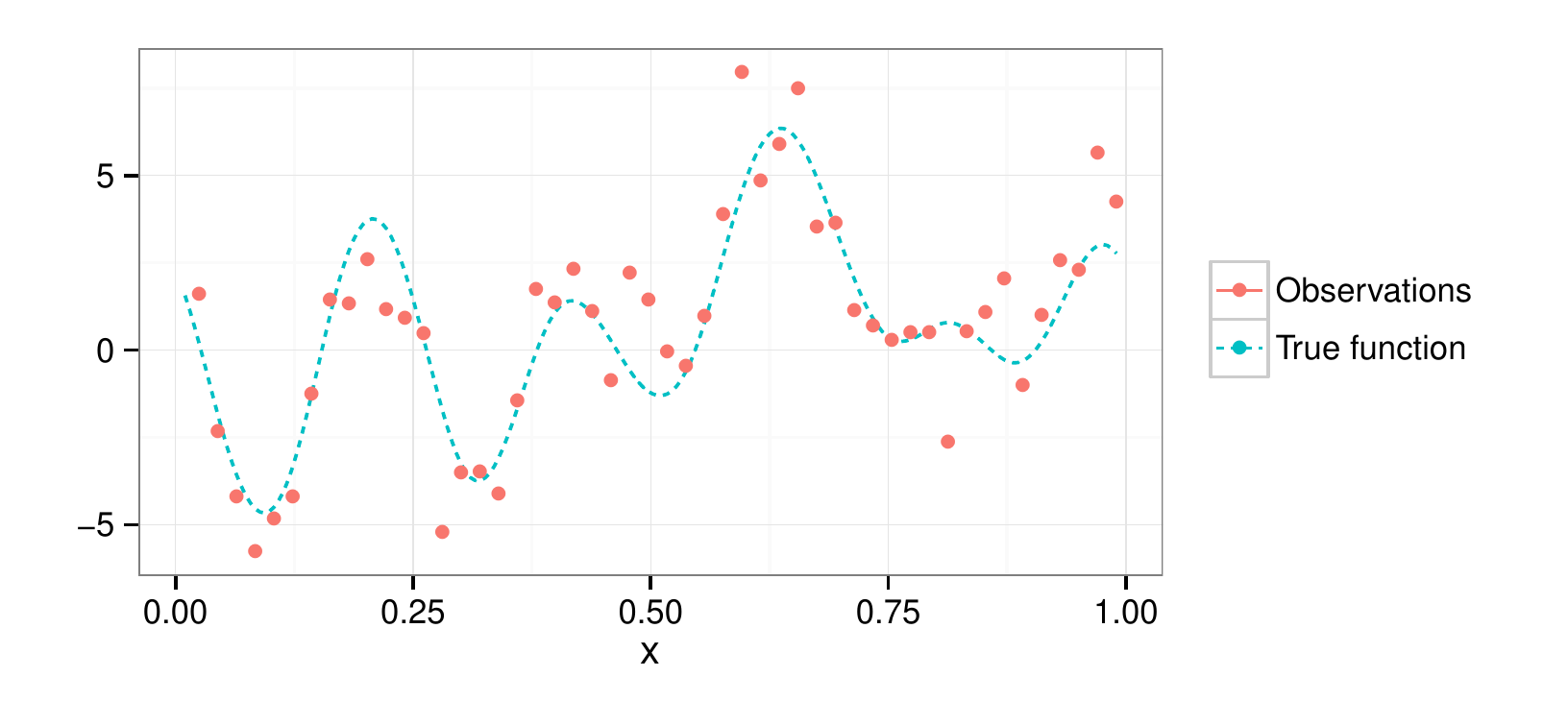}
   	\includegraphics[width = 0.49\textwidth,height = 0.12\textheight]{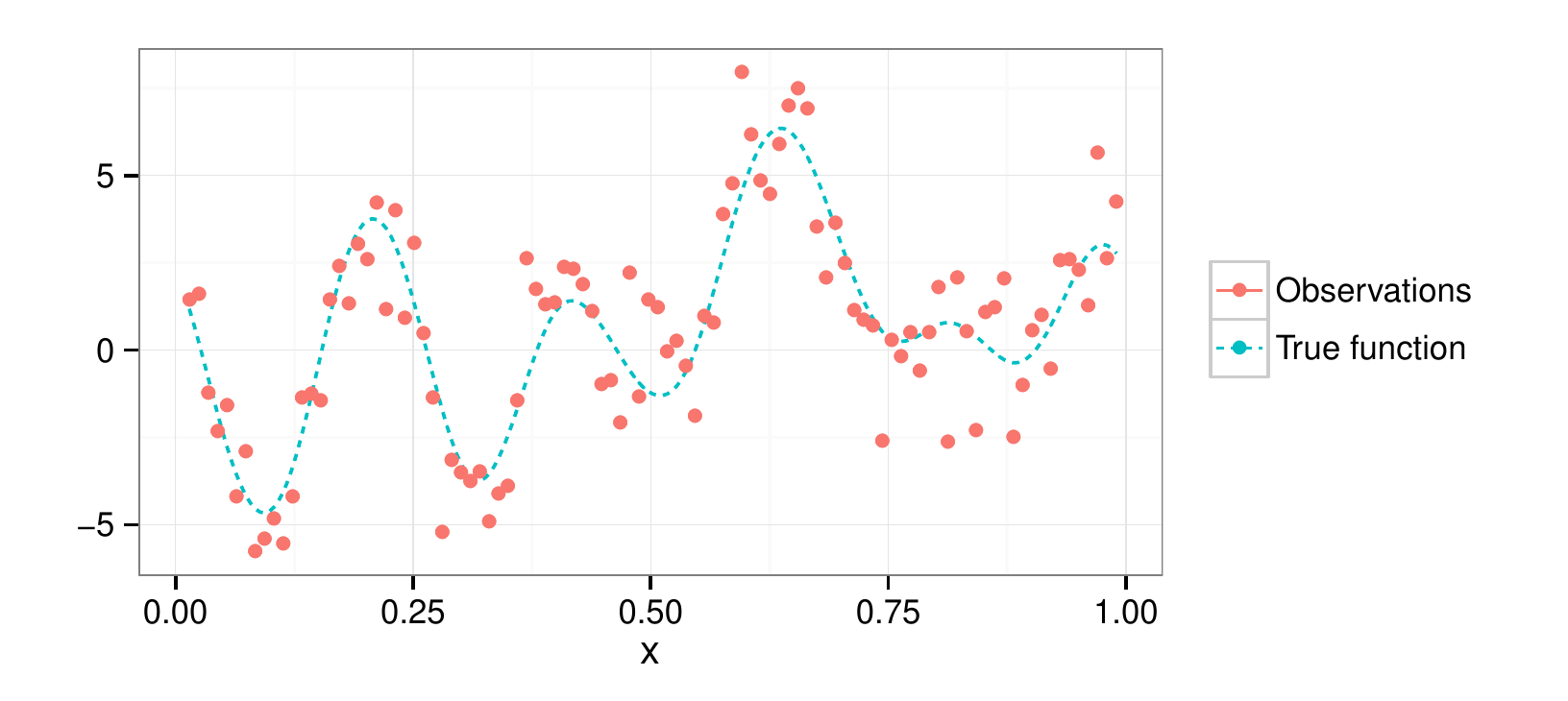}
  	\includegraphics[width = 0.49\textwidth,height = 0.12\textheight]{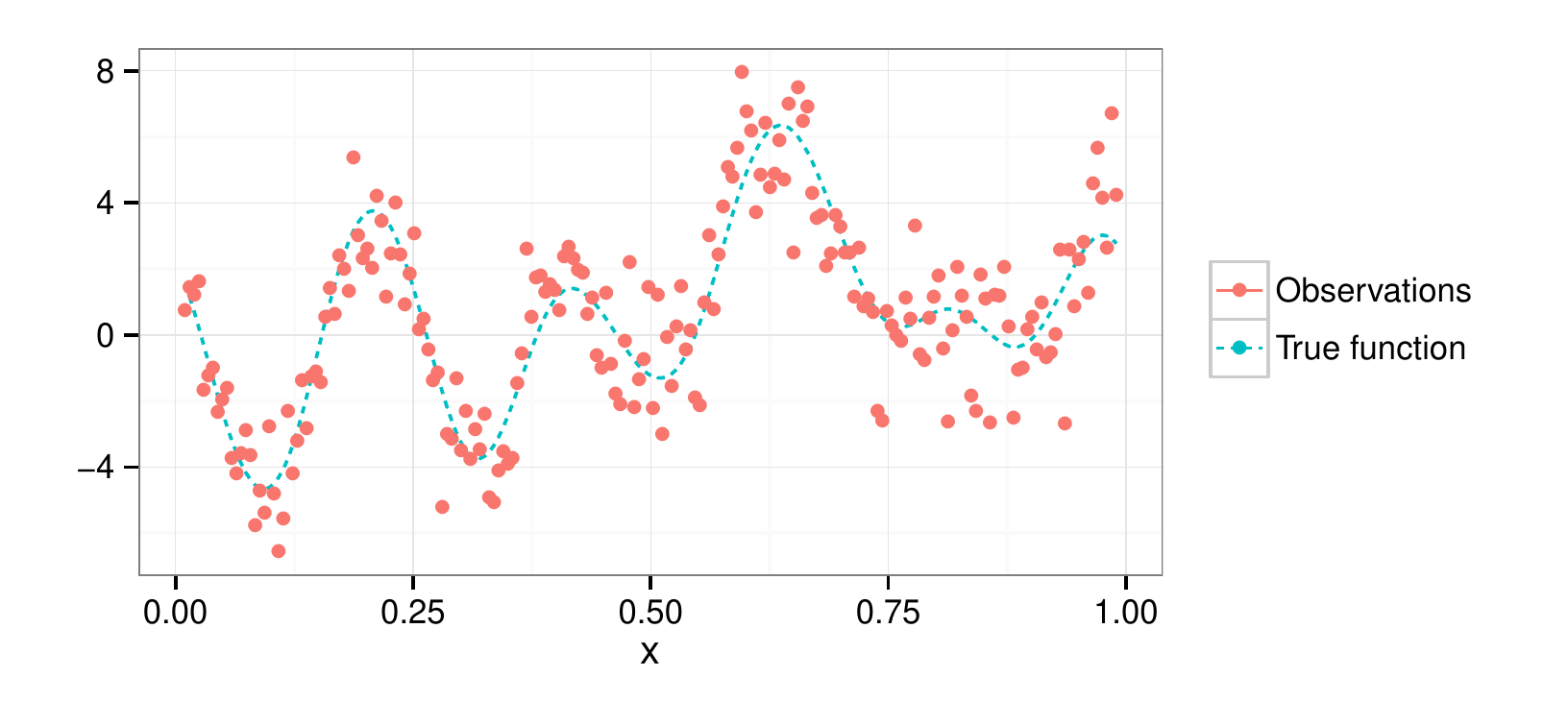}
 	\includegraphics[width = 0.49\textwidth,height = 0.12\textheight]{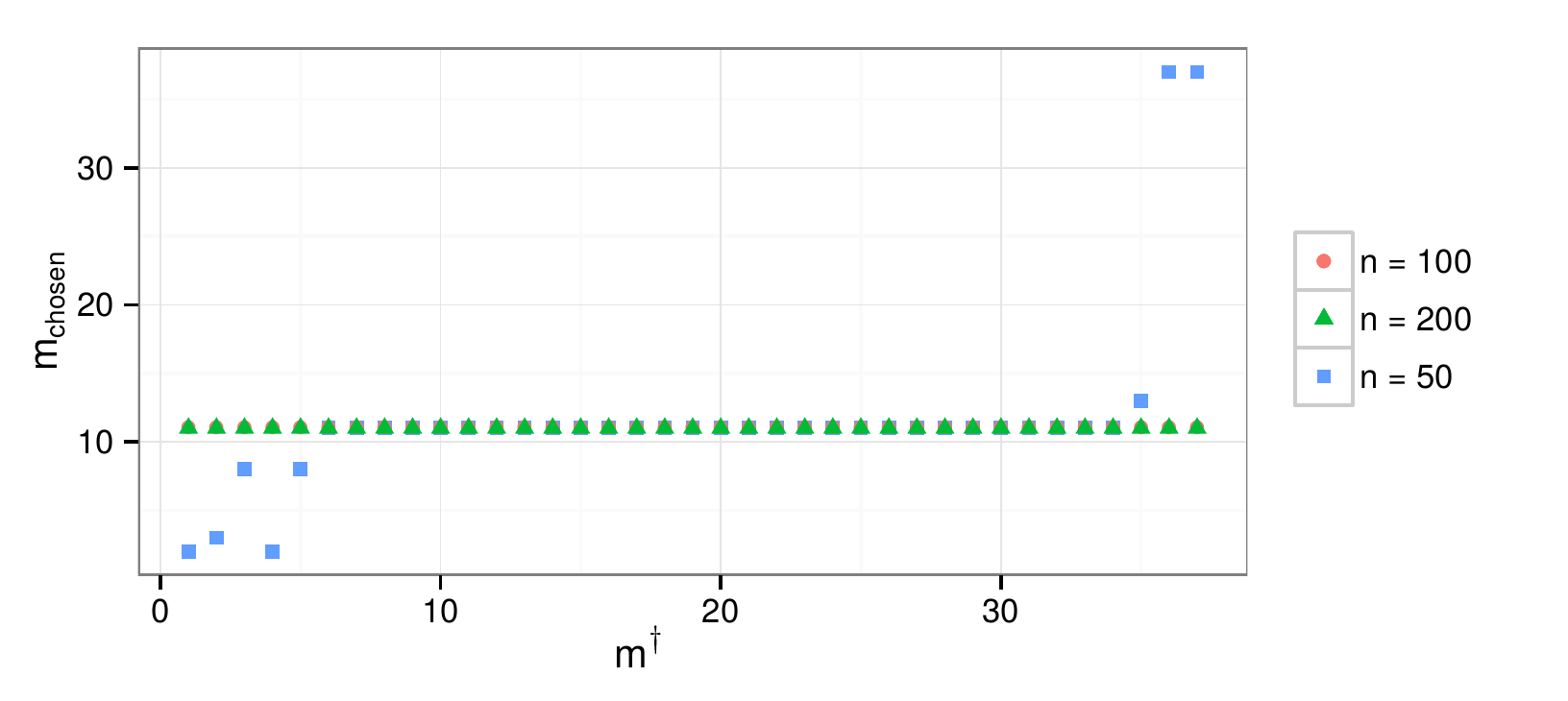}
	\end{center}
\caption{The first three plots show an exemplary function with \( n =50,100,200 \) observations. 
The right plot shows the \( \hat{\mm} \) chosen by the Bootstrap-SmA-Method as a function of the calibration dimension \( \mres \) and the number of observations. }
\label{fig:residualvarplot}
\end{figure}
\ifAoS{ % full
We also want to compare the true quantiles and their bootstrap substitute. 
Figure \ref{fig:ratio3dplot} plots the ratios of quantiles for all possible comparisons 
\( (\mm_{1},\mm_{2}) \) for the same function as before. Here we see that there is, as one would expect, still significant variation in the quantile ratios for small differences \( |\mm_{1} - \mm_{2}| \). Nonetheless the method works very well as seen in Fig. \ref{fig:residualvarplot}, but the variability in the ratios implies the possibility to stabilize the procedure even more by introducing some smoothing scheme for the quantiles.
\begin{figure}
 \begin{center}
  \includegraphics[width = 0.5\textwidth]{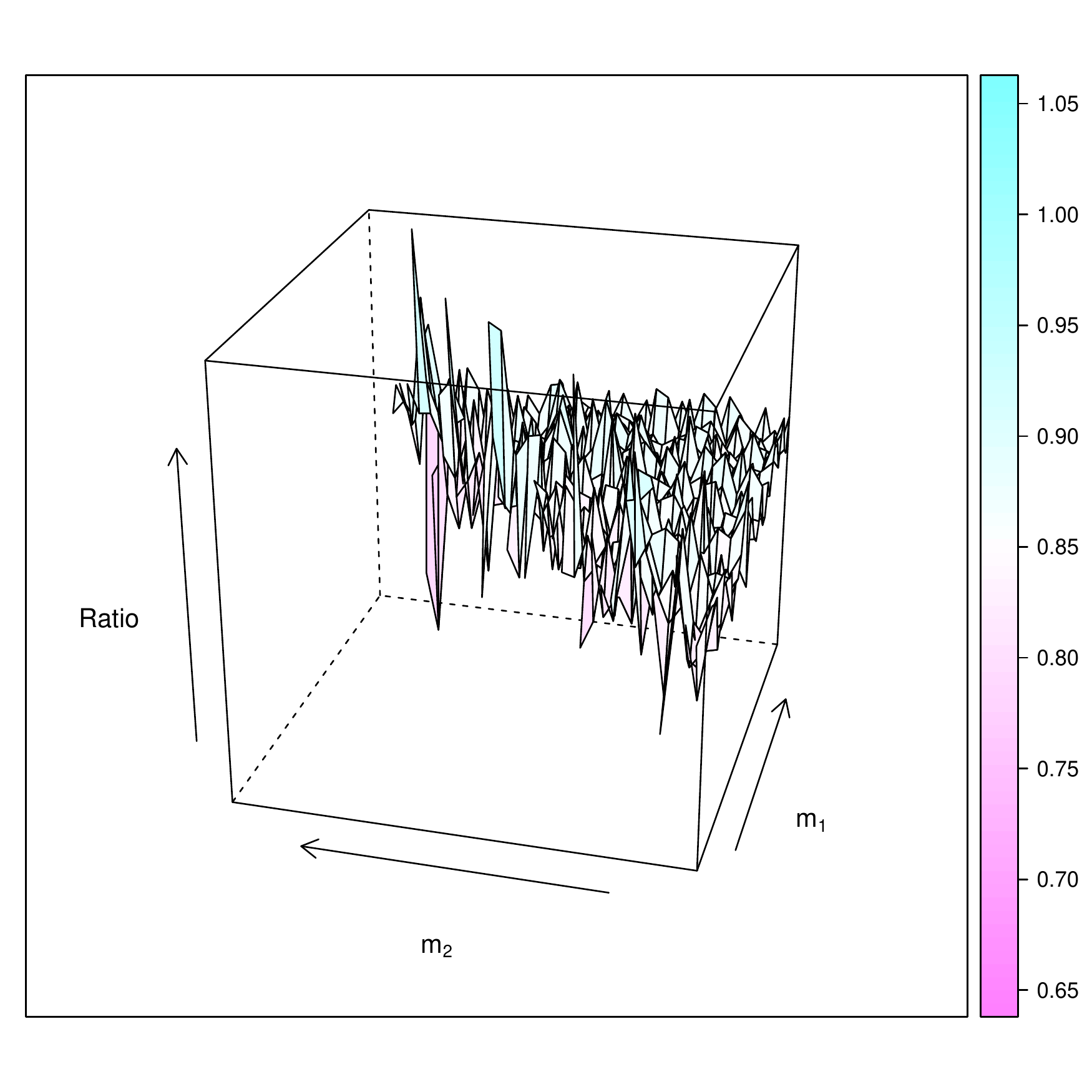}
\end{center}
\caption{Ratio of quantiles \( |\ZZbt_{\mm_{1}, \mm_{2}}/\ZZ_{\mm_{1},\mm_{2}}|^{2} \) for \( \mres = 20 \) and \( n = 200 \) with the data and true function as in Fig. \ref{fig:residualvarplot}.}
%of \( \accessiblestatistic{\mm_{1},\mm_{2}}^{2} \) and 
%\( {\accessiblestatisticbs{\mm_{1},\mm_{2}}}^{\mkern-35mu 2}\mkern20mu \) 
\label{fig:ratio3dplot}
\end{figure}
}{ % annals
}

%Again the dependency on the choice of \( \mres \) can be studied as seen in 
Figure \ref{fig:ratioplot} again demonstrates the dependence of the ratios on \( \mres \).
It is remarkable that the ratio is varying very slowly above \( \ms = 12 \).
\begin{figure}
 \begin{center}
  \includegraphics[width = 0.6\textwidth,height = 0.13\textheight]{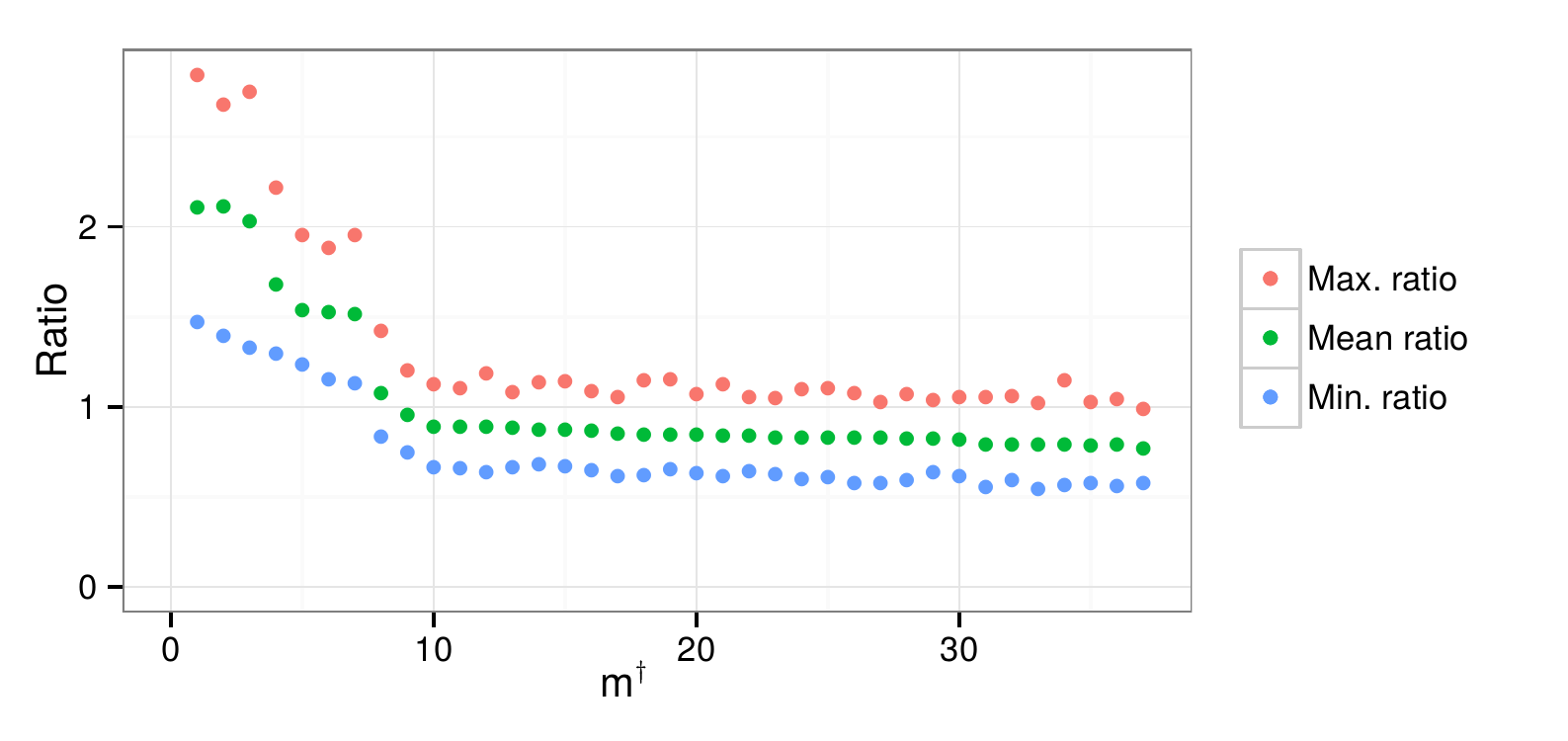}
\end{center}
\caption{Maximal, minimal and mean ratio of the bootstrap and theoretical tail functions at \( x=2 \), 
\( |\ZZbt_{\mm_{1}, \mm_{2}}/\ZZ_{\mm_{1},\mm_{2}}|^{2} \) 
%of \( \accessiblestatistic{\mm_{1},\mm_{2}}^{2} \) and 
%\( {\accessiblestatisticbs{\mm_{1},\mm_{2}}}^{\mkern-35mu 2}\mkern20mu \) 
as a function of \( \mres \).}
\label{fig:ratioplot}
\end{figure}
We also give the results on the simulation of \( n_{\text{hist}} = 100 \) repeated applications of the method to the same true underlying function observed with different realizations of the errors in Figure \ref{fig:histmethod}.
\begin{figure}
 \begin{center}
  \includegraphics[width = 0.49\textwidth,height = 0.12\textheight]{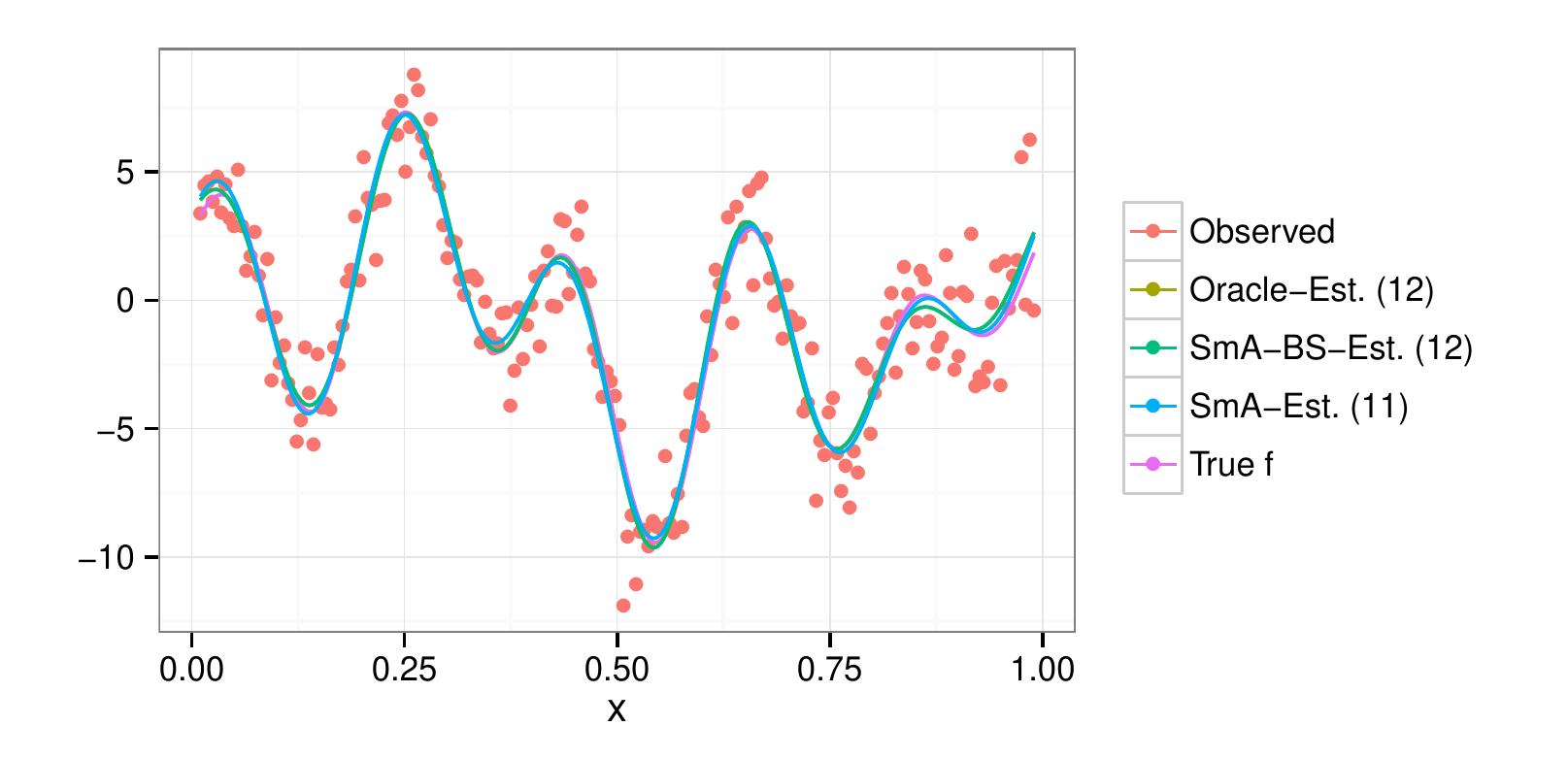}
  \includegraphics[width = 0.49\textwidth,height = 0.12\textheight]{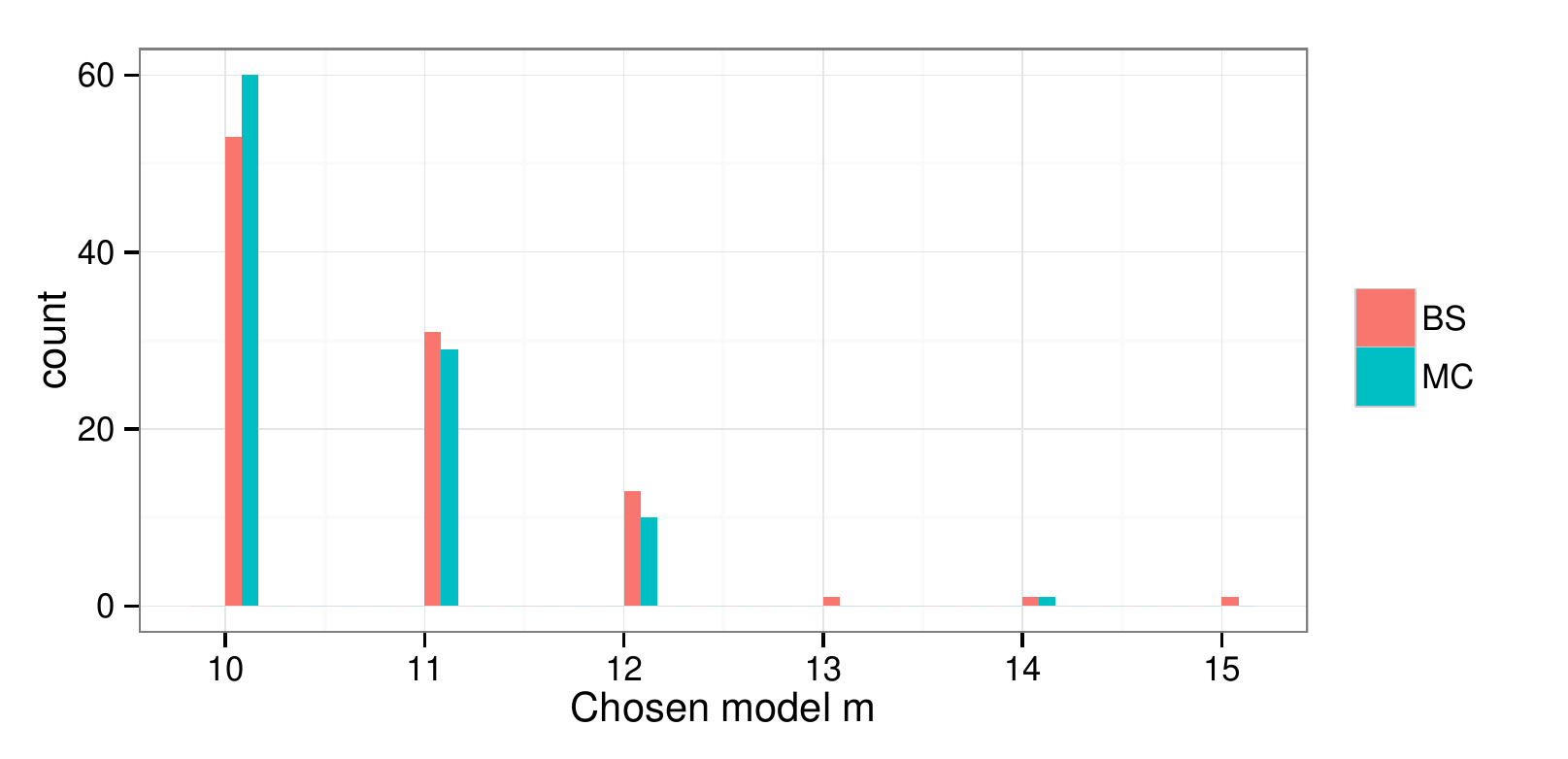}
\end{center}
\caption{In the left plot, the true function and observed values are plotted for one realization together with the oracle estimator, the known-variance SmA-Estimator (SmA-Est.) and the Bootstrap-SmA-Estimator (SmA-BS-Est.). The numbers in parentheses indicate the chosen model dimension. In the right plot, histograms for the selected model are given for the bootstrap  (BS) and the known-variance method (MC) for repeated observations of the same underlying function with a  simulation size  \( n_{\text{hist}} = 100 \).}
\label{fig:histmethod}
\end{figure}

The case of the estimation of the first derivative is similar. 
Figure \ref{fig:derivativeexample} shows the numerical results for estimation of the 
derivative in the same model as above. One can see that the bootstrap-version
of the SmA-procedure is again competitive with the procedure based on a known noise structure and the method does a good job of mimicking the oracle.
\begin{figure}
 \begin{center}
 \includegraphics[width = 0.49\textwidth,height = 0.13\textheight]{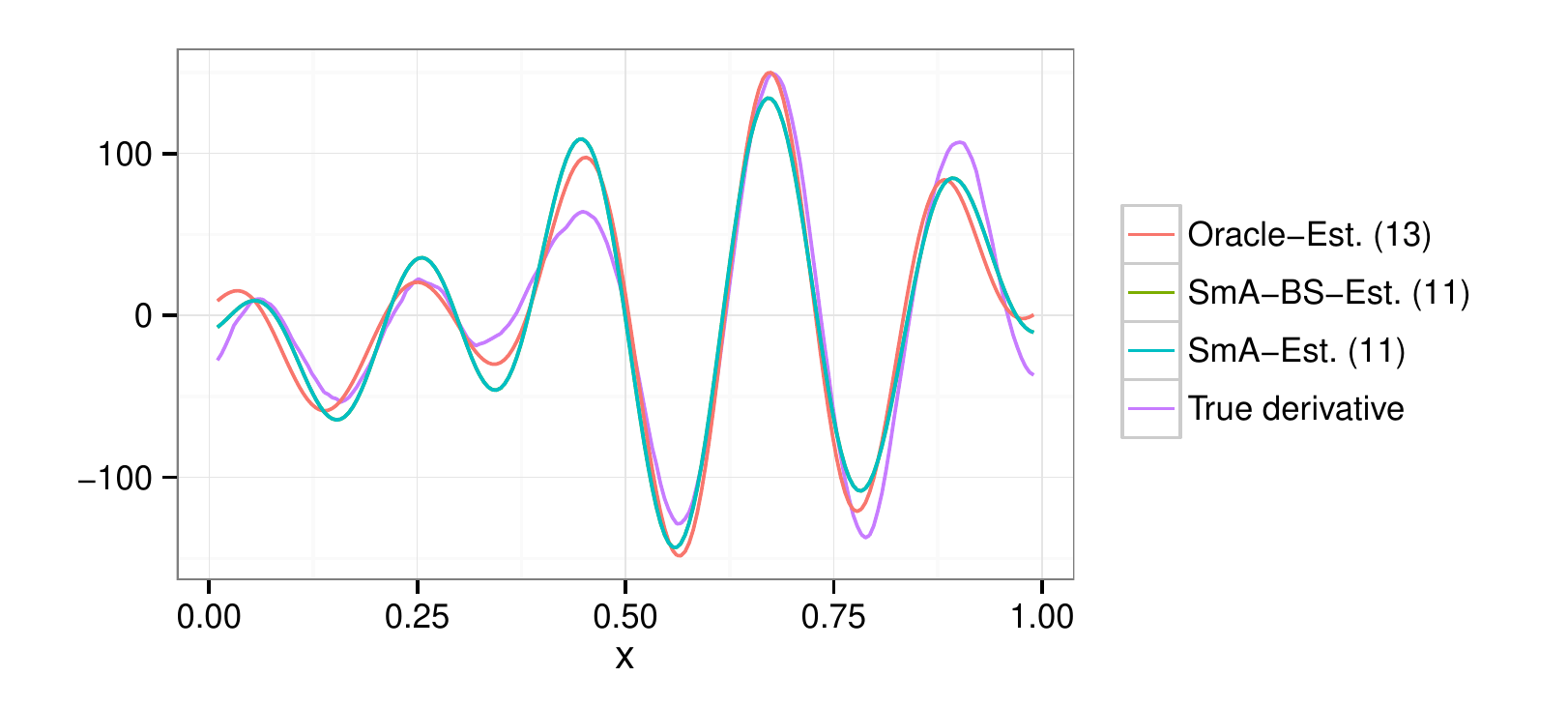}
 \includegraphics[width = 0.49\textwidth,height = 0.13\textheight]{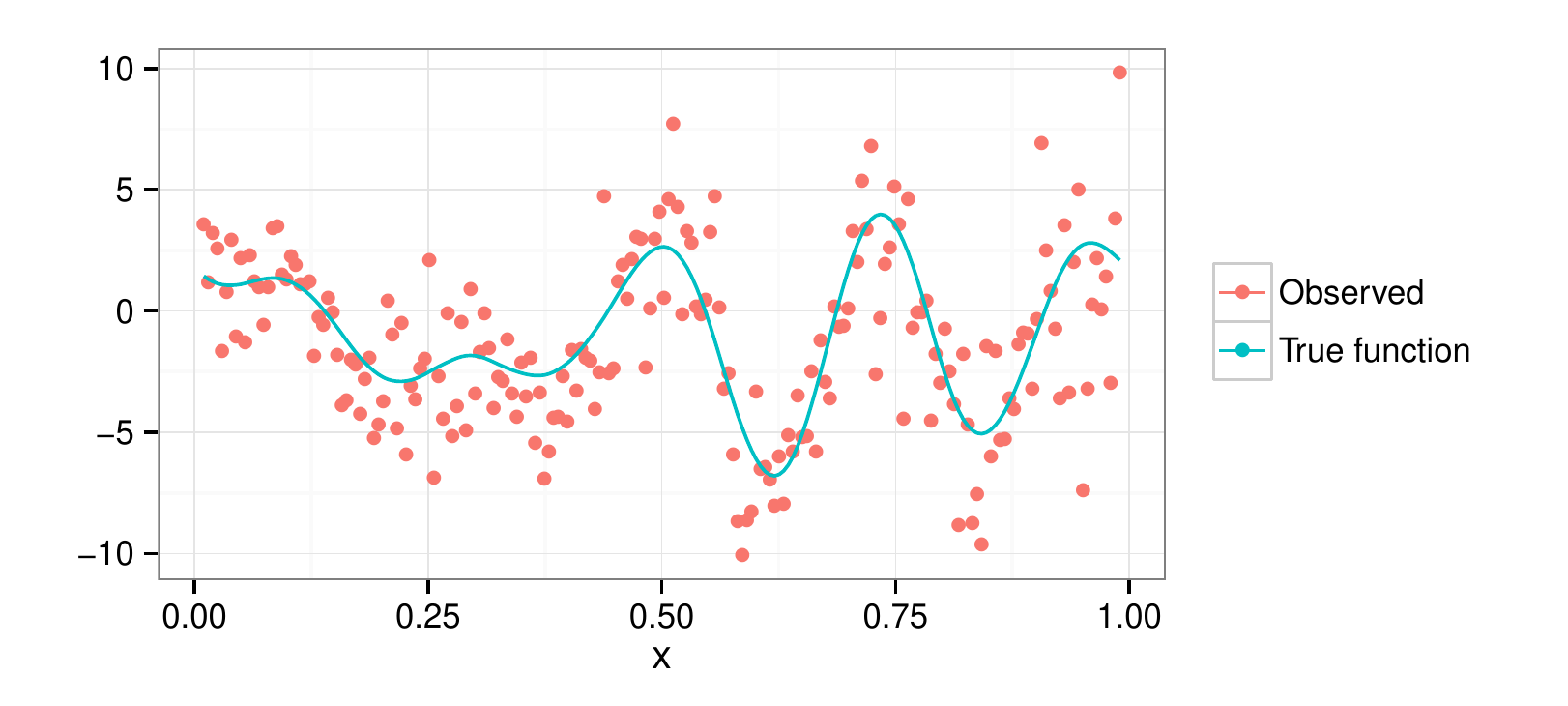}
 \includegraphics[width = 0.49\textwidth,height = 0.13\textheight]{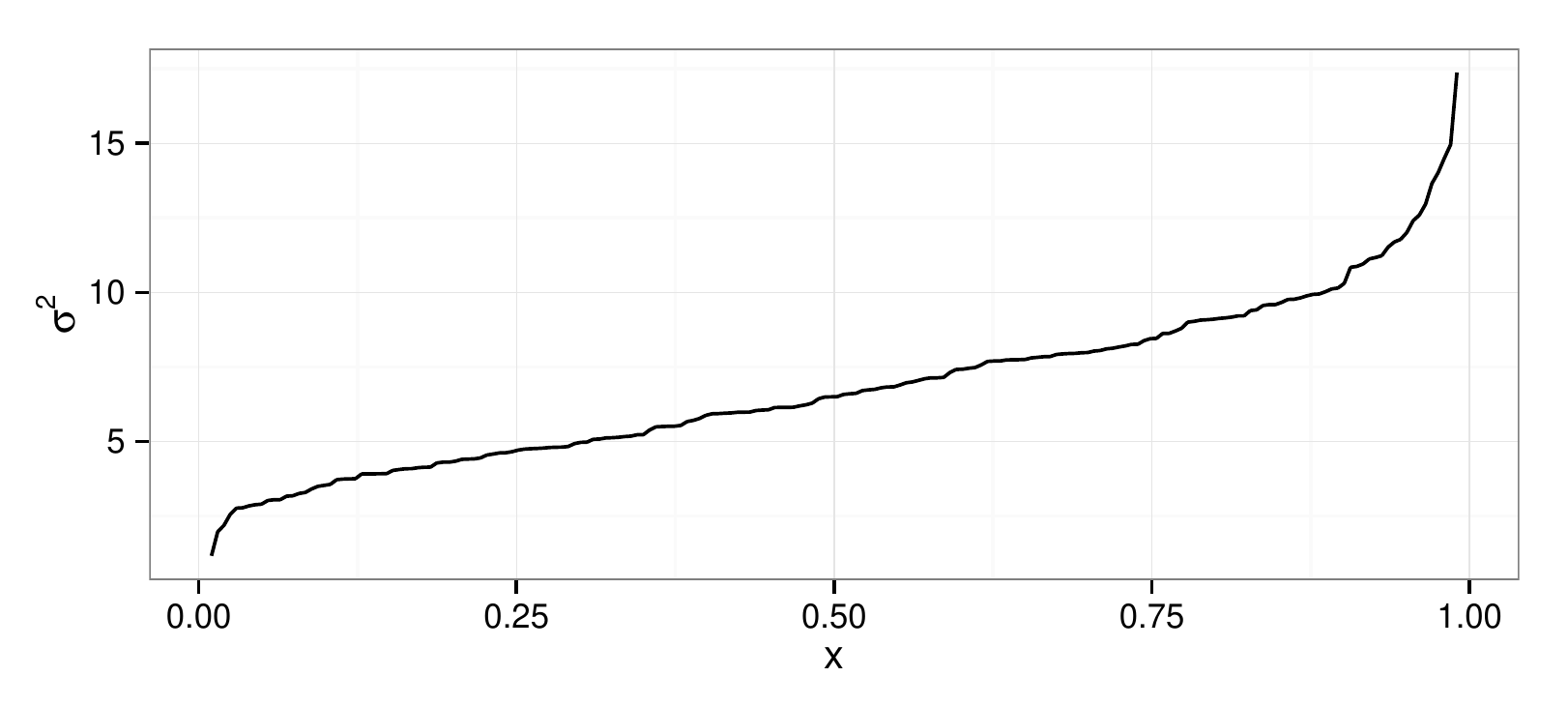}
 % lowvarheterofull.pdf: 423x216 pixel, 72dpi, 14.92x7.62 cm, bb=0 0 423 216
\end{center}
\caption{The upper left plot shows the true derivative, the oracle estimator, the known-variance SmA-Estimator (SmA-Est.) and the Bootstrap-SmA-Estimator (SmA-BS-Est.). 
The upper right plot shows the true function and the observations and in the lower plot one can find the standard deviation of the errors.}
\label{fig:derivativeexample}
\end{figure}

One can conclude that the proposed procedure is really universal and demonstrates 
a very good performance in various settings.

\appendix
\section{Proofs}
The appendix collects the proofs of announced results.

\subsection{Proof of Theorem~\ref{ToracleSmA}}
%\begin{proof}
The propagation property \eqref{Ppropmslin} claims that the oracle model \( \ms \)
will be accepted with high probability.
This yields that the selected model is not larger than \( \ms \), that is,
\( \hat{\mm} \leq \ms \) with a probability at least \( 1 - \ex^{-\xx} \).
Below we consider only this event. Let \( \mm \in \MMd(\ms) \).
Acceptance of \( \mm \) requires in particular that 
\( \Tmd_{\ms,\mm} \leq \ZZ_{\ms,\mm} \).
The representation 
\( \Tmd_{\ms,\mm} = \| \bias_{\ms,\mm} + \xiv_{\ms,\mm} \| \) implies
\begin{EQA}
%	&&
%	\nquad
	\P\bigl( 
		\Tmd_{\ms,\mm} 
		< \ZZ_{\ms,\mm} 
	\bigr)
%	\\
	& \leq &
	\P\bigl( 
		\| \xiv_{\ms,\mm}\| 
		> \| \bias_{\ms,\mm} \| - \ZZ_{\ms,\mm} 
	\bigr) .
\label{Pmsmmzzm}
\end{EQA} 
Under \eqref{MMxxdeflin} this yields
%Acceptance of the model \( \mm < \ms \) assumes that the check of 
%\( \Tmd_{\ms,\mm} \) does not fail.
\begin{EQA}
	\P\bigl( \mm \text{ is accepted} \bigr)
	& \leq &
	\P\bigl( 
		\bigl\| \bias_{\ms,\mm} + \xiv_{\ms,\mm} \bigr\|
		\leq \ZZ_{\ms,\mm}
	\bigr)
	\\
	& \leq &
	\P\bigl( 
		\bigl\| \xiv_{\ms,\mm} \bigr\|
		\geq \zq_{\ms,\mm}(\xxu)
	\bigr)
	\leq 
	\ex^{-\xxu} .
\label{Pmmacceptxx}
\end{EQA}
If the lower bound on the bias is fulfilled for all \( \mm \in \MMc \), 
then \eqref{Pmmacceptxx} helps to bound the probability of the event  
\( \{ \hat{\mm} \in \MMc \} \): 
\begin{EQA}
	\P\bigl( \hat{\mm} \in \MMc \bigr)
	& \leq &
	\sum_{\mm \in \MMc} \P\bigl( 
		\bigl\| \bias_{\ms,\mm} + \xiv_{\ms,\mm} \bigr\|
		< \ZZ_{\ms,\mm} 
	\bigr)
	\leq 
%	\\
%	& \leq &
	\sum_{\mm \in \MMc} \ex^{-\xxu} \leq \ex^{-\xx} .	
\label{PmmMMxxu}
\end{EQA}
Therefore, the probability that the SmA-selector picks up a value \( \mm > \ms \) 
or \( \mm \in \MMc \) is very small:
\begin{EQA}
	\P\Bigl( \hat{\mm} \in \MMu(\ms) \cup \MMc \Bigr) 
	&\leq & 
	2 \ex^{-\xx} .
%	\\
% 	\P\bigl( \hat{\mm} < \ms, \mm\in \MMc \bigr)
%	& \leq &
%	\ex^{-\xx} .
\label{PnotinMMxx}
\end{EQA}
It remains to study the case when \( \hat{\mm} = \mm \in \MMi = \MMd(\ms) \setminus \MMc \). 
We can use that \( \hat{\mm} \) is accepted, which implies by definition
\begin{EQA}
	\Tmd_{\ms,\mm} 
	& = & 
	\bigl\| 
		\tilde{\tarv}_{\mm} - \tilde{\tarv}_{\ms} 
	\bigr\|
	\leq
	\ZZ_{\ms,\mm} \, .
\label{Tmdmsmm1alp}
\end{EQA}
This yields \eqref{QLtmmtmshh}.
The bound \eqref{Proboralin} now follows by the triangle inequality.
%\end{proof}

\subsection{Proof of Proposition~\ref{Tpayment}}

Below we use the deviation bound 
\ifAoS{\eqref{Pxiv2dimAxx12} for a Gaussian quadratic form from Theorem~\ref{TexpbLGA}.}
{ %annals 
(C.2) for a Gaussian quadratic form from Theorem~C.1 in \suppSW.}
Note that similar results are available for non-Gaussian quadratic forms under 
exponential moment conditions; see e.g. \cite{SP2011}.
The result 
\ifAoS{ % full 
\eqref{Pxiv2dimAxx12}}
{ % annals
(C.2)}
combined with the Bonferroni correction \( \qq_{\md} = \log(|\MMu(\md)|) \leq \log(|\MM|) \) 
yields the following upper bound for the critical values \( \ZZ_{\mm,\md} \):
\begin{EQA}
	\ZZ_{\mm,\md}
	& \leq &
	\zq_{\mm,\md}(\xx + \qq_{\md}) + \alpb \dimA_{\mm,\md}^{1/2}
	\\
	&\leq &
	(1 + \alpb) \sqrt{\dimA_{\mm,\md}} 
		+ \sqrt{2 \supA_{\mm,\md} \, \{ \xx + \log(|\MMu(\md)|) \} } 
	\\
	&\leq &
	(1 + \alpb) \sqrt{\dimA_{\mm,\md}} 
		+ \sqrt{2 \supA_{\mm,\md} \, \{ \xx + \log(|\MM|) \} } .
\label{ZZmmmdubA}
\end{EQA}
For the payment for adaptation \( \ZZs_{\ms} \), the result \eqref{ZZmmmdubA} 
and the monotonicity condition 
\( \dimA_{\ms,\mm} \leq \dimA_{\ms,\mmmin} \leq \dimA_{\ms} \) and 
\( \supA_{\ms,\mm} \leq \supA_{\ms,\mmmin} \leq \supA_{\ms} \)
imply the following upper bound:
\begin{EQA}
	\ZZs_{\ms}
	&\leq &
	(1 + \alpb) \sqrt{\dimA_{\ms,\mmmin}} 
		+ \sqrt{2 \supA_{\ms,\mmmin} \, \{ \xx + \log(|\MMd(\ms)|) \} } 
%	\\
%	&\leq &
%	(1 + \alpb) \sqrt{\dimA_{\ms}} + \sqrt{\xx + \log(|\MMd(\ms)|) }
	\\
	&\leq &
	(1 + \alpb) \sqrt{\dimA_{\ms}} + \sqrt{2 \supA_{\ms} \, \{ \xx + \log(|\MM|) \} } \, 
\label{ZZbmsupxx}
\end{EQA}
which yields the claim.

\ifAoS{ % full version}
\subsection{Proof of Theorem~\ref{TSmApoly}}

%\begin{proof}
The result will be proved in two steps.
First we bound the risk on the set \( \hat{\mm} > \ms \):
\begin{EQA}
	\E \bigl\{ 
		\| \hat{\tarv} - \tarvs \|^{2} \Ind\bigl( \hat{\mm} > \ms \bigr) 
	\bigr\}
	& \leq &
	2 \alpd_{\ms} \riskt_{\ms} \, .
\label{EQLhts2Ige}
\end{EQA}
Then we consider the region \( \hat{\mm} < \ms \)
and prove an oracle inequality
\begin{EQA}
	\| \hat{\tarv} - \tilde{\tarv}_{\ms} \| \Ind\bigl( \hat{\mm} < \ms \bigr)
	& \leq &
	\ZZs_{\ms} 
\label{QLttmmtmsPLo}
\end{EQA}
and the oracle bound \eqref{EQLttts2ms}.
We start by proving \eqref{EQLhts2Ige}.
Let us fix \( \mm \in \MMu(\ms) \) and \( \mc \geq \mm \).
The definition \eqref{msdeflinPL} of the oracle \( \ms \) and the formula
\eqref{zqmmmdPL} for the critical value \( \ZZ_{\mc,\mmprec} \) implies for 
the test statistic \( \Tmd_{\mc,\mmprec} = \| \xiv_ {\mc,\mmprec} + \bias_ {\mc,\mmprec} \| \) \begin{EQA}
	\bigl\{ \Tmd_{\mc,\mmprec} > \ZZ_{\mc,\mmprec} \bigr\}
	& \subseteq &
	\bigl\{ \| \xiv_{\mc,\mmprec} \| > \zq_{\mc,\mmprec}(\xx_{\mmprec}) \bigr\} .
\label{Tmdmcmm1}
\end{EQA}
Now we can bound the risk of \( \hat{\tarv} \) on the set \( \hat{\mm} > \ms \).
We use that for \( \hat{\mm} = \mm > \ms \) in view of \eqref{biasmmmsPL}
\begin{EQA}
	\| \hat{\tarv} - \tarvs \|^{2} %\Ind(\hat{\mm} = \mm)
	&=& 
	\| \tilde{\tarv}_{\mm} - \tarvs \|^{2}
	=
	\| \xiv_{\mm} + \bias_{\mm} \|^{2}
	\\
	& \leq &
	2 \| \xiv_{\mm} \|^{2} + 2 \| \bias_{\mm} \|^{2}
	\leq
	2 \| \xiv_{\mm} \|^{2} + 2 \| \bias_{\ms} \|^{2}
\label{QLhtts2}
\end{EQA}
and it holds by \eqref{riskrmmPL2} and monotonicity \( \dimA_{\mm} > \dimA_{\ms} \)
\begin{EQA}
	&& \nquad
	\E \bigl\{ 
		\| \hat{\tarv} - \tarvs \|^{2} \Ind\bigl( \hat{\mm} > \ms \bigr) 
	\bigr\}
	\\
	& \leq &
	2 \sum_{\mm \in \MMu(\ms)} 
	\E \bigl\{\bigl( \| \xiv_{\mm} \|^{2} + \| \bias_{\ms} \|^{2} \bigr) 
		\Ind(\hat{\mm} = \mm) 
	\bigr\}
	\\
	& \leq &
	2 \sum_{\mm \in \MMu(\ms)} 
	\E \bigl\{\bigl( \| \xiv_{\mm} \|^{2} + \| \bias_{\ms} \|^{2} \bigr) 
		\Ind(\mmprec \text{ is rejected}) 
	\bigr\}
	\\
	& = &
	2 \sum_{\mm \in \MMu(\ms)} 
		\E \biggl[ (\| \xiv_{\mm} \|^{2} + \| \bias_{\ms} \|^{2}) 
		\Ind\Bigl( 
			\max_{\mc \in \MMu(\mm)} 
			\Bigl\{ \| \xiv_{\mc,\mmprec} \| - \zq_{\mc,\mmprec}(\xx_{\mm}) \Bigr\} > 0 
		\Bigr) 
	\biggr]
	\\
	& \leq &
	2\sum_{\mm \in \MMu(\ms)} \alp_{\mm} (\dimA_{\mm} + \| \bias_{\ms} \|^{2})
	\leq 
	2 \alpd_{\ms} \bigl( \dimA_{\ms} + \| \bias_{\ms} \|^{2} \bigr)
	=
	2 \alpd_{\ms} \riskt_{\ms} .
\label{riskmmgems}
\end{EQA}
Here we have used that \eqref{alpmmalpb} and \( \dimA_{\mm} \geq \dimA_{\ms} \) 
imply \( \sum_{\mm \in \MMu(\ms)} \alp_{\mm} \leq \alpd_{\ms} \).
This completes the proof of \eqref{EQLhts2Ige}.

In the situation when \( \hat{\mm} = \mm < \ms \), we can use the stability property:
as \( \mm \) is accepted, it holds
\begin{EQA}
	\| \tilde{\tarv}_{\mm} - \tilde{\tarv}_{\ms} \| \Ind(\hat{\mm} = \mm)
	& \leq &
	\ZZ_{\ms,\mm} \, ,
\label{QLttmmtmsPL}
\end{EQA}
which implies \eqref{QLttmmtmsPLo} by definition of \( \ZZs_{\ms} \).
This yields
\begin{EQA}
	\E \bigl\| \hat{\tarv} - \tarvs \bigr\|^{2}
	& \leq &
	2 \alpd_{\ms} \riskt_{\ms} 
	+ \E \bigl\{ 
		\bigl\| \hat{\tarv} - \tarvs \bigr\|^{2} \Ind(\hat{\mm} < \ms) 
	\bigr\}
	\\
	& \leq &
	2 \alpd_{\ms} \riskt_{\ms} 
	+ \E \bigl( \| \tilde{\tarv}_{\ms} - \tarvs \| + \ZZs_{\ms} \bigr)^{2}
	\\
	& \leq &
	2 \alpd_{\ms} \riskt_{\ms} 
	+ \bigl( \riskt_{\ms}^{1/2} + \ZZs_{\ms} \bigr)^{2}
\label{EQLttts2msf}
\end{EQA}
as required.
%\end{proof}

\subsection{Proof of Proposition~\ref{Talpmxxm}}

%\begin{proof}
Observe first that the choice \( \alp_{\mm} = (\dimA_{\mm}/\dimA_{\mmmin})^{-1-a} \) yields
\begin{EQA}
	\sum_{\mm \in \MMu(\ms)} \alp_{\mm} \dimA_{\mm}
	& \leq &
	\dimA_{\mmmin}^{1+a} 
	\sum_{\mm \in \MMu(\ms)} \dimA_{\mm}^{-a}
	\leq 
	\CONST \dimA_{\ms}^{-a} \, \dimA_{\mmmin}^{1+a} 
	=
	\CONST \alpd_{\ms} \dimA_{\ms} 
\label{summmalpdimAmm}
\end{EQA}
with \( \alpd_{\ms} = \CONST (\dimA_{\mmmin}/\dimA_{\ms})^{1+a} \).

%The next statement extends the bound of Lemma~\ref{Lxiv2LD}.
%\begin{lemma}
%\label{Lxivxiv1BBBB1}
For any random vector \( \xiv \) with \( \Var(\xiv) = \BB \) and \( \dimA = \tr(\BB) \)
and any random event \( \AA \),
it holds % some absolute constant \( \CONST \) and 
\begin{EQA}
	\E \Bigl[ 
		\dimA^{-1} \| \xiv \|^{2} 
		\Ind(\AA)
	\Bigr]
	& \leq &
	\bigl\{ 1 + \dimA^{-2} \Var(\| \xiv \|^{2}) \bigr\}^{1/2} \P^{1/2}(\AA) . 
\label{dimAdimA1}
\end{EQA}
Indeed, the Cauchy-Schwartz inequality implies
\begin{EQA}
	\E \Bigl\{
		\dimA^{-1} \| \xiv \|^{2} 
		\Ind(\AA)
	\Bigr\}
	& \leq &
	\E^{1/2} \bigl\{ \dimA^{-1} \| \xiv \|^{2} \bigr\}^{2} \,
	\P^{1/2} \bigl( \AA \bigr)
	\\
	& = &
	\bigl\{ 1 + \dimA^{-2} \Var(\| \xiv \|^{2}) \bigr\}^{1/2} \P^{1/2}(\AA) .
\label{dimAm1xi2E12P12}
\end{EQA}
Moreover, in the Gaussian case \( \xiv \sim \ND(0,\BB) \) with \( \| \BB \|_{\oper} \leq 1 \),
it holds \( \Var(\| \xiv \|^{2}) \leq 2 \dimA \). 
If \( \dimA \) is large then 
\( \Var(\| \xiv \|^{2}) / \dimA^{2} \) is small. 
In general \( \Var(\| \xiv \|^{2}) / \dimA^{2} \leq 2 \).

Result \eqref{dimAdimA1} and the choice \( \alp_{\mm} = \sqrt{3} \dimA_{\mm}^{-1-a} \) %from \eqref{alpmmdimAa}
allow to specify an upper bound on \( \xx_{\mm} \).
Namely, the choice \( \xx_{\mm} = \CONST \log(\dimA_{\mm}) \) ensures the propagation condition
\eqref{riskrmmPL}.
To see this, fix \( \mm \) and \( \mc \geq \mm \).
Let 
\begin{EQA}
	A'_{\mm}(\xx)
	& \eqdef &
	\Ind\Bigl( 
			\max_{\mc \in \MMu(\mm)} 
			\bigl\{ 
				\| \xiv_{\mc,\mm} \| - \sqrt{\dimA_{\mc,\mm}} 
				- \sqrt{2 \supA_{\mc,\mm} \, \{ \xx + \log(|\MM|) \}}
			\bigr\} > 0 
		\Bigr) 
\label{Ammxxmcmm122}
\end{EQA}
The arguments after Lemma~\ref{Lxiv2LD} with \( \xx_{\mmprec} = 2 (1+a) \log (\dimA_{\mm}) \)
and \eqref{dimAdimA1} imply
\begin{EQA}
	\E \Bigl[ 
		\dimA_{\mm}^{-1} \| \xiv_{\mm} \|^{2} 
		\Ind\{ A'_{\mmprec}(\xx_{\mmprec}) \}
	\Bigr]
	& \leq &
	\sqrt{3} \ex^{- (1+a) \log (\dimA_{\mm})}
	=
	\sqrt{3} \dimA_{\mm}^{- 1 - a}
\label{Epmm1cxim231a}
\end{EQA}
and by \eqref{zqmmmdPL}
\begin{EQA}
	\ZZ_{\mm,\md}
	& \leq &
	\sqrt{\dimA_{\mm,\md}} 
		+ \sqrt{2 \supA_{\mm,\md} \, \{ (1+a) \log (\dimA_{\md+1}) + \log(|\MM|) \}} \, . 
\label{ZZmmmdlogdimAmd}
\end{EQA}
This implies the upper bound \eqref{Zpaymentpower} on the payment for adaptation 
\( \ZZs_{\ms} \).
}{ %annals 
}

%\end{proof}

\subsection{Proof of Theorem~\ref{TGaussbootB}}
Any statement on the use of bootstrap-tuned parameters faces the same fundamental problem:
the bootstrap distribution is random and depends on the underlying sample.
When we use such values for the original procedure, we have to account for this dependence.
The statement of Theorem~\ref{TGaussbootB} is even more involved due to the presmoothing 
step and multiplicity correction \eqref{Pximdmmububst}.
The proof will be split into a couple of steps.
First we evaluate the effect of the presmoothing bias and variance and reduce the study to 
an artificial situation where one uses the errors \( \eps_{i} \) for resampling  in place of
the residuals \( \Yr_{i} \).
Then we compare \( \PPsi \) and \( \PPsib \) using the Pinsker inequality.

Below we write \( \Psi \) in place of \( \Psi_{\mmmax} \), where \( \mmmax \) is the largest
model in the collection. 
This does not conflict with our general setup, it is implicitly assumed that the largest
model coincides with the original one.
By \( \dimp \) we denote the corresponding parameter dimension, that is, \( \Psi \) 
is a \( \dimp \times n \) matrix.
Further, the feature matrix \( \Psi_{\mm} \) can be written as the product 
\( \Psi_{\mm} = \Pi_{\mm} \Psi \), where \( \Pi_{\mm} \) is the projector on the subspace 
of the feature space spanned by the features from the model \( \mm \):
\( \Pi_{\mm} = \Psi_{\mm}^{\T} \bigl( \Psi_{\mm} \Psi_{\mm}^{\T} \bigr)^{-1} \Psi_{\mm} \).
This allows to represent each estimator \( \tilde{\tarv}_{\mm} \) in the form
\begin{EQA}
	\tilde{\tarv}_{\mm}
	&=&
	\QL \tilde{\thetav}_{\mm}
	=
	\QL \Gam_{\mm} \Yv
	=
	\QL \bigl( \Psi_{\mm} \Psi_{\mm}^{\T} \bigr)^{-1} \Psi_{\mm} \Yv
	=
	\Sam_{\mm} \Psi \Yv
\label{ttrmQQSmP}
%\end{EQA}
%with 
%\begin{EQA}
	\\
	\Sam_{\mm}
	& \eqdef &
	\QL \bigl( \Psi_{\mm} \Psi_{\mm}^{\T} \bigr)^{-1} \Pi_{\mm} \, .
\label{SmQPmPmTm1}
\end{EQA}
This implies the following representation of the stochastic components \( \xiv_{\mm,\md} \):
\begin{EQA}
	\xiv_{\mm,\md}
	=
	\Sam_{\mm,\md} \Psi \epsv
	=
	\Sam_{\mm,\md} \score ,
	&\quad &
	\Sam_{\mm,\md}
	\eqdef 
	\Sam_{\mm} - \Sam_{\md} \, ,
\label{xivmmmdSmSmds}
\end{EQA}
where \( \score = \Psi \epsv \).
One can say that each stochastic vector \( \xiv_{\mm,\md} \) is a linear function
of the vector \( \score \).
A similar representation holds true in the bootstrap world:
\begin{EQA}
	\xivb_{\mm,\md}
	=
	\Sam_{\mm,\md} \Psi \diag(\Yvr) \Wb
	=
	\Sam_{\mm,\md} \scoreb ,
	&\quad &
	\scoreb
	\eqdef 
	\Psi \diag(\Yvr) \Wb .
\label{xivmmmdlinrp}
\end{EQA}
Here the original errors \(\epsv \) are replaced by their bootstrap surrogates 
\( \epsvb = \diag(\Yvr) \Wb\).
Therefore, it suffices to compare the distribution of \( \score = \Psi \epsv \) 
with the conditional distribution of \( \scoreb = \Psi \diag(\Yvr) \Wb \) given \( \Yv \).
Then the results will be automatically extended to any deterministic mapping of these 
two vectors.

Normality of the errors \( \eps_{i} \sim \ND(0,\sigma_{i}^{2}) \) implies 
that \( \score = \Psi \epsv \) is also normal zero mean:
\begin{EQA}
	\score
	& \sim &
%	\ND(0, \Psi \Sigma \Psi^{\T})
%	=
	\ND(0,\Varxi),
	\qquad
	\Varxi \eqdef \Psi \Sigma \Psi^{\T},
	\qquad
	\Sigma = \Var(\epsv)
	=
	\diag\bigl( \sigma_{1}^{2},\ldots,\sigma_{n}^{2} \bigr) .
\label{epsvND0Si}
\end{EQA} 
Similarly we can use standard normality of the bootstrap weights  
\( \wb_{i} \).
Given the data \( \Yv \), the vector \( \scoreb \) is conditionally normal zero mean 
with the conditional variance
\begin{EQA}
	\Varxib
	\eqdef 
	\Varb(\scoreb)
	& = &
	\Psi \diag\bigl( \Yr_{1}^{2}, \ldots, \Yr_{n}^{2} \bigr) \Psi^{\T}
	=
	\Psi \diag\bigl( \Yvr \cdot \Yvr \bigr) \Psi^{\T} .
\label{Varbscob1n}
\end{EQA}
Therefore, the problem is reduced to comparing two \( \dimp \)-dimensional Gaussian 
distributions with different covariance matrices. 
Equivalently, we have to bound the value \( \errSi = \sqrt{\tr(\BBB^{2})} \) for a random 
\( \dimp \times \dimp \) matrix \( \BBB \) given by
%\begin{EQA}
%	\DeltaPsi
%	& \eqdef & 
%	\bigl\| \Varxi^{-1/2} \bigl( \Varxib - \Varxi \bigr) \Varxi^{-1/2} \bigr\|_{\oper} \, .
%\label{Varxi12b12}
%\end{EQA}
\begin{EQA}
%	\errSi^{2}
%	& \eqdef & 
%	\tr(\BBB^{2}),
%	\qquad
	\BBB
	& \eqdef &
	\Varxi^{-1/2} \bigl( \Varxib - \Varxi \bigr) \Varxi^{-1/2}  \, .
\label{Varxi12b12}
\end{EQA}
%in the operator norm.
%
Define a \( \dimp \times n \) matrix \( \UV = \Varxi^{-1/2} \Psi \Sigma^{1/2} \) so that
\( %begin{EQA}
	\UV \UV^{\T} 
	=
	\Id_{\dimp} .
\label{UVUVTIdp}
\) %end{EQA} 
We will use the decomposition
\begin{EQA}
	\Sigma^{-1/2} \Yvr
	&=&
	\Sigma^{-1/2} (\Yv - \Pi \Yv)
	=
	\Sigma^{-1/2} (\epsv - \Pi \epsv) + \Sigma^{-1/2} (\fvs - \Pi \fvs)
	=
	\etav + \Bias
\label{Yvrdecomp}
\end{EQA}
with
\begin{EQ}[rcl]
	\etav
	\eqdef
	\Sigma^{-1/2} (\epsv - \Pi \epsv),
	&\quad &
	\Bias
	\eqdef
	\Sigma^{-1/2} \bigl( \fvs - \Pi \fvs \bigr) .
\label{tepsSi12PieB}
\end{EQ}
With %\( \gaussv = \Sigma^{-1/2} \epsv \),
the matrix \( \BBB \) can now be represented as 
\begin{EQA}[rclccl]
\label{BBBdefVm12VVV}
	\BBB 
	&=&
	\UV \diag\bigl\{ 
		(\etav + \Bias) \cdot (\etav + \Bias) - \Id_{n} 
	\bigr\} \UV^{\T}
	&&&
	\\
	&=&
	\UV \diag\bigl\{ 
		(\etav + \Bias) \cdot (\etav + \Bias) - \etav \cdot \etav 
	\bigr\} \UV^{\T} 
	& \qquad
	& \eqdef & \BBB_{1}
	\\
	&& \quad
	+ \, \UV \diag\bigl\{ 
		\etav \cdot \etav - \E (\etav \cdot \etav)
	\bigr\} \UV^{\T}
	& \qquad
	& \eqdef & \BBB_{2}
	\\
	&& \quad
	+ \, \UV \diag\bigl\{ 
		\E (\etav \cdot \etav) - \Id_{n} 
	\bigr\} \UV^{\T}
	& \qquad
	& \eqdef & \BBB_{3}
%	\\
%	&=&
%	\BBB_{1} + \BBB_{2} + \BBB_{3} .
\end{EQA}
The first term \( \BBB_{1} \) in this decomposition expresses the impact of the bias \( \Bias \) 
remaining after presmoothing,
the last two terms \( \BBB_{2} \) and \( \BBB_{3} \) measure the change of the noise covariance due to presmoothing.
The triangle inequality in the Frobenius norm \( \| \BBB \|_{\Fr} \eqdef \sqrt{\tr(\BBB^{2})} \)
and bounds from 
\ifAoS{ % full
Propositions~\ref{TUVeps1UV}, \ref{TUVeps2UV}, and \ref{TUVeps22UV} }
{ %annals
Propositions~D.2, D.4, D6}
with \( \UV \UV^{\T} = \Id_{\dimp} \) and \( \dimw = \dimwf = \tr(\UV \UV^{\T}) = \dimp \)
imply on a random set \( \Omega_{2}(\xx) = \Omega_{12}(\xx) \cup \Omega_{22}(\xx) \) with 
\( \P\bigl( \Omega_{2}(\xx) \bigr) \geq 1 - 2 \ex^{-\xx} \)
\begin{EQA}
	\| \BBB \|_{\Fr}
	& \leq &
	\| \BBB_{1} \|_{\Fr} + \| \BBB_{2} \|_{\Fr} + \| \BBB_{3} \|_{\Fr}
	\\
	& \leq &
	\errSi_{1}(\xx) + \errSi_{2}(\xx) + \errSi_{3}(\xx)
	\\
	&=&
	2 \sqrt{\dPsi^{2} \, \dimp \, (\xx + \log(n))}
	+ \sqrt{\supepsi^{2} \, \dimp} 
	+ \sqrt{\| \Bias \|_{\infty}^{4} \, \dimp}
%	+ \sqrt{\dPsi^{4} \| \Bias \|^{4} \dimp}
	+ 4 \, \dPsi^{2} \, \| \Bias \| \, \bigl( 1 + \sqrt{\xx} \bigr) .
\label{BBBFr124}
\end{EQA}
This proves \eqref{PPsiPPbTV} in view of 
\ifAoS{ % full
Pinsker's Lemma~\ref{KullbTVd}}
{ % annals
Lemma~E.1 in \suppSW}
with \( \bvs = \bvb = 0 \).

\subsection{Proof of Theorem~\ref{TGaussbootB2}}
The result of Theorem~\ref{TGaussbootB} justifies 
the bootstrap-phenomenon, namely it explains why the known bootstrap distribution can be used 
as a proxy for the unknown error distribution.
However, it cannot be applied directly to \eqref{multcorrBmix} because the quantities 
\( \zqb_{\mm,\md}(\xx) \) and \( \qqbt_{\md} \) are random and depend on the original data. 
This especially concerns the multiplicity correction \( \qqbt_{\md} \) which is based 
on the joint distribution of the vectors \( \xivb_{\mm,\md} \) from 
\eqref{xivbmmmddef} and is defined in \eqref{Pximdmmububst}.
The latter distribution is a random measure in the bootstrap world which is normal conditioned 
on the original sample. 
To cope with the problem of this cross-dependence, we apply the statement of Theorem~\ref{Tmultcorrb} in the Appendix.
The underlying idea is to use geometric arguments to sandwich the random probability in 
\eqref{Pximdmmububst} in two deterministic probabilities. 
Then the error of bootstrap approximation can again be bounded by using the Pinsker inequality.
The statement of Theorem~\ref{TGaussbootB2} can be derived from 
Theorem~\ref{Tmultcorrb} if an operator norm bound \( \| \BBB \|_{\oper} \) is available.
Note that Theorem~\ref{TGaussbootB} only requires a bound for the Frobenius norm. 
\ifAoS{ % full
By Proposition~\ref{TboundBBBinoper},}
{ % annals
By Proposition~D.9 in \suppSW,}
it holds with \( \tdn = \dPsi \), \( \xxn = \xx + \log(n) \), 
and \( \xx_{\dimp} = \xx + 2 \log(\dimp) \)
\begin{EQA}
	\| \BBB \|_{\oper}
	& \leq &
	\errSi_{\oper}(\xx),
	\\
	\errSi_{\oper}(\xx)
	& \eqdef &
	\| \Bias \|_{\infty}^{2} + \dPsi^{2} \| \Bias \| \sqrt{2 \xx} 
	+ 2 \dPsi \xx_{\dimp}^{1/2} + 2 \dPsi^{2} \xx_{\dimp}
	+ 2 \dPsip \xxn + \dPsip^{2} \xxn .
\label{BBBoperrSopxre}
\end{EQA}
The result of the theorem follows now by Theorem~\ref{Tmultcorrb}.

\subsection{Proof of Theorem~\ref{TdimAbtdimA}}
For a fixed pair \( \mm > \md \) from \( \MM \), consider 
\( \dimAb_{\mm,\md} = \Eb \| \xivb_{\mm,\md} \|^{2} \) and 
\( \dimA_{\mm,\md} = \E \| \xiv_{\mm,\md} \|^{2} \).
As \( \diag(\Yvr) \) and \( \Sigma \) are diagonal matrices, the definitions \eqref{xivbmmmddef} 
and \eqref{tepsSi12PieB} imply
\begin{EQA}
	\xivb_{\mm,\md} 
	&=&
	\Tam_{\mm,\md} \diag(\Yvr) \Wb  
	=
	\Tam_{\mm,\md} \Sigma^{1/2} \Sigma^{-1/2} \diag(\Yvr) \Wb 
	\\
	&=&
	\UV_{\mm,\md} \, \diag(\etav + \Bias) \Wb,
\label{xivbmmmdb1}
\label{xivbmmmdUVeB}
\end{EQA}
where \( \UV_{\mm,\md} \eqdef \Tam_{\mm,\md} \Sigma^{1/2} \).
It holds for \( \dimAb_{\mm,\md} \)
\begin{EQA}
	\dimAb_{\mm,\md}
	=
	\Eb \bigl\| \xivb_{\mm,\md} \bigr\|^{2}
	&=&
	\tr \Bigl( 
		\UV_{\mm,\md} \, \diag \bigl\{ (\etav + \Bias) \cdot (\etav + \Bias) \bigr\} \,
		\UV_{\mm,\md}^{\T} 
	\Bigr)
\label{dimAbmmmdUVT}
\end{EQA}
while \( \xiv_{\mm,\md} = \Tam_{\mm,\md} \Sigma^{1/2} \Sigma^{-1/2} \epsv \) and
\begin{EQA}
	\dimA_{\mm,\md}
	&=&
	\E \bigl\| \xiv_{\mm,\md} \bigr\|^{2}
	=
	\tr \bigl( \UV_{\mm,\md} \, \UV_{\mm,\md}^{\T} \bigr) .
\label{dimAmmmdEtrUVmd}
\end{EQA}
As we are interested in the ratio \( \dimAb_{\mm,\md}/ \dimA_{\mm,\md} \), one can assume
without loss of generality that \( \| \UV_{\mm,\md} \, \UV_{\mm,\md}^{\T} \|_{\oper} = 1 \)
and \( \dimA_{\mm,\md} \geq 1 \).
Now we again apply the decomposition \eqref{BBBdefVm12VVV}.
The bounds 
\ifAoS{ % full
\eqref{trBBB11} of Proposition~\ref{TUVeps1UV}, 
\eqref{trBBB21} of Proposition~\ref{TUVeps2UV}, 
and \eqref{trBBB32} of Proposition~\ref{TUVeps22UV} }
{ % annals
from Propositions~D.2, D.4, and D.6}
imply on a set \( \Omega_{\mm,\md}(\xx) \)
with \( \P\bigl( \Omega_{\mm,\md}(\xx) \bigr) \geq 1 - 3 \ex^{-\xx} \)
\begin{EQA}
	\biggl| \frac{\dimAb_{\mm,\md}}{\dimA_{\mm,\md}} - 1 \biggr|
	& \leq &
	\| \Bias \|_{\infty}^{2} + 4 \, \xx^{1/2} \, \tdn^{2} \, \| \Bias \|
	+ 4 \xx^{1/2} \, \tdn + 4 \, \xx \, \tdn^{2} 
	+ \supepsi \,  .
\label{pbmdpmdm1}
\end{EQA}
The choice of \( \xx = \xx_{\MM} = \xx + 2 \log\bigl( \bigl| \MM \bigr| \bigr) \) ensures 
a uniform bound for all pairs \( \mm > \md \) from \( \MM \).

%%%%%%%%%%%%%%%%%%%%%%%%%%%%%%%%%%%%%%%%%%%%%%%%%%%%%%%%%%%%%%%%
% !TEX root = script2014.tex

\section{Random multiplicity correction}
\label{Smultboot}
Suppose that \( \VPb \) is a random positive symmetric \( \dimp \times \dimp \) matrix
close to a deterministic matrix \( \VP \).
Below we use the operator norm for quantifying the difference between \( \VP \) and \( \VPb \):
% that
namely let with probability one 
\begin{EQA}
	\| \VP^{-1/2} \, \VPb \, \VP^{-1/2} - \Id_{\dimp} \|_{\oper}
	& \leq &
	\errV .
\label{VPbVPerrVop}
\end{EQA}
In what follows, \( \P  = \ND(0,\VP) \) is the normal measure on \( \R^{\dimp} \) with mean zero and covariance 
\( \VP \).
%\( \xiv \) is a Gaussian zero mean vector with \( \xiv \sim \ND(0,\VP) \) 
Similarly \( \Pb \) is a random measure on \( \R^{\dimp} \) which is conditionally on \( \VPb \) normal  with \( \Pb = \ND(0,\VPb) \).
%\( \xivb \) is conditionally on \( \VP \) normal vector with 
%\( \xivb \cond \VPb \sim \ND(0,\VPb) \).
Suppose that for each \( \mm \) from a given set \( \MM \) a linear mapping 
\( T_{\mm} \colon \R^{\dimp} \to \R^{\dimp_{\mm}} \) is fixed.
Given \( \xx \), define for each \( \mm \in \MM \) the corresponding tail function 
\( \zq_{\mm}(\xx) \) by
\begin{EQA}
	\P\bigl\{ \uv \colon \| T_{\mm} \uv \| \geq \zq_{\mm}(\xx) \bigr\}
	& = & 
	\ex^{- \xx} .
\label{PTmmxixqxx}
\end{EQA}
Also define a set \( \AA(\xx) \) as 
\begin{EQA}
	\AA(\xx) 
	& \eqdef &
	\biggl\{ \uv \colon
		\bigcap_{\mm \in \MM} \bigl\{ \| T_{\mm} \uv \| \leq \zq_{\mm}(\xx) \bigr\}  
	\biggr\}\, .
\label{AAxxdef}
\end{EQA}
Similarly define \( \zqb_{\mm}(\xx) \) by \eqref{PTmmxixqxx} with \( \Pb \) in place of 
\( \P \), \( \mm \in \MM \), and \( \AAb(\xx) \).
Note that all these objects are random because \( \Pb \) is random.
Finally, let \( \xxb_{\alp} \) be the random quantity providing 
\begin{EQA}
	\Pb\bigl( \AAb(\xxb_{\alp}) \bigr)
	&=&
	1 - \alp .
\label{PbAAbxxba1}
\end{EQA}
Below we try to address the question whether this random multiplicity correction 
based on \eqref{PbAAbxxba1} does a good job under \( \P \).
This question leads to analysis of value \( \P\bigl( \AAb(\xxb_{\alp}) \bigr) \):
the goal is in evaluating the difference 
\begin{EQA}[c]
	\P\bigl( \AAb(\xxb_{\alp}) \bigr) - (1 - \alp) .
\label{PxivAAbxxbalp}
\end{EQA}
The analysis is non-trivial because \( \AAb(\xx) \) and \( \xxb_{\alp} \) are random.

\begin{theorem}
\label{Tmultcorrb}
Let a random matrix \( \VPb \) satisfy \eqref{VPbVPerrVop} for 
%two deterministic matrix \( \VP^{-} \) and \( \VP^{+} \) such that \eqref{VPmVPbVPb} holds. 
a deterministic matrix \( \VP \) and \( \errV < 1/2 \).
Then it holds 
%with \( \errSi = \sqrt{\dimp} \, \errV \)
%for any Gaussian measure \( \P = \ND(0,\VP) \) with \( \VP^{-} \leq \VP \leq \VP^{+} \)
\begin{EQA}
	\bigl| \P\bigl( \AAb(\xxb_{\alp}) \bigr) - 1 + \alp \bigr|
	& \leq &
	\sqrt{\dimp} \, \errV .
\label{PAbxba1a}
\end{EQA}
\end{theorem}

\ifims{ %full
\begin{proof}
The key property of \( \Pb = \ND(0,\VPb) \) is that the random matrix \( \VPb \) 
concentrates around some deterministic matrix.
Below we use this property in the bracketing form:
\begin{EQA}[c]
	\VP^{-}
	\leq 
	\VPb 
	\leq 
	\VP^{+} 
\label{VPbVPerrV}
\end{EQA}
with 
\begin{EQA}
	\VP^{-}
	& \eqdef &
	(1 - \errV) \VP,
	\quad
	\VP^{+}
	\eqdef 
	(1 + \errV) \VP ,
	\quad
	\VP^{+} - \VP^{-} 
	= 
	2 \errV \VP. 
\label{VPbVPerrVopa}
\end{EQA}
In other words, 
the random matrix \( \VPb \) can be sandwiched in two deterministic matrices 
\( \VP^{-} \) and \( \VP^{+} \).
%Moreover, these two matrices are close to each other in the sense 
%\begin{EQA}
%	\distS(\VP^{+},\VP^{-}) 
%	&\leq &
%	2 \sqrt{\dimp} \, \errV ,
%\label{VPmVPbVPb}
%\end{EQA}
%where \( \errSi = \errV \sqrt{\dimp} \) and the distance \( \distS(\cdot,\cdot) \) is defined by
%\begin{EQA}
%	\distS(\VP,\VPc)
%	& \eqdef &
%	\tr^{1/2}\Bigl\{ \bigl( \VP^{-1/2} \, \VPc \, \VP^{-1/2} - \Id_{\dimp} \bigr)^{2} \Bigr\}
%	=
%	\bigl\| \VP^{-1/2} \, \VPc \, \VP^{-1/2} - \Id_{\dimp} \bigr\|_{Fr} \,  .
%\label{distSdefVPc}
%\end{EQA}
%For instance, for \( \VP^{-} \) and \( \VP^{+} \) from \eqref{VPbVPerrVopa},
%the bound \eqref{VPmVPbVPb} holds with \( \errSi \leq \errV \sqrt{\dimp} \).
%
For the proof of \eqref{PAbxba1a} we use the following well known property 
of the Gaussian distribution.

\begin{lemma}
\label{P1P2V1V2}
Let \( \P_{1} \sim \ND(0,\VP_{1}) \) and \( \P_{2} \sim \ND(0,\VP_{2}) \) with 
\( \VP_{1} \leq \VP_{2} \). 
Then for any centrally symmetric star-shaped set \( \AA \), it holds
\begin{EQA}
	\P_{1}(\AA)
	& \geq &
	\P_{2}(\AA) .
\label{P1AgP2A}
\end{EQA}
\end{lemma}

\begin{proof}
The statement is trivial in the univariate case, the general case is obtained by integration
over \( \AA \) in polar coordinates. 
\end{proof}

Introduce two Gaussian measures \( \P^{-} = \ND(0,\VP^{-}) \) and 
\( \P^{+} = \ND(0,\VP^{+}) \); see \eqref{VPbVPerrVopa}.
Let \( \zq^{-}_{\mm}(\xx) \) and \( \zq^{+}_{\mm}(\xx) \) be the corresponding tail functions, and
\( \AA^{-}(\xx) \) and \( \AA^{+}(\xx) \) - the corresponding sets.
The identities \eqref{VPbVPerrVopa} yield for each \( \xx \) the relation
\begin{EQA}
	\P^{+}\bigl( \AA^{+}(\xx) \bigr)
	&=&
	\P^{-}\bigl( \AA^{-}(\xx) \bigr) .
\label{AAumagmAAd}
\end{EQA}
%with
%\begin{EQA}
%	1 + \magnpm
%	& \eqdef &
%	\biggl( \frac{1 + \errV}{1 - \errV} \biggr)^{1/2} .
%\label{1magndef1pe1me}
%\end{EQA}
Lemma~\ref{P1P2V1V2} implies by \eqref{VPbVPerrVopa} for any \( \xx \)
\begin{EQA}
	\P^{+}(\AA(\xx))
	& \leq &
	\Pb(\AA(\xx))
	\leq 
	\P^{-}(\AA(\xx)) .
\label{PbAPbAPmA}
\end{EQA}

The key step of the proof is given by the next lemma where we
sandwich the random set \( \AAb(\xxb) \) in two specially constructed deterministic sets.

\begin{lemma}
\label{Lxxuxxm}
Define the deterministic values \( \xx^{-}_{\alp} \) and \( \xx^{+}_{\alp} \) by the equations
\begin{EQ}[rcl]
	\P^{+}\bigl( \AA^{-}(\xx^{+}_{\alp}) \bigr)
	&=&
	1 - \alp,
	\\
	\P^{-}\bigl( \AA^{+}(\xx^{-}_{\alp}) \bigr)
	&=&
	1 - \alp .
\label{P1alPm1alp}
\end{EQ}
Then 
\begin{EQA}[ccccc]
	\xx^{-}_{\alp}
	& \leq &
	\xxb_{\alp}
	& \leq & 
	\xx^{+}_{\alp}
	\\
	\AA^{-}(\xx^{-}_{\alp})
	& \subseteq &
	\AAb(\xxb_{\alp})
	& \subseteq &
	\AA^{+}(\xx^{+}_{\alp}) .
\label{AAmxxuAAbAAm}
\end{EQA}
\end{lemma}

\begin{proof}
By Lemma~\ref{P1P2V1V2} the following inequalities and inclusions hold true for any \( \xx \):
\begin{EQA}[rcccl]
	\zq^{-}_{\mm}(\xx)
	& \leq &
	\zqb_{\mm}(\xx)
	& \leq & 
	\zq^{+}_{\mm}(\xx), 
	\\
	\AA^{-}(\xx)
	& \subseteq &
	\AAb(\xx)
	& \subseteq &
	\AA^{+}(\xx) .
\label{zqmzqbAAmAAbAAu}
\end{EQA}
Now by definition \eqref{P1alPm1alp} in view of \eqref{PbAPbAPmA} and \eqref{zqmzqbAAmAAbAAu}
\begin{EQA}
	\Pb\bigl( \AAb(\xx^{+}_{\alp}) \bigr)
	& \geq &
	\P^{+}\bigl( \AAb(\xx^{+}_{\alp}) \bigr)
	\geq 
	\P^{+}\bigl( \AA^{-}(\xx^{+}_{\alp}) \bigr)
	=
	1 - \alp,
	\\
	\Pb\bigl( \AAb(\xx^{-}_{\alp}) \bigr)
	& \leq &
	\P^{-}\bigl( \AAb(\xx^{-}_{\alp}) \bigr)
	\leq 
	\P^{-}\bigl( \AA^{+}(\xx^{-}_{\alp}) \bigr)
	=
	1 - \alp .
\label{PbPuPm1al1al}
\end{EQA}
This yields by monotonicity of \( \Pb(\AAb(\xx)) \) in \( \xx \) that 
\( \xxb_{\alp} \) from \eqref{PbAAbxxba1} belongs to the interval 
\( [\xx^{-}_{\alp},\xx^{+}_{\alp}] \) and
\begin{EQA}
	\AA^{-}(\xx^{-}_{\alp})
	& \subseteq &
	\AAb(\xx^{-}_{\alp})
	\subseteq
	\AAb(\xxb_{\alp})
	\subseteq
	\AAb(\xx^{+}_{\alp})
	\subseteq
	\AA^{+}(\xx^{+}_{\alp}) .
\label{AAmAAuxxi}
\end{EQA}
This implies the result.
\end{proof}
%
%It remains to bound the difference 
Now we are prepared to finalize the proof.
The relations \eqref{AAmxxuAAbAAm} and \eqref{AAumagmAAd} imply
\begin{EQA}
%	\P^{+}\bigl( \AAb(\xxb_{\alp}) \bigr)
%	& \geq &
%	\P^{+}\bigl( \AA^{-}(\xx^{-}_{\alp}) \bigr) ,
%%	= 
%%	1 - \alp ,
%	\\
	\P^{+}\bigl( \AAb(\xxb_{\alp}) \bigr)
	& \leq &
	\P^{+}\bigl( \AA^{+}(\xx^{+}_{\alp}) \bigr) 
	= 
	\P^{-}\bigl( \AA^{-}(\xx^{+}_{\alp}) \bigr) .
%	\leq 
%	\magnpm + \P^{+}\bigl( \AA^{-}(\xx^{+}_{\alp}) \bigr)
%	\leq 
%	\magnpm + 1 - \alp .
\label{PuAAbxba1aPu1mae}
\end{EQA}
Furthermore, it holds by Pinsker' inequality Corollary~\ref{KullbTVoper} in view of  \eqref{VPbVPerrVop} and \eqref{P1alPm1alp}
\begin{EQA}
	\P^{-}\bigl( \AA^{-}(\xx^{+}_{\alp}) \bigr)
	& \leq & 
	\P^{+}\bigl( \AA^{-}(\xx^{+}_{\alp}) \bigr) + \sqrt{\dimp} \, \errV
	\leq 
	1 - \alp + \sqrt{\dimp} \, \errV .
\label{PuAuaPm1ae}
\end{EQA}
Similarly
\begin{EQA}
	\P^{-}\bigl( \AAb(\xxb_{\alp}) \bigr)
	& \geq &
	\P^{-}\bigl( \AA^{-}(\xx^{-}_{\alp}) \bigr)
	= 
%	\P^{+}\bigl( \AA^{-}(\xx^{-}_{\alp}) \bigr) - \sqrt{\dimp} \, \errV
	\P^{+}\bigl( \AA^{+}(\xx^{-}_{\alp}) \bigr)
%	=
%	1 - \alp - \sqrt{\dimp} \, \errV ,
	\\
%	\P^{-}\bigl( \AAb(\xxb_{\alp}) \bigr)
	& \geq &
%	\geq 
	\P^{-}\bigl( \AA^{+}(\xx^{-}_{\alp}) \bigr) - \sqrt{\dimp} \, \errV
	=
	1 - \alp - \sqrt{\dimp} \, \errV.
\label{PuAAbxba1aPu1mae}
\end{EQA}
This implies \eqref{PAbxba1a} for the measure \( \P \) .
\end{proof}
}{ % annals
The proof is given in \suppSW.}

\ifbook{}{
Now we consider a more general situation.
Let \( T_{\mm} \) be a family of test statistics in the real world,
and \( \Tb_{\mm} \)
}

%%%%%%%%%%%%%%%%%%%%%%%%%%%%%%%%%%%%%%%%%%%%%%%%%%%%%%%%%%%%%%%%

%%%%%%%%%%%%%%%%%%%%%%%%%%%%%%%%%%%%%%%%%%%%%%%%%%%%%%%%%%%%%%%%
\ifbook{}{
\input bootone
}
%%%%%%%%%%%%%%%%%%%%%%%%%%%%%%%%%%%%%%%%%%%%%%%%%%%%%%%%%%%%%%%%

\ifims{ %full

%%%%%%%%%%%%%%%%%%%%%%%%%%%%%%%%%%%%%%%%%%%%%%%%%%%%%%%%%%%%%%%%

\section{Deviation bounds for Gaussian law}
\label{SGaussdi}

This section collects some simple but useful facts about the properties of the multivariate standard 
normal distribution.
Many similar results can be found in the literature, we present the proofs to keep the presentation 
self-contained. 
%
%\subsection{Deviation bounds for Gaussian quadratic forms}
%This section collects some deviation bounds on the norm or quadratic form 
%of a standard normal vector. 
Everywhere in this section \( \gaussv \) means a standard normal vector in \( \R^{\dimp} \).
%

%\ifbook{}{
\begin{lemma}
\label{LGaussmom}
Let \( \mu \in (0,1) \).
Then for any vector \( \lambdav \in \R^{\dimp} \) with
\( \| \lambdav \|^{2} \le \dimp \) and any \( \rr > 0 \)
\begin{EQA}[c]
    \log \E \bigl\{
        \exp (\lambdav^{\T} \gaussv) \Ind\bigl( \| \gaussv \| > \rr \bigr)
    \bigr\}
    \le
    - \frac{1 - \mu}{2} \rr^{2} + \frac{1}{2 \mu} \| \lambdav \|^{2}
    + \frac{\dimp}{2} \log(\mu^{-1}) .
    \qquad
\label{lamgamrr}
\end{EQA}
Moreover, if \( \rr^{2} \ge 6 \dimp + 4 \xx \), then 
\begin{EQA}[c]
    \E \Bigl\{
        \exp (\lambdav^{\T} \gaussv) \Ind\bigl( \| \gaussv \| \le \rr \bigr)
    \Bigr\}
    \ge
    \ex^{ \| \lambdav \|^{2}/2} \bigl( 1 - \ex^{-\xx} \bigr) .
\label{exIdgamrr}
\end{EQA}
\end{lemma}

%\ifbook{}{
\begin{proof}
We use that for \( \mu < 1 \)
\begin{EQA}[c]
    \E \bigl\{
        \exp(\lambdav^{\T} \gaussv) \Ind\bigl( \| \gaussv \| > \rr \bigr)
    \bigr\}
    \le
    \ex^{- (1 - \mu) \rr^{2}/2}
        \E \exp \bigl\{ \lambdav^{\T} \gaussv + (1 - \mu) \| \gaussv \|^{2}/2 \bigr\} .
\label{explamgamrr}
\end{EQA}
It holds
\begin{EQA}
    \E \exp \bigl\{ \lambdav^{\T} \gaussv + (1 - \mu) \| \gaussv \|^{2}/2 \bigr\}
    &=&
    (2 \pi)^{-\dimp/2} \int
        \exp \bigl\{ \lambdav^{\T} \gaussv - \mu \| \gaussv \|^{2}/2 \bigr\} d \gaussv
	\\
    &=&
    \mu^{-\dimp/2} \exp\bigl( \mu^{-1} \| \lambdav \|^{2}/2 \bigr)
\label{exp1mummu}
\end{EQA}
and \eqref{lamgamrr} follows.

Now we apply this result with \( \mu = 1/2 \).
In view of
\( \E \exp(\lambdav^{\T} \gaussv) = \ex^{\| \lambdav \|^{2}/2} \),
\( \rr^{2} \geq 6 \dimp + 4 \xx \),
and \( 2 + \log(2) < 3 \), it follows
for \( \| \lambdav \|^{2} \le \dimp \)
\begin{EQA}
    && \nquad
    \ex^{- \| \lambdav \|^{2}/2}
    \E \bigl\{
        \exp (\lambdav^{\T} \gaussv) \Ind\bigl( \| \gaussv \| \le \rr \bigr)
    \bigr\}
    \\
    & \ge &
    1 - \exp\bigl( - \rr^{2}/4 + \dimp + (\dimp/2) \log(2) \bigr)
    \ge
    1 - \exp(- \xx)
\label{exIdgamrrd}
\end{EQA}
which implies \eqref{exIdgamrr}.
\end{proof}
%}

\begin{lemma}
For any \( \uv \in \R^{\dimp} \),
any unit vector \( \av \in \R^{\dimp} \), and any \( \zq > 0 \), it holds
%with \( \zq^{2} \geq  2 \dimp + 2 \| \gaussv \|^{2} \) 
\begin{EQA}
\label{gampxiz12}
	\P\bigl( \| \gaussv - \uv \| \geq \zq \bigr)
	& \leq &
	\exp\bigl\{ - \zq^{2}/4 + \dimp/2 + \| \uv \|^{2}/2 \bigr\} ,
	\\
	\E \bigl\{ | \gaussv^{\T} \av |^{2} \Ind\bigl( \| \gaussv - \uv \| \geq \zq \bigr) \bigr\}
	& \leq &
	(2 + |\uv^{\T} \av|^{2}) \exp\bigl\{ - \zq^{2}/4 + \dimp/2 + \| \uv \|^{2}/2 \bigr\} .
\label{gampxiz122}
\end{EQA}
\end{lemma}

%\ifbook{}{
\begin{proof}
By the exponential Chebyshev inequality, for any \( \lambda < 1 \) 
\begin{EQA}
	\P\bigl( \| \gaussv - \uv \| \geq \zq \bigr)
	& \leq &
	\exp \bigl( - \lambda \zq^{2}/2 \bigr) \E \exp\bigl( \lambda \| \gaussv - \uv \|^{2}/2 \bigr)
	\\
	&=&
	\exp\Bigl\{ 
		- \frac{\lambda \zq^{2}}{2} - \frac{\dimp}{2} \log(1 - \lambda) 
		+ \frac{\lambda}{2 (1-\lambda)} \| \uv \|^{2}
	\Bigr\} .
\label{gampxiz}
\end{EQA}
In particular, with \( \lambda = 1/2 \), this implies \eqref{gampxiz12}.
Further, for \( \| \av \| = 1 \)
\begin{EQA}
	\E \bigl\{ | \gaussv^{\T} \av |^{2} \Ind(\| \gaussv - \uv \| \geq \zq) \bigr\}
	& \leq &
	\exp \bigl( - \zq^{2}/4 \bigr) 
	\E \bigl\{ | \gaussv^{\T} \av |^{2} \exp\bigl( \| \gaussv - \uv \|^{2}/4 \bigr) \bigr\}
	\\
	& \leq &
	(2 + | \uv^{\T} \av |^{2}) \exp \bigl( - \zq^{2}/4 + \dimp/2 + \| \uv \|^{2}/2 \bigr) 
\label{Egaguq222}
\end{EQA}
and \eqref{gampxiz122} follows.
\end{proof}
%}

%%%%%%%%%%%%%%%%%%%%%%%%%%%%%%%%%%%%%%%%%%%%%%%%%%%%%%%%%%%%%%%%

%%%%%%%%%%%%%%%%%%%%%%%%%%%%%%%%%%%%%%%%%%%%%%%%%%%%%%%%%%%%%%%%
% !TEX root = script2014.tex

\section{Deviation bounds for Gaussian quadratic forms}
\label{SmatrBern}
\label{Sprobabquad}
This section collects some deviation bounds for Gaussian quadratic forms.
The next result explains the concentration effect of 
\( \gaussv^{\T} \BB \gaussv \)
for a standard Gaussian vector \( \gaussv \) and a symmetric matrix \( \BB \).
We use a version from \cite{laurentmassart2000}.

%For notational simplicity we assume that \( \BB \) is symmetric.
%Otherwise one should replace it with \( (\BB^{\T} \BB)^{1/2} \).

\begin{theorem}
\label{TexpbLGA}
\label{Lxiv2LD}
\label{Cuvepsuv0}
Let \( \gaussv \) be a standard normal Gaussian vector and \( \BB \) be symmetric positive.
Then with \( \dimA = \tr(\BB) \), \( \vA^{2} = \tr(\BB^{2}) \), and 
\( \supA = \| \BB \|_{\oper} \), it holds for each \( \xx \geq 0 \)
\begin{EQA}
	\P\bigl( \gaussv^{\T} \BB \gaussv > \dimA + 2 \vA \xx^{1/2} + 2 \supA \xx \bigr)
	& \leq &
	\ex^{-\xx} .
\label{Pxiv2dimAvp12}
\end{EQA}
This implies for any positive \( \BB \) %with \( \| \BB \|_{\oper} \leq 1 \)
\begin{EQA}
	\P\bigl( \| \BB^{1/2} \gaussv \| > \dimA^{1/2} + (2 \supA \xx)^{1/2} \bigr)
	& \leq &
	\ex^{-\xx} .
\label{Pxiv2dimAxx12}
\end{EQA}
Also
\begin{EQA}
	\P\bigl( \gaussv^{\T} \BB \gaussv < \dimA - 2 \vA \xx^{1/2} \bigr)
	& \leq &
	\ex^{-\xx} .
\label{Pxiv2dimAvp12m}
\end{EQA}
If \( \BB \) is symmetric but non necessarily positive then
\begin{EQA}
	\P\bigl( \bigl| \gaussv^{\T} \BB \gaussv - \dimA \bigr| > 2 \vA \xx^{1/2} + 2 \supA \xx \bigr)
	& \leq &
	2 \ex^{-\xx} .
\label{PxivTBBdimA2vp}
\end{EQA}
\end{theorem}

%\ifbook{}
{
\begin{proof}
Normalisation by \( \supA \) reduces the statement to the case with \( \supA = 1 \).
Further, the standard rotating arguments allow to reduce the Gaussian quadratic form 
\( \| \gaussv \|^{2} \) to the chi-squared form:
\begin{EQA}
	\gaussv^{\T} \BB \gaussv
	&=&
	\sum_{j=1}^{\dimp} \lambda_{j} \nu_{j}^{2}
\label{xiv2sj1p}
\end{EQA}
with independent standard normal r.v.'s \( \nu_{j} \).
Here \( \lambda_{j} \in [0,1] \) are eigenvalues of \( \BB \), and 
\( \dimA = \lambda_{1} + \ldots + \lambda_{\dimp} \), 
\( \vA^{2} = \lambda_{1}^{2} + \ldots + \lambda_{\dimp}^{2} \).
One can easily 
compute the exponential moment of \( (\gaussv^{\T} \BB \gaussv - \dimA)/2 \):
for each positive \( \mu < 1 \)
\begin{EQA}
	\log \E \exp\bigl\{ \mu (\gaussv^{\T} \BB \gaussv - \dimA)/2 \bigr\}
	&=&
	\frac{1}{2} \sum_{j=1}^{\dimp} \bigl\{ - \mu \lambda_{j} - \log(1 - \mu \lambda_{j}) \bigr\} .
\label{lEemux2p2}
\end{EQA}

\begin{lemma}
Let \( \mu \lambda_{j} < 1 \) and \( \lambda_{j} \leq 1 \).
Then 
\begin{EQA}
	\frac{1}{2} \sum_{j=1}^{\dimp} \bigl\{ - \mu \lambda_{j} - \log(1 - \mu \lambda_{j}) \bigr\}
	& \leq &
%	\sum_{j=1}^{\dimp} \frac{(\mu \lambda_{j})^{2}}{4 (1 - \mu)} 
%	\leq 
	\frac{\mu^{2} \vA^{2}}{4 (1 - \mu)} \, .
\label{jmu2v221mu}
\end{EQA}
\end{lemma}

\begin{proof}
In view of \( \mu \lambda_{j} < 1 \), it holds for every \( j \)
\begin{EQA}
	- \mu \lambda_{j} - \log(1 - \mu \lambda_{j}) 
	&=&
	\sum_{k=2}^{\infty} \frac{(\mu \lambda_{j})^{k}}{k}
	\\
	& \leq &
	\frac{(\mu \lambda_{j})^{2}}{2}
	\sum_{k=0}^{\infty} (\mu \lambda_{j})^{k}
	\leq 
	\frac{(\mu \lambda_{j})^{2}}{2 (1 - \mu \lambda_{j})} 
	\leq 
	\frac{(\mu \lambda_{j})^{2}}{2 (1 - \mu)},
\label{jmu2v221mup}
\end{EQA}
and thus
\begin{EQA}
	\frac{1}{2} \sum_{j=1}^{\dimp} \bigl\{ - \mu \lambda_{j} - \log(1 - \mu \lambda_{j}) \bigr\}
	& \leq &
	\sum_{j=1}^{\dimp} \frac{(\mu \lambda_{j})^{2}}{4 (1 - \mu)} 
	\leq 
	\frac{\mu^{2} \vA^{2}}{4 (1 - \mu)} \, .
\label{sjmu2v221mu}
\end{EQA}
\end{proof}
The next technical lemma is helpful.

\begin{lemma}
\label{Lmuvpxx}
For each \( \vA > 0 \) and \( \xx > 0 \), it holds
\begin{EQA}
	\inf_{\mu > 0} \biggl\{ 
		- \mu \bigl( \vA \xx^{1/2} + \xx \bigr) + \frac{\mu^{2} \vA^{2}}{4 (1 - \mu)} 
	\biggr\}
	& \leq &
	- \xx .
\label{infmuxxvp}
\end{EQA}
\end{lemma}

\begin{proof}
Let pick up 
\begin{EQA}
	\mu 
	&=& 
	1 - \frac{1}{2\xx^{1/2}/\vA + 1} = \frac{\xx^{1/2}}{\xx^{1/2} + \vA/2} , 
\label{mu12xx12vp1m1}
\end{EQA}
so that \( \mu / (1 - \mu) = 2 \xx^{1/2}/\vA \). Then
\begin{EQA}
	&& \nquad
	- \mu \bigl( \vA \xx^{1/2} + \xx \bigr) + \frac{\mu^{2} \vA^{2}}{4 (1 - \mu)}
	\\
	&=&
	- \mu \bigl( \vA \xx^{1/2} + \xx + \vA^{2}/4 \bigr)
	+ \frac{\mu \vA^{2}}{4 (1 - \mu)}
	\\
	&=&
	- \frac{\xx^{1/2}}{\xx^{1/2} + \vA/2} \bigl( \xx^{1/2} + \vA/2 \bigr)^{2} 
	+ \frac{2 \xx^{1/2} \vA }{4}
	=
	- \xx 
\label{mux2xv4x12}
\end{EQA}
and the result follows.
\end{proof}

Now we apply the Markov inequality 
\begin{EQA}
	&& \nquad
	\log \P\bigl( \gaussv^{\T} \BB \gaussv > \dimA + 2 \vA \xx^{1/2} + 2 \xx \bigr)
	=
	\log \P\bigl( (\gaussv^{\T} \BB \gaussv - \dimA) / 2 > \vA \xx^{1/2} + \xx \bigr)
	\\
	& \leq &
	\inf_{\mu > 0} \biggl\{ 
		- \mu \bigl( \vA \xx^{1/2} + \xx \bigr) 
		+ \log\E \exp\bigl\{ \mu (\gaussv^{\T} \BB \gaussv - \dimA)/2 \bigr\}
	\biggr\}
	\\
	& \leq &
	\inf_{\mu > 0} \biggl\{ 
		- \mu \bigl( \vA \xx^{1/2} + \xx \bigr) + \frac{\mu^{2} \vA^{2}}{4 (1 - \mu)}
	\biggr\}
	\leq 
	- \xx
\label{x2xv4x12}
\end{EQA}
and the first assertion \eqref{Pxiv2dimAvp12} follows.
The second statement follows from the first one by 
\( \tr(\BB^{2}) \leq \| \BB \|_{\oper} \tr(\BB) = \supA \, \dimA \).

Similarly for any \( \mu > 0 \)
\begin{EQA}
	\P\bigl( \gaussv^{\T} \BB \gaussv - \dimA < - 2 \vA \sqrt{\xx} \bigr)
	& \leq &
	\exp\bigl( - \mu \vA \sqrt{\xx} \bigr)
	\E \exp\Bigl( - \frac{\mu}{2} (\gaussv^{\T} \BB \gaussv - \dimA) \Bigr) .
%	\\
%	& \leq &
%	\exp\bigl( - \mu \vA \sqrt{\xx} \bigr)
%	\E \exp\Bigl( - \frac{\mu}{2} (\| \gaussv \|^{2} - \dimA) \Bigr) .
\end{EQA}
By \eqref{lEemux2p2}
\begin{EQA}
	\log \E \exp\bigl\{ - \mu (\gaussv^{\T} \BB \gaussv - \dimA)/2 \bigr\}
	&=&
	\frac{1}{2} \sum_{j=1}^{\dimp} 
		\bigl\{ \mu \lambda_{j} - \log(1 + \mu \lambda_{j}) \bigr\} .
\label{lEemux2p2m}
\end{EQA}
and 
\begin{EQA}
	\frac{1}{2} \sum_{j=1}^{\dimp} \bigl\{ \mu \lambda_{j} - \log(1 + \mu \lambda_{j}) \bigr\}
	&=&
	\frac{1}{2} \sum_{j=1}^{\dimp} \sum_{k=2}^{\infty} \frac{(- \mu \lambda_{j})^{k}}{k}
	\leq 
	\sum_{j=1}^{\dimp}\frac{(\mu \lambda_{j})^{2}}{4} 
	=
	\frac{\mu^{2} \vA^{2}}{4} .
\label{jmu2v221mum}
\end{EQA}
Here the choice \( \mu = 2 \sqrt{\xx} / \vA \) yields \eqref{Pxiv2dimAvp12m}.

One can put together the arguments used for obtaining the lower and the upper bound 
for getting a bound for a general 
quadratic form \( \gaussv^{\T} \BB \gaussv \), where \( \BB \) is symmetric but not necessarily 
positive.
\end{proof}
}
%
%\begin{corollary}
%Let \( \gaussv \) be standard normal in \( \R^{\dimp} \) and \( \BB \) be symmetric. 
%Then with \( \dimA = \tr(\BB) \), \( \vA^{2} = \tr(\BB^{2}) \), and 
%\( \supA = \supA_{\max}(\BB) \), it holds for each \( \xx \geq 0 \)
%\end{corollary}
%

Finally we apply this result to weighted sums of centered \( \gauss_{i}^{2} \).
\begin{corollary}
\label{Cuvepsuv}
For any unit vector \( \uv = (u_{i}) \in \R^{n} \) 
and standard normal r.v.'s \( \gauss_{i} \), it holds with 
\( \| \uv \|_{\infty} \eqdef \max_{i} |u_{i}| \)
%with \( \vB_{\uv}^{2} = \sum_{i} u_{i}^{4} \):
%and \( a^{*} = \max u_{i}^{2} = \| \uv \|_{\infty}^{2} \):
\begin{EQA}
	\P\biggl( 
		\biggl| \sum_{i=1}^{n} u_{i} (\gauss_{i}^{2} - 1) \biggr| 
		\geq  
		2 \xx^{1/2} + 2 \| \uv \|_{\infty} \xx
	\biggr)
	& \leq &
	2 \ex^{-\xx} .
\label{Puitei2m1v46x}
\label{Puitei2m1v46xx}
\end{EQA}
%Moreover, if  \( \| \uv \|_{\infty} \leq \td \), then
%\begin{EQA}
%	\P\biggl( 
%		\biggl| \sum_{i=1}^{n} u_{i}^{2} (\gauss_{i}^{2} - 1) \biggr| 
%		\geq  
%		2 \td \xx^{1/2} + 2 \td^{2} \, \xx
%	\biggr)
%	& \leq &
%	2 \ex^{-\xx} .
%\label{Puitei2m1v46xx}
%\end{EQA}
\end{corollary}

\begin{proof}
The statement follows directly from Theorem~\ref{Cuvepsuv0}.
It suffices to notice \( \vA^{2} = \| \uv \|^{2} = 1 \).
%One can also bound %for any unit vector \( \uv \)
%\begin{EQA}
%	\vB_{\uv}^{2}
%	&=&
%	\sum_{i=1}^{n} u_{i}^{4}
%	\leq 
%	\| \uv \|_{\infty}^{2} \sum_{i=1}^{n} u_{i}^{2}
%	=
%	\| \uv \|_{\infty}^{2} 
%	\leq 
%	\td^{2}.
%\label{vBuv242}
%\end{EQA}
%This implies \eqref{Puitei2m1v46xx}.
\end{proof}

As a special case, we present a bound for the chi-squared distribution 
corresponding to \( \BB = \Id_{\dimp} \).
Then \( \tr (\BB) = \dimp \), \( \tr(\BB^{2}) = \dimp \) and \( \supA(\BB) = 1 \).

\begin{corollary}
\label{Cchi2p}
Let \( \gaussv \) be a standard normal vector in \( \R^{\dimp} \).
Then
\begin{EQA}[lcl]
\label{Pxi2pm2px}
	\P\bigl( \| \gaussv \|^{2} \geq \dimp + 2 \sqrt{\dimp \xx} + 2 \xx \bigr)
	& \leq &
	\ex^{-\xx},
	\\
	\P\bigl( \| \gaussv \| \,\,  \geq \sqrt{\dimp} + \sqrt{2 \xx} \bigr)
	& \leq &
	\ex^{-\xx} ,
\label{Pxi2pm2px12}
	\\
	\P\bigl( \| \gaussv \|^{2} \leq \dimp - 2 \sqrt{\dimp \xx} \bigr)
	& \leq &
	\ex^{-\xx}	.
\label{Pxi2pm2px22}
\end{EQA}
\end{corollary} 

The previous results are mainly stated for a standard Gaussian vector \( \gaussv \in \R^{n} \).
Now we extend it to the case of a zero mean Gaussian vector \( \xiv \) with the 
\( n \times n \) covariance matrix \( \Covm = (\covm_{ij}) \) with 
\( \lambda_{\max}(\Covm) \leq \supAB \).
Given a unit vector \( \uv = (u_{1},\ldots,u_{n})^{\T} \in \R^{n} \), consider the quadratic form
\begin{EQA}
	\GQF
	&=&
	\sum_{i=1}^{n} u_{i} \xi_{i}^{2} .
\label{S1nuiei2d}
\end{EQA}
We aim at bounding \( \GQF - \E \GQF \).
To apply the result of Theorem~\ref{Cuvepsuv0} represent \( \GQF \)
as \( \gaussv^{\T} \BB \gaussv \) with \( \BB \) depending on 
\( \uv \) and \( \Covm \).
More precisely, let \( \xiv = \Covm^{1/2} \gaussv \) for a standard Gaussian vector 
\( \gaussv \in \R^{n} \).
Then with \( \Uv = \diag(u_{1},\ldots,u_{n}) \), it holds
\begin{EQA}
	S
	&=&
	\tr\bigl( \Uv \xiv \xiv^{\T} \bigr)
	=
	\tr\bigl( \Uv \Covm^{1/2} \gaussv \gaussv^{\T} \Covm^{1/2} \bigr)
	=
	\tr\bigl( \BB \gaussv \gaussv^{\T} \bigr)
	=
	\gaussv^{\T} \BB \gaussv
\label{•}
\end{EQA}
with \( \BB = \Covm^{1/2} \Uv \Covm^{1/2} \).
Therefore, the bound \( \| \Covm \|_{\oper} \leq \supAB \) implies 
\begin{EQA}
	\supA
	&=&
	\supA(\BB)
	=
	\| \Covm^{1/2} \Uv \Covm^{1/2} \|_{\oper}
	\leq 
	\supAB \, \| \uv \|_{\infty} \, ,
	\\
	\vA^{2}
	&=&
	\tr(\BB^{2})
	=
	\tr\bigl( \Covm^{1/2} \Uv \Covm \Uv \Covm^{1/2} \bigr)
	\leq 
	\supAB \tr\bigl( \Uv \Covm \Uv \bigr)
	\leq 
	{\supAB}^{2} \| \uv \|^{2} = {\supAB}^{2}.
\label{supAvp2BB2}
\end{EQA}
Now the general results of Theorem~\ref{Cuvepsuv0} implies the result 
similar to Corollary~\ref{Cuvepsuv}.
 
\begin{corollary}
\label{CuvepsuvnG}
For any unit vector \( \uv = (u_{i}) \in \R^{n} \), \( \| \uv \| = 1 \),
and normal zero mean vector \( \xiv \sim \ND(0,\Covm) \) in \( \R^{n} \) with 
\( \| \Covm \|_{\oper} \leq \supAB \), it holds 
\begin{EQA}
	\P\biggl( 
		\biggl| \sum_{i=1}^{n} u_{i} (\xi_{i}^{2} - \E \xi_{i}^{2}) \biggr| 
		\geq  
		2 \supAB \, \xx^{1/2} + 2 \supAB \, \| \uv \|_{\infty} \xx
	\biggr)
	& \leq &
	2 \ex^{-\xx} .
\label{Puitei2m1v46xnG}
\end{EQA}
\end{corollary}

It is worth noting that the identity \( \| \uv \| = 1 \) implies 
\( \| \uv \|_{\infty} \leq 1 \).
Moreover, in typical situations,
\( \| \uv \|_{\infty} \asymp n^{-1/2} \), and the leading term in the bounds of 
Corollaries~\ref{Cuvepsuv} and \ref{CuvepsuvnG} is \( 2 \supAB \, \xx^{1/2} \).

\def\tdPi{\tdn}
\section{Sums of random matrices}
Here we present a number of deviation bounds for a sum of random matrices.

\subsection{Matrix Bernstein inequality}
This section collects some useful facts about deviation of stochastic matrices
from their mean. 
We mainly use the arguments from the book \cite{Tropp2014}.
The main step of the proof is the following Master bound.

\begin{theorem}[Master bound]
\label{MastBound}
Assume that  \(\Sv_{1},\dots,\Sv_{n}\) are independent Hermitian matrices of the same size and  
\(\Zv = \sum_{i=1}^{n} \Sv_{i} \). 
Then 
\begin{EQA}
	\E \supA^+_{\max}(\Zv)
	& \leq &
	\inf_{\theta>0}\frac{1}{\theta} 
	\log \tr\exp\left(\sum_{i=1}^{n}\log\E\ex^{\theta \Sv_{i}}\right),
\label{MaCheEx}
	\\
	\P\{\supA_{\max}^+(\Zv) \geq \zq\} 
	& \leq & 
	\inf_{\theta>0} \ex^{-\theta \zq} \, 
		\tr\exp\left(\sum_{i=1}^{n}\log\E \ex^{\theta \Sv_{i}}\right),
\label{MaCheProb}
\end{EQA}
where \( \supA^+_{\max}(\Zv) \) denotes the algebraically largest eigenvalue of \(\Zv\).
\end{theorem}

\ifbook{For the proof see e.g. \cite{Tropp2014}.}
{
\begin{proof} 
By the Markov inequality 
\begin{EQA}
	\P\{\supA^+_{\max}(\Zv) \geq \zq\} 
	&\leq&
 	\inf_{\theta} \ex^{-\theta \zq} \E \exp(\theta \supA^+_{\max}(\Zv)).
\end{EQA}
Recall the spectral mapping theorem: for any function \( f\colon \R \to\R \) 
and Hermitian matrix \( A \) eigenvalues of \(f(A)\) are equal to eigenvalues of \(A\). 
Thus 
\begin{EQA}
	\exp(\theta \supA_{\max}(\Zv))
 	&= &
 	\exp( \supA^+_{\max}(\theta\Zv) ) 
 	= 
 	\supA^+_{\max}\bigl( \exp(\theta \Zv ) \bigr) 
	\leq 
	\tr \, \ex^{\theta \Zv} .
\end{EQA}
Therefore,
\begin{EQA}
	\P\{\supA^+_{\max}(\Zv) \geq \zq\}
	&\leq &
	\inf_{\theta} \ex^{-\theta \zq} \E\tr\exp(\theta \Zv),
\label{Tr3.2.1}
\end{EQA}
and \eqref{MaCheProb} follows.

To prove \eqref{MaCheEx} fix \(\theta\). 
Using the spectral mapping theorem one can get that 
\begin{EQA}
	\E \supA^+_{\max}(\Zv) 
	& = &
	\frac{1}{\theta}\E \supA^+_{\max}(\theta \Zv)  
	= 
	\frac{1}{\theta} \log\E\exp \bigl( \supA^+_{\max}(\theta \Zv) \bigr)
	=
	\frac{1}{\theta} \log\E \supA^+_{\max}\bigl( \exp (\theta \Zv) \bigr) \, . 
\label{Tr3.2.2}
\end{EQA}
Thus we get 
\begin{EQA}
	\E\supA^+_{\max}(\Zv)
	& \leq &
	\frac{1}{\theta} \log \tr \E \exp(\theta \Zv).
\label{Tr3.2.2.1}
\end{EQA}
The final step in proving the master inequalities is to bound from above 
\( \E\tr\exp\left(\sum_{i=1}^{n} \Sv_{i} \right)\). 
To do this we use Jensen's inequality for the convex function  
\( \tr \exp(H + \log(X)) \) (in matrix \(X\)), where \(H\) is deterministic Hermitian  matrix.  For a  random Hermitian matrix \( X \) one can write
\begin{EQA}
	\E \tr\exp(H + X) 
	&=& 
	\E \tr \exp(H + \log \ex^{X})
	\leq   
	\tr \exp(H + \log \E \ex^{X}).
\label{Tr3.4.1}
\end{EQA}
Denote by \( \E_{i} \) the conditional expectation with respect to  random matrix \( X_{i} \).  To bound \( \E \tr\exp\left(\sum_{i=1}^{n} \Sv_{i}\right) \)  we use \eqref{Tr3.4.1} for the sum of independent Hermitian matrices by taking the conditional expectations with respect to \( i \)-th matrix:
\begin{EQA}
	\E\tr\exp\left(\sum_{i=1}^{n} \Sv_{i}\right) 
 	&= &
	\E\E_{n}\tr\exp\left(\sum_{i=1}^{n-1} \Sv_{i} + \Sv_{n}\right) 
	\\ 
 	& \leq & 
	\E\tr\exp\left(\sum_{i=1}^{n-1} \Sv_{i} + \log(\E_{n}\exp(\Sv_{n}))\right) 
	\\
	& \leq &  
	\tr\exp\left(\sum_{i=1}^{n}\log\E \ex^{\theta \Sv_{i}} \right).
\label{Tr3.4.1}
\end{EQA}
To complete the prove of the Master's theorem combine  \eqref{Tr3.2.1} and \eqref{Tr3.2.2.1} with \eqref{Tr3.4.1}.
\end{proof}
}
\medskip

The same result applied to \( - \Zv \) yields the bound for the operator norm \( \| \Zv \| \):
\begin{EQA}
	\P\{ \| \Zv \|_{\oper} \geq \zq\} 
	& \leq & 
	\inf_{\theta>0} \ex^{-\theta \zq} \, 
		\tr\exp\left(\sum_{i=1}^{n}\log\E \ex^{\theta \Sv_{i}}\right)  \\
		& & +  \inf_{\theta>0} \ex^{-\theta \zq} \tr\exp\left(\sum_{i=1}^{n}\log\E \ex^{-\theta \Sv_{i}}\right) .
\label{MaCheProb2}
\end{EQA}

\ifbook{} %See again \cite{Tropp2014}.}
{
\begin{theorem}[Bernstein inequality for a sum of random Hermitian  matrices]
\label{BernSqTh}
\label{TmatrBern}
Let \(\Zv = \sum_{i=1}^{n} \Sv_{i}\), where \(\Sv_{i}\), \(i=1,\dots,n\) are independent, random, Hermitian  matrices of the dimension \(d\times d\) and 
\begin{EQA}
 	\supA^+_{\max}(\Sv_{i}) 
	& \leq & 
	R.
\end{EQA}
Denote \( \vp^{2} = \vp^{2}(\Zv) = \|\E (\Zv^{2}) \|_{\oper} \). 
Then 
\begin{EQA}
 	\E\supA^+_{\max}(\Zv)
	& \leq  & 
	\sqrt{2\vp^{2}\log(d)}+\frac{1}{3}R\log(d),
\label{BernSqE}
	\\
	\P \bigl\{ \supA^+_{\max}(\Zv) \geq \zq \bigr\} 
	& \leq & 
	d \exp \left( \frac{-\zq^{2}/2}{\vp^{2} + R \zq/3} \right).
\label{BernSqP}
\end{EQA}
\end{theorem}

\begin{proof}
Note that 
\begin{EQA}
	\vp^{2} 
	&=& 
	\left\|\sum_{i=1}^{n}\E \Sv_{i}^{2}\right\|_{\oper} \, .
\end{EQA}
For the sake of simplicity let \( \vp^{2} = 1 \). 
Denote 
\begin{EQA}
	g(\theta)
 	&= &
 	\frac{\theta^{2}/2}{1-R\theta/3}.
\end{EQA}
Apart the Master inequalities, we use the following lemma:

\begin{lemma}
\label{PreBern}
Let \(\Zv\) be a random Hermitian  matrix \(\E \Zv = 0 \), \( \supA^+_{\max}(\Zv) \leq R\), then for \(0< \theta<3/R \) the following inequalities hold
\begin{EQA}
	\E \ex^{\theta \Zv}
	&\leq &
 	\exp\left(\frac{\theta^{2}/2}{1-R\theta/3} \E(\Zv^{2})\right),
	\\
	\log \E \ex^{\theta \Zv}
	&\leq&
 	\frac{\theta^{2}/2}{1-R\theta/3} \E(\Zv^{2}).
\end{EQA}
\end{lemma}
 
\begin{proof}
Decompose the exponent in the following way
\begin{EQA}
	\ex^{\theta \Zv}
 	&= &
	\Id + \theta \Zv + (\ex^{\theta \Zv}-\theta \Zv - \Id) 
	= 
	\Id + \theta \Zv + \Zv \cdot f(\Zv)\cdot \Zv,
\end{EQA}
where 
\begin{EQA}
 	f(x) 
	& = &
	\frac{\ex^{\theta x}-\theta x - 1}{x^{2}}, \quad \text{for} \quad x\not=0, \quad 
	f(0) = \frac{\theta^{2}}{2}.
\end{EQA}
One can check that the function \(f(x)\) is non-decreasing, thus for \(x\leq R\), one has 
\(f(x)\leq f(R)\). 
By the matrix transfer rule 
\( f(\Zv)\leq f(R) \Id \) and 
\begin{EQA}
 	\E \ex^{\theta \Zv} 
	& = & 
	\Id + f(R) \E \Zv^{2}.
\end{EQA}
In order to estimate \(f(R)\)  use \(q!\geq 2\cdot 3^{q-2}\) to get 
\begin{EQA}
  	f(R) 
	& = &
	\frac{\ex^{\theta R}-\theta R - \Id}{R^{2}} 
	= 
	\frac{1}{R^{2}}\sum_{q=2}^{\infty}\frac{(\theta R)^q}{q!}
	\leq 
	\theta^{2} \sum_{q=2}^{\infty} \frac{(R\theta)^{q-2}}{3^{q-2}} 
	= \frac{\theta^{2}/2}{1-R\theta/3}.
\end{EQA}
To get the result of the Lemma note that \(1+a\leq \ex^{a}\).
\end{proof}

To prove \eqref{BernSqE} and \eqref{BernSqP} we apply the Master inequalities and 
Lemma \ref{PreBern}:
\begin{EQA}
	\E\supA^+_{\max}(\Zv)
	& \leq &
	\inf_{\theta>0}\frac{1}{\theta} 
	\log \tr\exp \left( \sum_{i=1}^{n} \log \E\exp(\theta \Sv_{i})\right) 
	\\
	& \leq  &
	\inf_{0<\theta<3/R}\frac{1}{\theta} 
	\log\tr\exp \left( g(\theta) \sum_{i=1}^{n} \E \Sv_{i}^{2}\right) 
	\\
	& \leq  &
	\inf_{0<\theta<3/R}\frac{1}{\theta}\log\tr\exp \left( g(\theta) \E \Zv^{2}\right) 
	\\
	& \leq & 
	\inf_{0<\theta<3/R}\frac{1}{\theta}\log d \exp \left( g(\theta) \|\E \Zv^{2}\|_{\oper} \right)
	\\ 
	& \leq &
	\inf_{0<\theta<3/R}\left\{\frac{\log(d)}{\theta} + \frac{\theta/2}{1-R\theta/3}\right\}. 
\end{EQA}
Minimizing the right hand side in \(\theta\) one can get \eqref{BernSqE}.

The second inequality can be obtained in the same manner:
\begin{EQA}
	\P\{\supA^+_{\max}(\Zv) \geq \zq \}
	& \leq  & 
	\inf_{\theta>0}\ex^{-\theta \zq}  \tr\exp \left( \sum_{i=1}^{n} \log \E\exp(\theta \Sv_{i})\right) 
	\\
	& \leq  & 
	\inf_{0<\theta<3/R}\ex^{-\theta \zq}  \tr \exp \left( g(\theta) \E \Zv^{2}\right) 
	\\
	& \leq & 
	\inf_{0<\theta<3/R}\ex^{-\theta \zq}  d \exp \left( g(\theta) \|\E \Zv^{2}\|_{\oper}\right)
	\\ 
	& \leq & 
	\inf_{0<\theta<3/R}\ex^{-\theta \zq}  d \exp \left( g(\theta)\right).
\end{EQA}
Here instead of minimizing the right hand side in \(\theta\) we have used 
\(\theta = \zq/{(1+R \zq/3)}\).
\end{proof}

\begin{theorem}[Bernstein inequality for a sum of random Hermitian  matrices]
\label{TmatrBern}
Let \(\Zv = \sum_{i=1}^{n}\Sv_{i}\), where \(\Sv_{i}\), \(i=1,\dots,n\) are  independently distributed random matrices of the size \(d_{1}\times d_{2}\) and 
\begin{EQA}
 	\|\Sv_{i}\|_{\oper} 
 	& \leq &
  	R. 
\end{EQA}
Denote \(\vp^{2} = \vp^{2}(\Zv) 
= \max\left\{\|\E(\Zv^{*} \Zv)\|_{\oper}, \|\E(\Zv \Zv^{*})\|_{\oper}\right\} \). 
Then 
\begin{EQA}
 	\E \| \Zv \|_{\oper}
	& \leq & 
	\sqrt{2\vp^{2}\log(d_{1}+d_{2})} + \frac{1}{3} R\log(d),
\label{BernSqE12}
  	\\
	\P\{\| \Zv \|_{\oper} \geq \zq \} 
	& \leq & 
	(d_{1}+d_{2}) \exp \left( \frac{-\zq^{2}/2}{\vp^{2} + R \zq/3}\right).
\label{BernSqP12}
\end{EQA}
\end{theorem}

\begin{proof}
Use the following hint: define the matrix
\begin{EQA}
 	H(\Zv) 
 	&=& 
 	\begin{pmatrix} 
		0 & \Zv 
		\\ 
		\Zv^{*} & 0 
	\end{pmatrix}.
\end{EQA}
It can be easily seen that \(\vp^{2} =  \|H(\Zv)^{2}\|_{\oper}\), and 
\(\|\Zv\|_{\oper} = \supA^+_{\max}\bigl( H(\Zv) \bigr) \), 
thus applying Proposition \ref{BernSqTh} to \( H(\Zv) \) the statements \eqref{BernSqE} and \eqref{BernSqP} are straightforward.
\end{proof}
}

\subsection{Matrix deviation bounds}
The next result provides a deviation bound for a matrix-valued quadratic forms. 

\begin{proposition}[Deviation bound for matrix quadratic forms]
\label{CUvepsUv}
Consider a \( \dimp \times n \) matrix \( \UV \) such that
\begin{EQA}
\label{OrthoU}
	\UV \UV^{\T}
	&=&
	\Id_{\dimp} .
\end{EQA}
Let the columns \( \UVcol_{1},\dots,\UVcol_{n} \in \R^{\dimp} \) of the matrix \( \UV \) satisfy
\begin{EQA}
\label{BoundRow}
 	\|\UVcol_{i}\|
 	&\leq &
 	\tdn
\end{EQA}
for a fixed constant \( \tdn \).
For a random vector \( \gaussv = (\gauss_{1},\ldots,\gauss_{n})^{\T} \) with independent standard Gaussian components, define
\begin{EQA}
	\Zv 
	& \eqdef &
	\UV \diag\bigl\{ \gaussv \cdot \gaussv - 1 \bigr\} \UV^{\T}
	=
	\sum_{i=1}^{n} (\gauss_{i}^{2} - 1) \UVcol_{i}\UVcol_{i}^{\T}.
\label{ZvUVecem1UVT}
\end{EQA}
Then with \( \xxp = \xx + \log (2\dimp) \)
\begin{EQA}
	\P\Bigl( 
		\| \Zv \|_{\oper}
		\geq  
		2 \tdn \sqrt{ \xxp } + 2 \tdn^{2} \xxp 
	\Bigr)
	& \leq &
	\ex^{-\xx}.
\label{Zvop2tdxx2p}
\end{EQA}
%\end{lemma}
\end{proposition}

\begin{proof}
From the Master bound \eqref{MaCheProb2} %(see Proposition \ref{MastBound})
\begin{EQA}
\label{matrix_chernoff}
	\P\bigl( 
		\| \Zv \|_{\oper}
		\geq 
		\zq
	\bigr)
	& \leq & \inf_{\theta> 0} \ex^{-\theta \zq} 
		\tr \exp \left(
		\sum_{i=1}^{n} \log \E\exp( \theta (\gauss_{i}^{2}-1) \UVcol_{i} \UVcol_{i}^{\T})\right)\\
		& +& 
		\inf_{\theta> 0} \ex^{-\theta \zq} 
		\tr \exp \left( 
		\sum_{i=1}^{n} \log \E\exp( \theta (-\gauss_{i}^{2}+1) \UVcol_{i} \UVcol_{i}^{\T})
	\right).
\end{EQA}
Now we use the following general fact:

\begin{lemma}
\label{LlEexiUv}
If \( \chi \) is a random variable and \( \Pi \) is a projector in \( \R^{\dimp} \), then 
\begin{EQA}
	\log \E \exp(\chi \Pi) 
	&=&
	\log \bigl( \E \ex^{\chi} \bigr) \Pi .
\label{logEexxiUv}
\end{EQA}
\end{lemma}
\begin{proof}
The result \eqref{logEexxiUv} can be easily obtained by applying twice the spectral mapping theorem. 
\end{proof}

This result yields, in particular, for any unit vector \( \UVcol \in \R^{\dimp} \)
\begin{EQA}
	\log \E \exp\bigl( \chi \UVcol \UVcol^{\T} \bigr)
	&=&
	\log \bigl( \E \ex^{\chi} \bigr) \UVcol \UVcol^{\T} .
\label{logEexxiuvuvT}
\end{EQA}
Moreover, for any vector \( \UVcol \in \R^{\dimp} \), 
the normalized product \( \UVcol \UVcol^{\T} / \| \UVcol \|^{2} \) is a rank-one projector, 
and hence,
\begin{EQA}
	\log \E \exp\bigl( \chi \UVcol \UVcol^{\T} \bigr)
	&=&
	\log \bigl( \E \ex^{\chi \| \UVcol \|^{2}} \bigr) \frac{\UVcol \UVcol^{\T}}{\| \UVcol \|^{2}} \, .
\label{logEexxiuv2uv2}
\end{EQA}
With \( \Uv_{i} \eqdef \UVcol_{i} \UVcol_{i}^{\T} / \| \UVcol_{i} \|^{2} \) and 
\( \chi_{i} = \theta (\gauss_{i}^{2} - 1) \), we derive
\begin{EQA}
	\log \E\exp\bigl\{ \theta (\gauss_{i}^{2} - 1) \UVcol_{i} \UVcol_{i}^{\T} \bigr\}
	&=&
	\log \E \exp\bigl\{ \theta (\gauss_{i}^{2} - 1) \| \UVcol_{i} \|^{2} \bigr\} \Uv_{i}
	\\
	&=&
	\log\left(  
		\frac{\exp\bigl( - \| \UVcol_{i} \|^{2}\theta\bigr)}{\sqrt{1 - 2\|\UVcol_{i}\|^{2} \theta}}
	\right) \Uv_{i}
	\\
	&=&
	\biggl\{ 
		- \|\UVcol_{i}\|^{2}\theta - \frac{1}{2} \log(1 - 2 \theta \| \UVcol_{i} \|^{2}) 
	\biggr\}
	\Uv_{i}
\label{logEtei2m1i}
\end{EQA}
and 
\begin{EQA}
	\log \E\exp\bigl\{ \theta (-\gauss_{i}^{2}+1) \UVcol_{i} \UVcol_{i}^{\T} \bigr\}
	&=&
	\log \E \exp\left( \theta (-\gauss_{i}^{2}+1) \| \UVcol_{i} \|^{2} \right) \Uv_{i}
	\\
	&\le&- \|\UVcol_{i}\|^{2}\theta \Uv_{i}
		\\
	&\leq&
	\bigl\{ 
		- \|\UVcol_{i}\|^{2}\theta - \frac{1}{2} \log(1 - 2 \theta \| \UVcol_{i} \|^{2}) 
	\bigr\}
	\Uv_{i}.
\label{logEtei2m1ilower}
\end{EQA}
Then it follows by \eqref{matrix_chernoff}
\begin{EQA}
	&& \nquad
	\P\bigl( \| \Zv \|_{\oper} \geq \zq \bigr)
	\\
	& \leq &
	2 \inf_{\theta > 0 } \ex^{-\theta \zq} 
	\tr	\exp \biggl\{ 
		\sum_{i=1}^{n} \frac{\UVcol_{i}\UVcol_{i}^{\T}}{\|\UVcol_{i}\|^{2}} 
		\Bigl\{ - \|\UVcol_{i}\|^{2}\theta - \frac{1}{2} \log(1 - 2 \theta \| \UVcol_{i} \|^{2}) 
		\Bigr\} 
	\biggr\}.
	\qquad
\label{MastEqLog}
\end{EQA}	
Denote \( \etav = (\eta_{1},\dots,\eta_{n})^{\T}\), where 
\begin{EQA}
	\eta_{i}
	&=&  
	- \theta - \frac{\log(1 - 2\|\UVcol_{i}\|^{2} \theta)}{2\|\UVcol_{i}\|^{2}}.
\end{EQA}
The use of %\( \UV^{\T} \UV = \Id \) 
\eqref{jmu2v221mup} and \eqref{BoundRow} yields for \(\theta<(2\tdn^{2})^{-1}\) 
\begin{EQA}
	\eta_{i}
	&=&
	\frac{1}{2 \| \UVcol_{i} \|^{2}} 
	\bigl\{ 2 \theta \| \UVcol_{i} \|^{2} - \log(1 - 2 \theta \|\UVcol_{i}\|^{2}) \bigr\}
	\\
	& \leq &
	\frac{\bigl( 2 \theta \| \UVcol_{i} \|^{2} \bigr)^{2}}
		 {4 \| \UVcol_{i} \|^{2}(1 - 2 \theta \tdn^{2})}
	\leq 
	\frac{\theta^{2} \tdn^{2}}{(1 - 2 \theta \tdn^{2})} .
\label{ei12ui241m2}
\end{EQA}
%\begin{EQA}
%	\eta_{i} 
%	& \leq & 
%	2\|\UVcol_{i}\|^{4} \theta^{2}.
%\end{EQA}
Then by \eqref{MastEqLog} and \( \UV \UV^{\T} = \Id_{\dimp} \) using \( \mu = 2 \theta \tdn^{2} \)
\begin{EQA}
\label{}
	\P\bigl( 
		\| \Zv \|_{\oper}
		\geq  
		\zq
	\bigr)
	& \leq &
	2 \inf_{\theta > 0 } \ex^{-\theta \zq} 
	\tr	\exp \bigl\{ \UV \diag(\etav) \UV^{\T} \bigr\}
	\leq 
	2 \inf_{\theta > 0 } \ex^{-\theta \zq} 
	\tr	\exp \bigl\{ \| \etav \|_{\infty} \Id_{\dimp} \bigr\}
	\\
	& \leq &
	2 \dimp \inf_{\theta > 0} \exp\biggl\{ 
		- \theta \zq + \frac{\theta^{2} \tdn^{2}}{1 - 2 \theta \tdn^{2}} 
	\biggr\} 
	=
	2 \dimp \inf_{\mu > 0} \exp\biggl\{ 
		- \mu \frac{\zq}{2 \tdn^{2}} + \frac{\mu^{2} \tdn^{-2}}{1 - \mu} 
	\biggr\} .
\end{EQA}	
%Therefore, 
Lemma~\ref{Lmuvpxx} helps to bound 
%similarly to \eqref{x2xv4x12} 
for \( \xxp = \xx + \log (2\dimp) \) and \( \zq = 2 \tdn \xxp^{1/2} + 2 \tdn^{2} \xxp \) 
that  
\begin{EQA}
	\inf_{\mu > 0} \exp\biggl\{ 
		- \mu \frac{\zq}{2 \tdn^{2}} + \frac{\mu^{2} \tdn^{-2}}{1 - \mu} 
	\biggr\}
	& = &
	\inf_{\mu > 0} \biggl\{ 
		- \mu \bigl( \tdn^{-1} \xxp^{1/2} + \xxp \bigr) + \frac{\mu^{2} \tdn^{-2}}{4(1 - \mu)} 
	\biggr\}
	\leq 
	- \xxp \, .
\label{inft0tt21m2t2}
\end{EQA}
Therefore,
%\( t = 2 \tdn^{2}\bigl\{ \log (\dimp) + \xx \bigr\} 
%+ 2 \tdn \bigl\{ \log (\dimp) + \xx \bigr\}^{1/2} \)
\begin{EQA}
	\P\biggl( 
		\|\Zv \|_{\oper}
		\geq  
		2 \tdn \sqrt{ \xxp} 
		+ 2\tdn^{2} \xxp
	\biggr)
	& \leq &
	2 \dimp \, \ex^{ - \xxp}
	=
	\ex^{-\xx}
\end{EQA}
as required.
\end{proof}

\begin{proposition}[Deviation bound for matrix Gaussian sums]
\label{CUvepsB}
Let vectors \(\UVcol_{1},\dots,\UVcol_{n}\) in \( \R^{\dimp} \) satisfy
\begin{EQA}
\label{BoundRowmat}
 	\|\UVcol_{i}\|
 	&\leq &
 	\tdn
\end{EQA}
for a fixed constant \(\tdn\).
Let \(\gauss_{i}\) be independent standard Gaussian, \( i=1,\ldots,n \).  
Then for each vector \( \Bias = (b_{1},\ldots,b_{n})^{\T} \in \R^{n} \), 
the matrix \( \Zv_{1} \) with
\begin{EQA}
	\Zv_{1}
	& \eqdef &
	\sum_{i=1}^{n} \gauss_{i} b_{i} \UVcol_{i} \UVcol_{i}^{\T}
\end{EQA}
fulfills
\begin{EQA}
	\P\biggl( 
		\| \Zv_{1} \|_{\oper}
		\geq  
		\tdn^{2} \| \Bias \| \sqrt{2 \xx} 
	\biggr)
	& \leq &
	2 \ex^{-\xx}.
\end{EQA}
%\end{lemma}
\end{proposition}

\begin{proof}
As \( \gauss_{i} \) are i.i.d. standard normal and \( \E \ex^{a \gauss_{i}} = \ex^{a^{2}/2} \) for 
\( |a| < 1/2 \),
%\begin{EQA}
%	\log \E \exp\bigl( \theta \gauss_{i} b_{i} \| \UVcol_{i} \|^{2} \bigr)
%	&=&
%	\theta^{2} b_{i}^{2} \| \UVcol_{i} \|^{4}/2
%\label{logEexteibi}
%\end{EQA}
%and ,
it follows from the Master inequality and Lemma~\ref{LlEexiUv} 
\begin{EQA}
\label{matrix_chernoff1}
	\P\bigl( 
		\| \Zv_{1} \|_{\oper} \geq \zq
	\bigr)
	& \leq &
	2 \inf_{\theta > 0} \ex^{-\theta \zq} 
		\tr \exp \biggl\{ 
		\sum_{i=1}^{n} \log \E\exp( \theta \gauss_{i} b_{i} \UVcol_{i} \UVcol_{i}^{\T}) 
	\biggr\} 
	\\
	& \leq &
	2 \inf_{\theta > 0} \ex^{-\theta \zq} 
	\tr \exp \biggl\{ 
		\sum_{i=1}^{n} \frac{\theta^{2} b_{i}^{2} \| \UVcol_{i} \|^{4}}{2} \,\,
		\frac{\UVcol_{i}\UVcol_{i}^{\T}}{\|\UVcol_{i}\|^{2}} 
	\biggr\} .
\end{EQA}
Moreover, as \( \| \UVcol_{i} \| \leq \tdn \) and 
\( \Uv_{i} = \UVcol_{i} \UVcol_{i}^{\T} / \| \UVcol_{i} \|^{2} \) is a rank-one projector with 
\( \tr \Uv_{i} = 1 \), it holds
\begin{EQA}
	\tr \exp \biggl\{ 
		\frac{\theta^{2}}{2} \sum_{i=1}^{n} b_{i}^{2} \| \UVcol_{i} \|^{4} \Uv_{i}
	\biggr\}
	& \leq &
	\exp \tr\biggl( \frac{\theta^{2} \tdn^{4}}{2} \sum_{i=1}^{n} b_{i}^{2} \Uv_{i} \biggr)
	=
	\exp \frac{\theta^{2} \tdn^{4} \| \Bias \|^{2}}{2} \, .
\label{trexdPs422}
\end{EQA}
This implies
for \( \zq = \tdn^{2} \| \Bias \| \sqrt{2 \xx} \)
\begin{EQA}
	\P\bigl( 
		\| \Zv_{1} \|_{\oper} \geq \zq
	\bigr)
	& \leq &
	2 \inf_{\theta > 0} \exp\biggl(
		- \theta \zq + \frac{1}{2} \theta^{2} \tdn^{4} \| \Bias \|^{2} 
	\biggr)
	=
	2 \ex^{-\xx} 
\label{PZott0t12}
\end{EQA}
and the assertion follows.
\end{proof}

\subsection{Matrix valued quadratic forms}
Let \( \xiv = (\xi_{1},\ldots,\xi_{n})^{\T} \in \R^{n} \) be a Gaussian zero mean vector with the covariance matrix 
\( \Covm \) such that \( \| \Covm \|_{\oper} = \lambda_{\max}(\Covm) \leq \supAB \).
Let also \( \UV \) be a \( \dimp \times n \) matrix with columns 
\( \UVcol_{1},\ldots,\UVcol_{n} \in \R^{\dimp} \) such that 
\begin{EQ}[rcl]
	\tr(\UV \UV^{\T})
	= 
	\sum_{i=1}^{n} \| \UVcol_{i} \|^{2}
	& \leq &
	\dimw,
%	\\
%	\tr(\UV \UV^{\T})^{2}
%	& \leq &
%	\dimw,
	\\
	\max_{i} \| \UVcol_{i} \|^{2}
	& \leq &
	\tdn .	
\label{UVconds0}
\end{EQ}
A typical situation we have in mind is when 
\begin{EQA}
	\sum_{i=1}^{n} \UVcol_{i} \UVcol_{i}^{\T} 
	&=&
	\Id_{\dimp} .
%	\qquad
%	\| \UVcol_{i} \| \leq \tdn \, .
\label{sum1nuviuviTM}
\end{EQA}
Then \eqref{UVconds0} is satisfied with \( \dimw = \dimp \).
Moreover, it also holds \( \tr\bigl\{ (\UV \UV^{\T})^{2} \bigr\} = \tr(\UV \UV^{\T}) = \dimp \).
Consider the \( \dimp \times \dimp \) random matrix \( \BBB_{0} \) given by
\begin{EQA}
	\BBB_{0}
	& \eqdef &
	\sum_{i=1}^{n} \UVcol_{i} \UVcol_{i}^{\T} (\xi_{i}^{2} - \E \xi_{i}^{2}) .
\label{BBBdefsum1nevieEe}
\end{EQA}
In the case of \( \Covm = \Id_{n} \), Proposition~\ref{CUvepsUv} provides a bound 
for the operator norm of \( \BBB_{0} \). 
Below we extend this result to the case of a general matrix \( \Covm \) and 
establish similar bounds for quadratic forms of a non-centered vectors \( \xiv + \Bias \).
Also we evaluate the nuclear and Frobenius norms of this matrix.
We begin with establishing a bound on the Frobenius norm of \( \BBB_{0} \).
% trace \( \tr (\BBB_{0}^{2}) \). 

\begin{proposition}
\label{TmatrixQF}
Let vectors \( \UVcol_{i} \in \R^{\dimp} \) fulfill \eqref{UVconds0}.
Let also \( \xiv \sim \ND(0,\Covm) \) be a zero mean Gaussian vector 
\( \| \Covm \|_{\oper} \leq \supAB \).
Then the random matrix \( \BBB_{0} \) from \eqref{BBBdefsum1nevieEe} satisfies
\begin{EQA}
	\P\Bigl( 
		\tr(\BBB_{0}^{2}) > 2 \supAB \, \tdn^{2} \, \dimw \, (\xxn^{1/2} + \dPsis \, \xxn) 
	\Bigr)
	& \leq &
	\ex^{-\xx} ,
\label{PtrBBB22sAb}
\end{EQA}
where \( \dPsis \leq 1 \) and \( \xxn = \xx + \log(n) \).
% \( \dPsis \eqdef \max_{k} \| \UVcol_{k} \|_{\infty} \leq 1 \).
\end{proposition}

\begin{proof}
Denote \( \eta_{i} = \xi_{i}^{2} - \E \xi_{i}^{2} \) and 
\( c_{ij} = \UVcol_{i}^{\T} \UVcol_{j} \).
Then 
\begin{EQA}
	\tr (\BBB_{0}^{2})
	&=&
	\sum_{i=1}^{n} \sum_{j=1}^{n} 
		\eta_{i} \eta_{j} \tr \bigl( \UVcol_{i} \UVcol_{i}^{\T} \UVcol_{j} \UVcol_{j}^{\T} \bigr) 
	=
	\sum_{i=1}^{n} \sum_{j=1}^{n} c_{ij}^{2} \eta_{i} \eta_{j} \, .
\label{trBBB2e}
\end{EQA}
The \( n \times n \) matrix \( \Cv = (c_{ij}^{2}) \) is obviously symmetric positive. 
Therefore, one can represent it in the form
\( \Cv = \Uv \Mv \Uv^{\T} \) for a diagonal matrix \( \Mv = \diag(\mu_{1},\ldots,\mu_{n}) \) and 
an orthonormal \( n \times n \) matrix \( \Uv = (\uv_{1},\ldots,\uv_{n}) \) 
whose columns \( \uv_{k} \) are orthonormal vectors in \( \R^{n} \).
Therefore, for the vector \( \etav = (\eta_{1},\ldots,\eta_{n})^{\T} \), 
\begin{EQA}
	\tr (\BBB_{0}^{2})
	&=&
	\etav^{\T} \Cv \etav
	=
	\etav^{\T} \Uv \Mv \Uv^{\T} \etav
	=
	\sum_{k=1}^{n} \mu_{k} | \uv_{k}^{\T} \etav |^{2} .
\label{trBBB22LUz2}
\end{EQA}
Further, one can bound each \( \uv_{k}^{\T} \etav \) by the result of 
Corollary~\ref{CuvepsuvnG}: 
for any \( \xxn > 0 \)
\begin{EQA}
	\P\Bigl( 
		|\uv_{k}^{\T} \etav| 
		> 
		2 \supAB \bigl( \xxn^{1/2} + \| \uv_{k} \|_{\infty} \, \xxn \bigr) 
	\Bigr)
	& \leq &
	\ex^{-\xxn} .
\label{uvkzv112xxn}
\end{EQA} 
%with a probability at least \( 1 - \ex^{-\xxn} \).
The choice \( \xxn = \xx + \log(n) \) and \eqref{trBBB22LUz2} imply
\begin{EQA}
	\P\biggl( 
		\tr(\BBB_{0}^{2}) 
		> 
		\sum_{k=1}^{n} 2 \supAB \mu_{k} 
			\bigl(  \xxn^{1/2} + \| \uv_{k} \|_{\infty} \, \xxn \bigr) 
	\biggr)
	& \leq &
	\ex^{-\xx} .
\label{Psumk1plk2supAB}
\end{EQA}
Also by construction and \eqref{UVconds0}
\begin{EQA}
	\sum_{k=1}^{n} \mu_{k}
	&=&
	\tr (\Cv)
	=
	\sum_{i=1}^{n} c_{ii}^{2}
	=
	\sum_{i=1}^{n} \| \UVcol_{i} \|^{4}
	\leq 
	\tdn^{2} \tr \biggl( \sum_{i=1}^{n} \UVcol_{i} \UVcol_{i}^{\T} \biggr)
	=
	\tdn^{2} \dimw \, .
\label{sumk1nmukdPr}
\end{EQA}
The result \eqref{PtrBBB22sAb} uses a very rough bound \( \| \uv_{k} \|_{\infty} \leq \dPsis \)
for some constant \( \dPsis \leq 1 \).
In typical situation one can refine it to \( \| \uv_{k} \|_{\infty} \leq \CONST \tdn \).
\end{proof}

\medskip
Now we consider a slightly more general situation with a bias component. 
Given a bias vector \( \Bias \) in \( \R^{n} \),
a \( \dimp \times n \) matrix \( \UV \),
and a stochastic Gaussian zero mean vector \( \xiv \), define 
a random \( \dimp \times \dimp \) matrix 
\begin{EQA}
	\BBB_{1}
	& \eqdef &
	\UV \diag\bigl\{ 
		(\xiv + \Bias) \cdot (\xiv + \Bias) - \xiv \cdot \xiv 
	\bigr\} \UV^{\T} .
\label{BBBdefVm12VVVa}
\end{EQA}
The next result bounds the values \( \tr(\BBB_{1}) \) and \( \tr(\BBB_{1}^{2}) \).

\begin{proposition}
\label{TUVeps1UV}
Suppose that a Gaussian vector \( \xiv \sim \ND(0,\Covm) \) satisfies %\eqref{supepsdefa} 
\begin{EQ}[rcl]
	\supeps
	& \eqdef &
	\| \Covm - \Id_{n} \|_{\oper},
%	\\
%	\supepsi
%	& \eqdef & 
%	\max_{i} |\E \xi_{i}^{2} - 1|.
\label{supepsdefa}
\end{EQ}
Let also \( \UV \UV^{\T} \leq \Id_{\dimp} \) and
the vectors \( \UVcol_{i} \) - columns of \( \UV \) - satisfy for some 
\( \dimwf \leq \dimw \)
\begin{EQ}[ccccl]
	\tr(\UV \UV^{\T})
	&=& 
	\sum_{i=1}^{n} \| \UVcol_{i} \|^{2}
	& \leq &
	\dimw,
	\\
	\tr\bigl\{ (\UV \UV^{\T})^{2} \bigr\}
	&=&
	\sum_{i,j=1}^{n} \bigl| \UVcol_{i}^{\T} \UVcol_{j} \bigr|^{2}
	& \leq &
	\dimwf,
	\\
	&&
	\max_{i} \| \UVcol_{i} \|
	& \leq &
	\tdn .	
\label{UVconds}
\end{EQ}
%\eqref{UVUVTIdp}.
Then on a random set \( \Omega_{10}(\xx) \) with 
\( \P\bigl( \Omega_{10}(\xx) \bigr) \geq 1 - 2 \ex^{-\xx} \), it holds
\begin{EQA}
\label{trBBB10}
	\bigl\| \BBB_{1} \bigr\|_{\oper}
	& \leq &
	\| \Bias \|_{\infty}^{2} + \tdn^{2} \| \Bias \| \sqrt{2 \xx} 
%\label{trBBB4412}
\end{EQA}
and on a random set \( \Omega_{11}(\xx) \) with 
\( \P\bigl( \Omega_{11}(\xx) \bigr) \geq 1 - \ex^{-\xx} \), it holds
\begin{EQA}
\label{trBBB11}
	\bigl| \tr(\BBB_{1}) \bigr|
	& \leq &
	\dimw \, \| \Bias \|_{\infty}^{2}
	+ 4 \, \xx^{1/2} \, \tdn^{2} \, \| \Bias \| .
%\label{trBBB4412}
\end{EQA}
Further, \( \| \BBB_{1} \|_{\Frob} =	\sqrt{\tr(\BBB_{1}^{2})} \) fulfills
on a random set \( \Omega_{12}(\xx) \) with 
\( \P\bigl( \Omega_{12}(\xx) \bigr) \geq 1 - \ex^{-\xx} \) 
\begin{EQA}
	\| \BBB_{1} \|_{\Frob}
	& \leq &
	\errSi_{1}(\xx) ,
	\\
	\errSi_{1}(\xx)
	& \eqdef &
	\sqrt{\| \Bias \|_{\infty}^{4} \, \dimwf}
%	+ \sqrt{\tdn^{4} \| \Bias \|^{4} \dimp}
	+ 4 \tdn^{2} \| \Bias \| \bigl( 1 + \sqrt{\xx} \bigr) .
\label{errSiuBcuL}
\end{EQA}
\end{proposition}

\begin{proof}
We use the representation
\begin{EQA}
	\BBB_{1}
	&=&
	\UV \diag\bigl\{ (\xiv + \Bias) \cdot (\xiv + \Bias) \bigr\} \UV^{\T}
	- \UV \diag\bigl\{ \xiv \cdot \xiv \bigr\} \UV^{\T}
	\\
	&=&
	\underbrace{
		\UV \diag\bigl\{ \Bias \cdot \Bias \bigr\} \UV^{\T}
	}_{\BBB_{11}}
	+ \underbrace{
		2 \, \UV \diag\bigl\{ \xiv \cdot \Bias \bigr\} \UV^{\T} 
	}_{\BBB_{12}} \, .
\label{BBB12}
\end{EQA}
It obviously holds with \( \| \BBB_{1} \|_{\Frob} \eqdef \sqrt{\tr(\BBB_{1}^{2})} \)
\begin{EQ}[ccc]
	\bigl| \tr(\BBB_{1}) \bigr|
	& \leq &
	\bigl| \tr(\BBB_{11}) \bigr| + \bigl| \tr(\BBB_{12}) \bigr| \, ,
	\\
	\| \BBB_{1} \|_{\Frob}
	& \leq &
	\| \BBB_{11} \|_{\Frob} + \| \BBB_{12} \|_{\Frob} \, .
\label{errSiB12312}
\end{EQ}
We proceed with each \( \BBB_{1m} \) for \( m = 1,2 \) separately 
starting from \( \BBB_{11} \).

\paragraph{Bounds for \( \BBB_{11} \):}
The bias term \( \BBB_{11} \) can be estimated as follows:
\begin{EQ}[rcl]
	\bigl\| \BBB_{11} \bigr\|_{\oper}
	&=&
	\biggl\| \sum_{i=1}^{n} \UVcol_{i} \UVcol_{i}^{\T} \bias_{i}^{2} \biggr\|_{\oper}
	\leq 
	\| \Bias \|_{\infty}^{2} \bigl\| \UV \UV^{\T} \bigr\|_{\oper}
	\leq 
	\| \Bias \|_{\infty}^{2},
	\\
	\tr(\BBB_{11})
	&=&
	\tr \biggl( \sum_{i=1}^{n} \UVcol_{i} \UVcol_{i}^{\T} \bias_{i}^{2} \biggr)
	\leq 
	\| \Bias \|_{\infty}^{2} \tr\bigl( \UV \UV^{\T} \bigr)
	\leq 
	\dimw \, \| \Bias \|_{\infty}^{2},
	\\
	\tr(\BBB_{11}^{2})
	&=&
	\tr \biggl( \sum_{i=1}^{n} \UVcol_{i} \UVcol_{i}^{\T} \bias_{i}^{2} \biggr)^{2}
%	=
%	\sum_{i,j=1}^{n} c_{ij}^{2} \bias_{i}^{2} \bias_{j}^{2}
	\leq 
	\| \Bias \|_{\infty}^{4} \tr\bigl\{ \bigl( \UV \UV^{\T} \bigr)^{2} \bigr\}
	\leq 
	\dimwf  \,\| \Bias \|_{\infty}^{4} .
\label{trBBB3bi2bj2}
\end{EQ}
Another way of bounding the value \( \tr(\BBB_{11}^{2}) \) is based on the condition
\( \| \UVcol_{i} \| \leq \tdn \). 
Then for any unit vector \( \gammav \in \R^{\dimp} \), it holds 
\( |\gammav^{\T} \UVcol_{i}| \leq \tdn \) and
\begin{EQA}
	\gammav^{T} \BBB_{11} \gammav
	&=&
	\sum_{i=1}^{n} |\UVcol_{i}^{\T} \gammav|^{2} \bias_{i}^{2}
	\leq 
	\tdn^{2} \sum_{i=1}^{n} \bias_{i}^{2}
	=
	\tdn^{2} \| \Bias \|^{2} .
\label{gTBBB3gt}
\end{EQA}
Therefore, \( \| \BBB_{11} \|_{\oper} \leq \tdn^{2} \| \Bias \|^{2} \) and hence,
\begin{EQA}
	\tr(\BBB_{11}^{2})
	& \leq &
	\bigl( \tdn^{2} \| \Bias \|^{2} \bigr)^{2} \dimp
	=
	\tdn^{4} \| \Bias \|^{4} \dimp .
\label{trBBB3pdP4B4}
\end{EQA}
%cf. \eqref{uvbiasepsuvB}.
Note, however, that the bound \eqref{trBBB3bi2bj2} is typically more accurate:
the value \( \tdn^{2} \) is of order \( \dimw/n \) and 
\( \| \Bias \|^{2} \asymp n \| \Bias \|_{\infty}^{2} \), so that 
\( \tdn^{4} \| \Bias \|^{4} \dimp \gg \dimwf \| \Bias \|_{\infty}^{4} \) for 
\( \dimp \) large. 

\paragraph{Bounds for \( \BBB_{12} \):}
Proposition~\ref{CUvepsB} provides a bound for \( \BBB_{12} \) in the operator norm for 
standard Gaussian \( \xiv \):
\begin{EQA}
	\P\biggl( 
		\| \BBB_{12} \|_{\oper}
		\geq  
		\tdn^{2} \| \Bias \| \sqrt{2 \xx} 
	\biggr)
	& \leq &
	2 \ex^{-\xx}.
\end{EQA}
Now we bound \( \tr(\BBB_{12}) \) and \( \tr(\BBB_{12}^{2}) \).
By definition,
\begin{EQA}
	\tr(\BBB_{12})
	&=&
	2 \tr \biggl( \sum_{i=1}^{n} \UVcol_{i} \UVcol_{i}^{\T} \xi_{i} \bias_{i} \biggr)
	=
	2 \, \sum_{i=1}^{n} \| \UVcol_{i} \|^{2} \xi_{i} \bias_{i}
	=
	2 \uv^{\T} \xiv ,
\label{BBB4si1nbi}
\end{EQA}
where \( \uv \) is the vector in \( \R^{n} \) with the entries 
\( u_{i} = \| \UVcol_{i} \|^{2} \bias_{i} \).
As \( \Var(\xiv) = \Covm \) with \( \| \Covm \|_{\oper} \leq 2 \), 
\( \uv^{\T} \xiv \) is a Gaussian zero mean random variable whose variance satisfies
\begin{EQA}
	\Var(\uv^{\T} \xiv)
	& \leq &
	\uv^{\T} \Covm \uv
	\leq 
	2 \| \uv \|^{2}
	=
	2 \sum_{i=1}^{n} \| \UVcol_{i} \|^{4} \bias_{i}^{2} 
	\leq 
	2 \tdn^{4} \| \Bias \|^{2} .
\label{Varuvepsvr4dB2}
\end{EQA}
Here we have used \eqref{supepsdefa} and \( \| \UVcol_{i} \| \leq \tdn \).
Therefore, on a random set \( \Omega_{12}(\xx) \) with 
\( \P\bigl( \Omega_{12}(\xx) \bigr) \geq 1 - \ex^{-\xx} \),
\begin{EQA}
	\bigl| \tr(\BBB_{12}) \bigr|
	& \leq &
	2 \sqrt{2} \, \tdn^{2} \, \| \Bias \| \, \zq_{1}(\xx)
	\leq 
	4 \, \xx^{1/2} \tdn^{2} \, \| \Bias \| ,
\label{trBBB4u2s2}
\end{EQA}
where \( \zq_{1}(\xx) \leq \sqrt{2\xx} \) is given by 
\( \P(|\gauss| > \zq_{1}(\xx)) = \ex^{-\xx} \)
for a standard normal \( \gauss \).

It remains to bound \( \tr(\BBB_{12}^{2}) \).
Because of cross-dependence of the \( \xi_{i} \)'s, we cannot directly apply the result
of Proposition~\ref{CUvepsB}.
Instead we use the following representation:
%with \( c_{ij} = \UVcol_{i}^{\T} \UVcol_{j} \)
\begin{EQA}
	\tr(\BBB_{12}^{2})
	&=&
	4 \sum_{i,j=1}^{n} (\UVcol_{i}^{\T} \UVcol_{j})^{2} 
		\bias_{i} \bias_{j} \xi_{i} \xi_{j} \, .
\label{trBBB4sumij1n}
\end{EQA}
Denote by \( \Cv_{1} \) the \( n \times n \) matrix with the entries 
\( (\UVcol_{i}^{\T} \UVcol_{j})^{2} \bias_{i} \bias_{j} \) for \( i,j=1,\ldots,n \).
The use of \( \xiv = \Covm^{1/2} \gaussv \) with \( \Covm = \Var(\xiv) \) and a standard normal
\( \gaussv \in \R^{n} \) yields
\begin{EQA}
	\tr (\BBB_{12}^{2}) 
	&=& 
	4 \xiv^{\T} \Cv_{1} \xiv 
	=
	4 \gaussv^{\T} \Covm^{1/2} \Cv_{1} \Covm^{1/2} \gaussv 
	=
	4 \gaussv^{\T} \Cv_{2} \gaussv ,
\label{trBBB42eTC2Ce}
\end{EQA}
where \( \Cv_{2} = \Covm^{1/2} \Cv_{1} \Covm^{1/2} \).
Now the bound of Proposition~\ref{TexpbLGA} on Gaussian quadratic forms is well applicable.
It holds
\begin{EQA}
	\dimA(\Cv_{2})
	&=&
	\tr(\Cv_{2})
	\leq 
	2 \tr (\Cv_{1})
	=
	2 \sum_{i=1}^{n} \| \UVcol_{i} \|^{4} \bias_{i}^{2} 
%	\\
%	& \leq &
	\leq 
	2 \tdn^{4} \| \Bias \|^{2} .
\label{dimACv242}
\end{EQA}
Similarly for any unit vector \( \uv \in \R^{n} \), it holds by 
\( |\UVcol_{i}^{\T} \UVcol_{j}| \leq \tdn^{2} \)
\begin{EQA}
	\uv^{\T} \Cv_{1} \uv
	&=&
	\sum_{i,j=1}^{n} u_{i} u_{j} \bias_{i} \bias_{j} (\UVcol_{i}^{\T} \UVcol_{j})^{2}
	\leq 
	\tdn^{4} \biggl( \sum_{i=1}^{n} u_{i} \bias_{i} \biggr)^{2}
	\leq 
	\tdn^{4} \, \| \uv \|^{2} \, \| \Bias \|^{2}
	=
	\tdn^{4} \, \| \Bias \|^{2}
\label{uTCv1uvd4B2}
\end{EQA}
yielding \( \supA_{\max}(\Cv_{1}) \leq \tdn^{4} \, \| \Bias \|^{2} \) and 
\begin{EQA}
	\supA(\Cv_{2})
	& \eqdef &
	\supA_{\max}(\Cv_{2})
	\leq 
	2 \tdn^{4} \, \| \Bias \|^{2} .
\label{supACv21pe42}
\end{EQA}
Proposition~\ref{TexpbLGA} implies on a random set of probability at least \( 1 - \ex^{-\xx} \)
\begin{EQA}
	\sqrt{\tr(\BBB_{12}^{2})}
	& = &
	2 \sqrt{\gaussv^{\T} \Cv_{2} \gaussv}
	\\
	& \leq &
	2 \sqrt{\dimA(\Cv_{2})} + 2 \sqrt{2 \supA(\Cv_{2}) \xx}
	\\
	&=&
	2 \sqrt{2} \, \tdn^{2} \, \| \Bias \| \, \bigl( 1 + \sqrt{2 \xx} \bigr) 
	\\
	& \leq &
	4 \, \tdn^{2} \, \| \Bias \| (1 + \sqrt{\xx}) .
\label{strBBB421e12}
\end{EQA}
Putting all bounds together yields by \eqref{errSiB12312} the statements of the proposition.
%on a set \( \Omega(\xx) \) with 
%\( \P\bigl( \Omega(\xx) \bigr) \geq 1 - 3 \ex^{-\xx} \)
%\begin{EQA}
%	\sqrt{\tr(\BBB^{2})}
%	& \leq &
%	\sqrt{\tr(\BBB_{1}^{2})} + \sqrt{\tr(\BBB_{2}^{2})} 
%	+ \sqrt{\tr(\BBB_{3}^{2})} + \sqrt{\tr(\BBB_{4}^{2})} 
%	\\
%	& \leq &
%	\sqrt{4 \, \tdn^{2} \, \dimw \, \xxn}
%	+ \sqrt{\supepsi^{2} \dimwf} 
%	+ \sqrt{\| \Bias \|_{\infty}^{4} \dimwf}
%%	+ \sqrt{\tdn^{4} \| \Bias \|^{4} \dimp}
%	+ 4 \tdn^{2} \| \Bias \| \bigl( 1 + \sqrt{\xx} \bigr) .
%\label{errSiuBcu}
%\end{EQA}
%This implies \eqref{errSiuBcuL}. 
\end{proof}

\medskip
Now we evaluate the effect of presmoothing in the stochastic component. 
Let %\( \gaussv \in \R^{\dimp} \) be standard normal while 
\( \xiv \) be normal zero mean vector 
in \( \R^{\dimp} \) with a covariance matrix \( \Covm \) satisfying \eqref{supepsdefa}.
The goal is to bound the values \( \tr(\BBB_{2}) \) and \( \tr(\BBB_{2}^{2}) \) for the difference
\begin{EQA}
	\BBB_{2}
	& \eqdef &
	\UV \diag\bigl\{ \xiv \cdot \xiv - \E (\xiv \cdot \xiv) \bigr\} \UV^{\T} . 
\label{BBB2defxxr}
\end{EQA}

\begin{proposition}
\label{TUVeps2UV}
Suppose that %\( \gaussv \sim \ND(0,\Id_{\dimp}) \) while
\( \xiv \sim \ND(0,\Covm) \) with \( \Covm \) satisfying \eqref{supepsdefa}.
Let also \( \UV \UV^{\T} \leq \Id_{\dimp} \) and
the vectors \( \UVcol_{i} \) - columns of \( \UV \) - satisfy \eqref{UVconds} for some 
\( \dimwf \leq \dimw \).
Then on a random set \( \Omega_{22}(\xx) \) with 
\( \P\bigl( \Omega_{22}(\xx) \bigr) \geq 1 - \ex^{-\xx} \), it holds
for \( \| \BBB_{2} \|_{\Frob} = \sqrt{\tr(\BBB_{2}^{2})} \)  
\begin{EQA}
	\| \BBB_{2} \|_{\Frob}
	& \leq &
	\errSi_{2}(\xx) 
	\eqdef 
	2 \, \sqrt{\xxn \, \dimw \, \tdn^{2} } \, . %+ \sqrt{\supepsi^{2} \, \dimwf}\,  .
\label{errSiuBcuL}
\end{EQA}
Moreover, on a random set \( \Omega_{21}(\xx) \) with 
\( \P\bigl( \Omega_{21}(\xx) \bigr) \geq 1 - 2 \ex^{-\xx} \), it holds 
\begin{EQA}
\label{trBBB21}
	\bigl| \tr(\BBB_{2}) \bigr| 
	& \leq &
	4 \sqrt{\xx \, \dimw \, \tdn^{2}} + 4 \, \xx \, \tdn^{2} . %+ \supepsi \, \dimw  .
%\label{trBBB4412}
\end{EQA}
\end{proposition}

\begin{proof}
%\paragraph{Bounds for \( \tr(\BBB_{21}^{2}) \):}
%This part is most involved because \( \BBB_{21} \) is a matrix valued quadratic form of \( \xiv \).
%
By \eqref{supepsdefa}, the covariance matrix \( \Covm = \Var(\xiv) \) fulfills 
%In general, there exists a fixed constant \( \supeps \leq 1/2 \) such that
\begin{EQA}
	\| \Covm \|_{\oper}
	& \leq &
	1 + \supeps \leq 2.
\label{lmaCverAB}
\end{EQA}
Now, %we bound \( \tr(\BBB_{21}^{2}) \).
by Proposition~\ref{TmatrixQF}, for \( \xxn = \xx + \log(n) \) and \( \dPsis \leq 1 \),
it holds on a random set \( \Omega'_{22}(\xx) \) with 
\( \P\bigl( \Omega'_{22}(\xx) \bigr) \geq 1 - \ex^{-\xx} \)
\begin{EQA}
	\tr(\BBB_{2}^{2})
	& \leq &
	2 (1 + \supeps) \, \tdn^{2} \, \dimw \, (\xxn^{1/2} + \dPsis \, \xxn) .
\label{trBBB12supAB}
\end{EQA}
Here \( \dPsis \leq 1 \), usually \( \dPsis \ll 1 \), and \( \supeps \leq 1 \), 
so we simplify the bound to 
\begin{EQA}
	\tr(\BBB_{2}^{2})
	& \leq &
	4 \, \xxn \, \dimw \, \tdn^{2} .
\label{trBBB12supABn}
\end{EQA}

%\paragraph{Bounds for \( \tr(\BBB_{22}^{2}) \):}
%Here we can use the same arguments with \( \supeps = 0 \) yielding on 
%a random set \( \Omega''_{22}(\xx) \) with 
%\( \P\bigl( \Omega''_{22}(\xx) \bigr) \geq 1 - \ex^{-\xx} \)
%\begin{EQA}
%%	\bigl| \tr(\BBB_{22}) \bigr|
%%	& \leq &
%%	2 \xx^{1/2} \, \dimw \, \tdn^{2} + 2 \, \xx \, \tdn^{2} ,
%%	\\
%	\tr(\BBB_{22}^{2})
%	& \leq &
%	2 \, \tdn^{2} \, \dimw \, \xxn .
%\label{trBBB22supABn}
%\end{EQA}
%
%
%
%\paragraph{Bounds for \( \tr(\BBB_{23}^{2}) \):}

Now we bound the trace \( \tr(\BBB_{2}) \).
By definition, it holds for columns \( \UVcol_{i} \in \R^{\dimq} \) of \( \UV \)
\begin{EQA}
	\tr(\BBB_{2})
	&=&
	\sum_{i=1}^{n} \tr \bigl( \UVcol_{i} \UVcol_{i}^{\T} \bigr) 
		(\xi_{i}^{2} - \E \xi_{i}^{2})
	=
	\sum_{i=1}^{n} \| \UVcol_{i} \|^{2}	(\xi_{i}^{2} - \E \xi_{i}^{2}),
\label{trBBB1i1nE}
\end{EQA}
and by Corollary~\ref{CuvepsuvnG}, it holds on a random set \( \Omega_{21}(\xx) \) with
\( \P\bigl( \Omega_{21}(\xx) \bigr) \geq 1 - 2 \ex^{-\xx} \)
\begin{EQA}
	\bigl| \tr(\BBB_{2}) \bigr|
	& \leq &
	4 \xx^{1/2} \biggl( \sum_{i=1}^{n} \| \UVcol_{i} \|^{4} \biggr)^{1/2} 
	+ 4 \xx \max_{i} \| \UVcol_{i} \|^{2} .
\label{1trBBB11}
\end{EQA}
This implies in view of \( \| \UVcol_{i} \| \leq \tdn \) and \eqref{UVconds}
\begin{EQA}
	\sum_{i=1}^{n} \| \UVcol_{i} \|^{4}
	& \leq &
	\max_{i} \| \UVcol_{i} \|^{2} \sum_{i=1}^{n} \| \UVcol_{i} \|^{2}
%	\leq 
%	\tdn^{2} \tr\biggl( \sum_{i=1}^{n} \UVcol_{i} \UVcol_{i}^{\T} \biggr)
	\leq 
	\dimw \, \tdn^{2} .
\label{sumi1nui4}
\end{EQA}
Therefore, it holds on \( \Omega_{21}(\xx) \)
\begin{EQA}
	\bigl| \tr(\BBB_{2}) \bigr|
	& \leq &
	4 \sqrt{\xx \, \dimw \, \tdn^{2}} + 4 \, \xx \, \tdn^{2} 
\label{trBBB1u}
\end{EQA}
as required.
\end{proof}

Finally we bound the deterministic term 
\begin{EQA}
	\BBB_{3}
	& \eqdef &
	\UV \diag\bigl\{ \E (\xiv \cdot \xiv) - \Id_{n} \bigr\} \UV^{\T} .
\label{BBB3def1n}
\end{EQA}

\begin{proposition}
\label{TUVeps22UV}
Suppose that %\( \gaussv \sim \ND(0,\Id_{\dimp}) \) while
\( \xiv \sim \ND(0,\Covm) \) with 
% \( \Covm \) satisfying \eqref{supepsdefa}.
\begin{EQA}
	\supepsi
	& \eqdef & 
	\max_{i} |\E \xi_{i}^{2} - 1|.
\label{supepsdeff}
\end{EQA}
Let also \( \UV \UV^{\T} \leq \Id_{\dimp} \) and
the vectors \( \UVcol_{i} \) - columns of \( \UV \) - satisfy \eqref{UVconds} for some 
\( \dimwf \leq \dimw \).
Then it holds for the matrix \( \BBB_{3} \) from \eqref{BBB3def1n} 
\begin{EQA}
\label{trBBB32}
%	\bigl| \tr(\BBB_{21}) \bigr| + 
	\bigl| \tr(\BBB_{3}) \bigr|
	& \leq &
	\supepsi \, \dimw ,
	\\
	\tr(\BBB_{3}^{2}) 
	& \leq &
	\errSi_{3}(\xx) \eqdef \supepsi^{2} \, \dimwf   .
%\label{trBBB4412}
\end{EQA}
\end{proposition}

\begin{proof}
By direct calculus, it holds
in view of \eqref{UVconds}
%\( \UV \UV^{\T} = \Id_{\dimp} \),
%\( \| \UVcol_{i} \|\leq \tdn \), 
and \( |\E \xi_{i}^{2} - 1| \leq \supepsi \)
\begin{EQA}
	\tr(\BBB_{3}^{2})
	&=&
	\tr \biggl( \sum_{i=1}^{n} \UVcol_{i} \UVcol_{i}^{\T} (\E \xi_{i}^{2} - 1) \biggr)^{2}
%	\\
%	& \leq &
	\leq 
	\supepsi^{2} \, \tr\bigl\{ \bigl( \UV \UV^{\T} \bigr)^{2} \bigr\}
%	=
%	\supepsi^{2} \tr (\Id_{\dimp}) 
	\leq 
	\supepsi^{2} \, \dimwf ;
\label{trBBB2cij2}
\end{EQA}
Further,
\begin{EQA}
	\bigl| \tr(\BBB_{3}) \bigr|
	& \leq &
	\tr \biggl( \sum_{i=1}^{n} \UVcol_{i} \UVcol_{i}^{\T} \bigl| \E \xi_{i}^{2} - 1 \bigr| 
	\biggr)
	\leq 
	\supepsi \tr (\UV \UV^{\T}) 
	\leq 
	\supepsi \, \dimw .
\label{trBBB23}
\end{EQA}
This yields the assertion.
\end{proof}

Now we present a bound in the operator norm of the matrix \( \BBB \)
with 
\begin{EQA}
	\BBB 
	& \eqdef &
	\UV \diag\bigl\{ 
		(\xiv + \Bias) \cdot (\xiv + \Bias) - \Id_{n} 
	\bigr\} \UV^{\T} \, .
\label{BBBdef145}
\end{EQA}

%\ifims{ % full
\begin{proposition}
\label{TboundBBBinoper}
Assume the conditions of Proposition~\ref{TUVeps1UV},
and let the rows \( \SPiS_{i}^{\T} \) of \( \SPiS \eqdef \Sigma^{-1/2} \Pi \Sigma^{1/2} \) satisfies 
\( \| \SPiS_{i} \| \leq \tdPi \).
Then on a random set \( \Omega_{\oper}(\xx) \) with 
\( \P\bigl( \Omega_{\oper}(\xx) \bigr) \geq 1 - 6 \ex^{-\xx} \),
it holds with \( \xxn = \xx + \log(n) \) and \( \xx_{\dimp} = \xx + 2 \log(\dimp) \)
for \( \BBB \) from \eqref{BBBdef145}
\begin{EQA}
	\| \BBB \|_{\oper}
	& \leq &
	\errSi_{\oper}(\xx),
	\\
	\errSi_{\oper}(\xx)
	& \eqdef &
	\| \Bias \|_{\infty}^{2} + \tdn^{2} \| \Bias \| \sqrt{2 \xx} 
	+ 2 \tdn \xx_{\dimp}^{1/2} + 2 \tdn^{2} \xx_{\dimp}
	+ 2 \tdPi \xxn + \tdPi^{2} \xxn .
\label{BBBoperrSopx}
\end{EQA}
\end{proposition}

\begin{proof}
We use the following decomposition:
%a slightly changed version of the decomposition \eqref{BBBdefVm12VVV}:
\begin{EQA}[rclccl]
	\BBB 
%	&=&
%	\UV \diag\bigl\{ 
%		(\xiv + \Bias) \cdot (\xiv + \Bias) - \Id_{n} 
%	\bigr\} \UV^{\T}
%	\\
	&=&
	\UV \diag\bigl\{ 
		(\xiv + \Bias) \cdot (\xiv + \Bias) - \xiv \cdot \xiv 
	\bigr\} \UV^{\T} 
	& \qquad
	& \eqdef &
	\BBB_{1}
	\\
	&& \quad
	+ \,\, \UV \diag\bigl\{ 
		\xiv \cdot \xiv - \gaussv \cdot \gaussv
	\bigr\} \UV^{\T}
	& \qquad
	& \eqdef &
	\BBB_{4}
	\\
	&& \quad
	+ \,\, \UV \diag\bigl\{ 
	 	\gaussv \cdot \gaussv - \Id_{n} 
	\bigr\} \UV^{\T}
	& \qquad
	& \eqdef &
	\BBB_{5}
%	\\
%	&=&
%	\BBB_{1} + \BBB_{4} + \BBB_{5} .
\label{BBBdefggVVV}
\end{EQA}
Here \( \gaussv \eqdef \Sigma^{-1/2} \epsv \) is a standard Gaussian vector.
Obviously
\begin{EQA}
	\| \BBB \|_{\oper}
	& \leq &
	\| \BBB_{1} \|_{\oper} + \| \BBB_{4} \|_{\oper} + \| \BBB_{5} \|_{\oper} \, .
\label{BBBopub}
\end{EQA}
The value \( \| \BBB_{1} \|_{\oper} \) is already evaluated in \eqref{trBBB10} of 
Proposition~\ref{TUVeps1UV}:
\begin{EQA}
	\| \BBB_{1} \|_{\oper}
	& \leq &
	\| \Bias \|_{\infty}^{2} + \tdn^{2} \| \Bias \| \sqrt{2 \xx} 
\label{BBB1op2x}
\end{EQA}
on a random set \( \Omega_{10}(\xx) \) with 
\( \P\bigl( \Omega_{10}(\xx) \bigr) \geq 1 - 2 \ex^{-\xx} \).
The matrix \( \BBB_{5} \) can be bounded by a version of matrix Bernstein inequality 
\eqref{Zvop2tdxx2p} in Proposition~\ref{CUvepsUv}:
one a set \( \Omega_{5}(\xx) \) with \( \P\bigl( \Omega_{5}(\xx) \bigr) \geq 1 - 2 \ex^{-\xx} \)
\begin{EQA}
	\| \BBB_{5} \|_{\oper}
	& \leq & 
	2 \tdn \sqrt{ \xxp} + 2 \tdn^{2} \xxp .
\label{BBB5vop2tdxx2p}
\end{EQA}
It remains to bound the value \( \| \BBB_{4} \|_{\oper} \).
By definition, with \( \SPiS = \Sigma^{-1/2} \Pi \Sigma^{1/2} \) 
\begin{EQA}
	\xiv 
	&=& 
	\Sigma^{-1/2} (\epsv - \Pi \epsv)
	=
	\gaussv - \SPiS \gaussv .
\label{xivdefSm12gXg}
\end{EQA} 
This obviously implies 
\begin{EQA}
	\xiv \cdot \xiv - \gaussv \cdot \gaussv
	&=&
	(\SPiS \gaussv) \cdot (\SPiS \gaussv) - 2 (\SPiS \gaussv) \cdot \gaussv
\label{xicxigcgXX}
\end{EQA}
and
\begin{EQA}
	\| \BBB_{4} \|_{\oper}
	& \leq &
	\bigl\| \UV \diag\bigl\{ (\SPiS \gaussv) \cdot (\SPiS \gaussv) \bigr\} \UV^{\T} \bigr\|_{\oper}
	+ 2 \bigl\| \UV \diag\bigl\{ (\SPiS \gaussv) \cdot \gaussv \bigr\} \UV^{\T} \bigr\|_{\oper}
\label{BBB4ople}
\end{EQA}
The condition \( \UV \UV^{\T} \leq \Id_{\dimp} \) helps to bound
\begin{EQA}
	\bigl\| \UV \diag\bigl\{ (\SPiS \gaussv) \cdot (\SPiS \gaussv) \bigr\} \UV^{\T} \bigr\|_{\oper}
	& \leq &
	\| \SPiS \gaussv \|_{\infty}^{2} \,\, \| \UV \UV^{\T} \|_{\oper}
	\leq 
	\| \SPiS \gaussv \|_{\infty}^{2} .
\label{UVXigcXg}
\end{EQA}
Similarly 
\begin{EQA}
	\bigl\| \UV \diag\bigl\{ (\SPiS \gaussv) \cdot \gaussv \bigr\} \UV^{\T} \bigr\|_{\oper}
	& \leq &
	\| \SPiS \gaussv \|_{\infty} \,\, \| \gaussv \|_{\infty} \, \| \UV \UV^{\T} \|_{\oper}
	\leq 
	\| \SPiS \gaussv \|_{\infty} \,\, \| \gaussv \|_{\infty} .
\label{UVXigcgi}
\end{EQA}
It is well known that the sup-norm of a standard Gaussian vector \( \gaussv \) can be bounded as
\begin{EQA}
	\| \gaussv \|_{\infty}
	& \leq &
	\sqrt{2 \xxn}
\label{gausinfln}
\end{EQA}
with \( \xxn = \xx + \log(n) \)
on a set of probability \( 1 - \ex^{-\xx} \). 
Further, if each row \( \SPiS_{i}^{\T} \) of \( \SPiS \) satisfies 
\( \| \SPiS_{i} \| \leq \tdPi \), then
the scalar product 
\( \SPiS_{i}^{\T} \gaussv \) is a normal zero mean r.v. with the variance 
\begin{EQA}
	\Var\bigl( \SPiS_{i}^{\T} \gaussv \bigr) 
	&=&
	\| \SPiS_{i} \|^{2}
	\leq 
	\tdPi^{2} 
\label{VaruvTSPi}
\end{EQA}
and 
\begin{EQA}
	\P\bigl( |\SPiS_{i}^{\T} \gaussv| > \tdPi \zq_{1}(\xx) \bigr)
	& \leq &
	\ex^{-\xx}
\label{PSPiSTgatd1}
\end{EQA}
with \( \zq_{1}(\xx) \leq \sqrt{2 \xx} \)
yielding 
\begin{EQA}
	\P\bigl( \| \SPiS \gaussv \|_{\infty} > \tdPi \sqrt{2 \xxn} \bigr)
	& \leq &
	\ex^{-\xx} .
\label{PSPiSTgatdinf}
\end{EQA}
Summing together results in the bound
\begin{EQA}
	\| \BBB_{4} \|_{\oper}
	& \leq &
	2 \tdPi \xxn + \tdPi^{2} \xxn
\label{BBB4oper2}
\end{EQA}
on a set \( \Omega_{4}(\xx) \) with \( \P\bigl( \Omega_{4}(\xx) \bigr) \geq 1 - 2 \ex^{-\xx} \).
\end{proof}
%}{ % annals
%}

\ifbook{}{
\subsection{Empirical covariance matrix}
Let \( \epsv_{i} \) be independent centered random vectors in \( \R^{\dimp} \).
Consider the empirical covariance matrix 
\begin{EQA}
	\hat{\Varxi}
	& \eqdef &
	\sum_{i} \epsv_{i} \epsv_{i}^{\T} .
\label{VPrdefT}
\end{EQA}
The interesting question is how well this matrix stabilizes around its expectation
\begin{EQA}
	\Varxi
	& \eqdef &
	\E \hat{\Varxi}
	=
	\sum_{i=1}^{n} \Var(\epsv_{i})
	=
	\Var\biggl( \sum_{i=1}^{n} \epsv_{i} \biggr).
\label{VPEhVPi1n}
\end{EQA}
We implicitly assume that \( \Varxi \) is large and hence each \( \Varxi^{-1/2}\epsv_{i} \) is small.
Below we focus on the relative difference 
\begin{EQA}
	\Zv
	& \eqdef &
	\Varxi^{-1/2} (\hat{\Varxi} - \Varxi) \Varxi^{-1/2}
	=
	\Varxi^{-1/2} \, \hat{\Varxi} \, \Varxi^{-1/2} - \Id_{\dimp} \, .
\label{ZvVP12VPr}
\end{EQA}
This allows to reduce the study to the case \( \Varxi \equiv \Id_{\dimp} \) considered below.
Note that each vector \( \epsv_{i} \) is replaced by \( \Varxi^{-1/2} \epsv_{i} \).
For bounding \( \| \Zv \|_{\oper} \), one can apply
the Bernstein inequality provided that each vector \( \epsv_{i} \) is bounded 
by a fixed small constant \( \rei \):
\begin{EQA}
	\| \epsv \|_{\infty}
	=
	\max_{i \leq n}
	\| \epsv_{i} \|
	& \leq &
	\rei .
\label{epsvi2rr}
\end{EQA}
Then \( \Sv_{i} = \epsv_{i} \epsv_{i}^{\T} - \E \epsv_{i} \epsv_{i}^{\T} \) obviously fulfills
\begin{EQA}
	\| \Sv_{i} \|_{\oper}
	=
	\| \epsv_{i} \epsv_{i}^{\T} - \E \epsv_{i} \epsv_{i}^{\T} \|_{\oper}
	\leq 
	\| \epsv_{i} \|^{2}
	& \leq &
	\rei^{2} .
\label{epsiepsiToprr2}
\end{EQA}
Define
\begin{EQA}
	\vpi^{2}
	&=&
	\left\|\sum_{i=1}^{n} \E \Sv_{i}^{2} \right\|_{\oper} 
	=
	\left\|
		\sum_{i=1}^{n} 
			\Bigl\{ \E \bigl( \epsv_{i} \epsv_{i}^{\T} \bigr)^{2}
				- \bigl( \E \epsv_{i} \epsv_{i}^{\T} \bigr)^{2} 
			\Bigr\}
	\right\|_{\oper} \, .
\label{vp2sESvi2}
\end{EQA}
A rough bound on \( \vpi \) can be obtained by using again \( \| \epsv_{i} \| \leq \rei \) and 
\( \Varxi = \Id_{\dimp} \):
\begin{EQA}
	\vpi^{2}
	& \leq &
	\rei^{2} \left\|\sum_{i=1}^{n} \E \epsv_{i} \epsv_{i}^{\T} \right\|_{\oper}
	=
	\rei^{2}
\label{vp2rr2i1nr2}
\end{EQA}
Now the result of Proposition~\ref{BernSqTh} implies for \( \Zv = \hat{\Varxi} - \Id_{\dimp} \)
 with \( \Sv_{i} = \epsv_{i} \epsv_{i}^{\T} - \E \epsv_{i} \epsv_{i}^{\T} \)
\begin{EQA}
	\P\{\| \hat{\Varxi} - \Id_{\dimp} \|_{\oper} \geq \vpi \, \zq \} 
	& \leq & 
	2 \dimp \exp \left( \frac{-\zq^{2}/2}{1 + \rei^{2} \vpi^{-1} \zq /3}\right).
\label{PVprIpopp}
\end{EQA}
The bound is particularly informative if 
the ratio \( \rei^{2}/\vpi \) is small. 
The use of the rough upper bound \eqref{vp2rr2i1nr2}
yields for each \( \zq > 0 \)
\begin{EQA}
	\P\{\| \hat{\Varxi} - \Id_{\dimp} \|_{\oper} \geq \rei \, \zq \} 
	& \leq & 
	2 \dimp \exp \left( \frac{- \zq^{2}/2}{1 + \zq \rei /3}\right).
\label{PVPrIpt221t3}
\end{EQA} 
One can pick up \( \zq \approx \sqrt{2 \xx \log (2\dimp)} \) to ensure a sensible deviation bound
about \( \ex^{- \xx} \) for moderate values of \( \xx \).

\begin{theorem}
\label{Tmatrdevboundteps}
Let random vectors \( \epsv_{1},\ldots,\epsv_{n} \) in \( \R^{\dimp} \) 
be independent zero mean, \( \Varxi \) is given by \eqref{VPEhVPi1n}, and 
\( \tilde{\epsv}_{i} \eqdef \Varxi^{-1/2} \epsv_{i} \) fulfill
\begin{EQA}
	\| \tilde{\epsv}_{i} \|
	=
	\| \Varxi^{-1/2} \epsv_{i} \|
	& \leq &
	\rei \text{ a.s. },
	i=1,\ldots,n,
\label{•}
\end{EQA}
for a constant \( \rei < \infty \).
Then with \( \vpi^{2} \) defined by
\begin{EQA}
	\vpi^{2}
	& \eqdef &
	\bigl\| \E \bigl(\Varxi^{-1/2} \hat{\Varxi} \Varxi^{-1/2} \bigr)^{2} \bigr\|_{\oper}
	=
	\left\| 
		\sum_{i=1}^{n} 
			\Bigl\{ \E \bigl( \tilde{\epsv}_{i} \tilde{\epsv}_{i}^{\T} \bigr)^{2}
				- \bigl( \E \tilde{\epsv}_{i} \tilde{\epsv}_{i}^{\T} \bigr)^{2} 
			\Bigr\}
	\right\|_{\oper} \, ,
\label{vpi2VhVV12op}
\end{EQA}
it holds \( \vpi^{2} \leq \rei^{2} \) and for any \( \zq \geq 0 \)
\begin{EQA}
	\P\Bigl\{
		\| \Varxi^{-1/2} \, \hat{\Varxi} \, \Varxi^{-1/2} - \Id_{\dimp} \|_{\oper} \geq \vpi \, \zq 
	\Bigr\} 
	& \leq & 
	2 \dimp \exp \left( \frac{-\zq^{2}/2}{1 + \rei^{2} \vpi^{-1} \zq /3} \right),
	\\
	\P\Bigl\{ 
		\| \Varxi^{-1/2} \, \hat{\Varxi} \, \Varxi^{-1/2} - \Id_{\dimp} \|_{\oper} \geq \rei \, \zq 
	\Bigr\} 
	& \leq & 
	2 \dimp \exp \left( \frac{-\zq^{2}/2}{1 + \zq \rei /3} \right).
\label{PVprIpoppfin}
\end{EQA}
\end{theorem}
The result can be easily extended to the case when each \( \epsv_{i} \) is not bounded 
but can be bounded with a high probability.

\begin{theorem}
\label{Tmatrepsvunbounded}
Let for some \( \xx > 0 \) there exists a constant \( \rei(\xx) \) such that
\begin{EQA}
	\P\bigl( \| \Varxi^{-1/2} \epsv \|_{\infty} > \rei(\xx) \bigr)
	& \leq &
	\ex^{- \xx} .
\label{Pmineirx}
\end{EQA}
Then with \( \vpi \) from \eqref{vp2sESvi2}
\begin{EQA}
	\P\{\| \Varxi^{-1/2} \, \hat{\Varxi} \, \Varxi^{-1/2} - \Id_{\dimp} \|_{\oper} \geq \vpi \, \zq \} 
	& \leq & 
	2\dimp \exp \left( \frac{-\zq^{2}/2}{1 + \rei^{2}(\xx) \vpi^{-1} \zq /3}\right) + \ex^{-\xx}.
\label{PVprIpoppxx}
\end{EQA}
\end{theorem}

As a practical corollary, one deduces for moderate \( \xx \) that
\( \zq(\xx) \approx \sqrt{(2 + \alpha) \xx \log (2\dimp)} \) for some small \( \alpha > 0 \) ensures
\begin{EQA}
	\P\Bigl\{ 
		\| \Varxi^{-1/2} \, \hat{\Varxi} \, \Varxi^{-1/2} - \Id_{\dimp} \|_{\oper} 
		\geq \vpi \, \zq(\xx) 
	\Bigr\} 
	& \leq & 
	2 \ex^{-\xx} .
\label{PVPrIdpvz2x}
\end{EQA}

\tobedone{The i.i.d. case}
}
%%%%%%%%%%%%%%%%%%%%%%%%%%%%%%%%%%%%%%%%%%%%%%%%%%%%%%%%%%%%%%%%

%%%%%%%%%%%%%%%%%%%%%%%%%%%%%%%%%%%%%%%%%%%%%%%%%%%%%%%%%%%%%%%%
% !TEX root = script2014.tex

\section{Gaussian comparison via KL-divergence and Pinsker's inequality}
Suppose that two \( \dimp \)-dimensional zero mean Gaussian vectors 
\( \xiv \sim \ND(0,\Varxi) \) and \( \bxiv \sim \ND(0,\Varxib) \) are given.
Let also \( T \) map \( \R^{\dimp} \) to \( \R^{\Meta} \) and \( \Xv = T(\xiv) \)
and \( \Yv = T(\bxiv) \).
We aim to bound the distance between distributions of \( \Xv \) and \( \Yv \)
under the conditions
\begin{EQA}
	\| \Varxi^{-1/2} \Varxib \Varxi^{-1/2} - \Id_{\dimp} \|_{\oper}
	& \leq &
	\rd
	\leq 
	1/2,
	\quad
	\\
	\tr\bigl( \Varxi^{-1/2} \Varxib \Varxi^{-1/2} - \Id_{\dimp} \bigr)^{2}
	& \leq &
%	\leq 
	\errSi^{2}
\label{dSiSic2}
\end{EQA}
for some \( \rd \leq 1/2 \) and \( \errSi \geq 0 \).
%Obviously \( 2 \errSi^{2} \leq \dimp \rd^{2} \).
The next lemma bounds from above the Kullback-Leibler divergence between 
two normal distributions.

\begin{lemma}
\label{KullbTVd}
Let \( \P_{0} = \ND(\bvs,\Varxi) \) and \( \P_{1} = \ND(\bvb,\Varxib) \) for
some non-degenerated matrices \( \Varxi \) and \( \Varxib \). 
If 
\begin{EQA}
	\| \Varxi^{-1/2} \Varxib \Varxi^{-1/2} - \Id_{\dimp} \|_{\oper} 
	& \leq &
%	\rd \leq 
	1/2 ,
	\\
	\tr \Bigl\{ \bigl( \Varxi^{-1/2} \Varxib \Varxi^{-1/2} - \Id_{\dimp} \bigr)^{2} \Bigr\}
	& \leq &
	\errSi^{2},
\label{DPcbIdeps}
\end{EQA}
then 
\begin{EQA}
	\kullb(\P_{0},\P_{1})
	&=&
	- \E_{0} \log \frac{d\P_{1}}{d\P_{0}}
%	\\
%	& \leq &
 	\leq 
	\frac{\errSi^{2}}{2} + \frac{1}{2} (\bvs - \bvb)^{\T} \Varxib (\bvs - \bvb)  .
%	\leq  
%	\rd^{2} \, \dimp + (1 + \rd) \| \betav \|^{2} .
\label{kullbP0P1}
\end{EQA}
For any measurable set \( A \subset \R^{\dimp} \), it holds 
%with \( \gammav \sim \ND(0,\Id_{\dimp}) \)
\begin{EQA}[c]
	\bigl| \P_{0}(A) - \P_{1}(A) \bigr|
%	=
%	\bigl| \P\bigl( \gammav \in A \bigr)
%	- \P\bigl( \DD (\gammav - \betav) \in A \bigr) \bigr| 
	\leq 
	\sqrt{\kullb(\P_{0},\P_{1}) / 2}.
\end{EQA}
\end{lemma}

\begin{proof}
The  change of variables \( \uv = \Varxi^{-1/2} (\xv - \bvs) \) reduces the general case 
to the situation when \( \P_{0} \) is standard normal in \( \R^{\dimp} \) while 
\( P_{1} = \ND(\betav, \BB) \) with 
\( \betav = \Varxi^{1/2} (\bvb - \bvs) \) and \( \BB \eqdef \Varxi^{-1/2} \Varxib \Varxi^{-1/2} \) 
%\( \betav \eqdef \DPc (\thetav_{1} - \thetav_{0}) \)
\begin{EQA}[c]
	2 \log \frac{d\P_{1}}{d\P_{0}}(\gammav)
	=
	\log \det (\BB) 
	- (\gammav - \betav)^{\T} \BB (\gammav - \betav)  
	+ \| \gammav \|^{2}
\end{EQA}
with \( \gammav \) standard normal and
\begin{EQA}
	2 \kullb(\P_{0},\P_{1})
	&=&
	- 2 \E_{0} \log \frac{d\P_{1}}{d\P_{0}}
%	\\
%	&=&
	=
	- \log \det (\BB) 
	+ \tr (\BB - \Id_{\dimp}) 
	+ \betav^{\T} \BB \betav  .
\label{2kullbP0P1}
\end{EQA}
Let \( a_{j} \) be the \( j \)th eigenvalue of \( \BB - \Id_{\dimp} \). 
The condition \( \| \BB - \Id_{\dimp} \|_{\oper} \leq 1/2 \) yields 
\( |a_{j}| \le 1/2 \) and 
\begin{EQA}
	2 \kullb(\P_{0},\P_{1})
	&=&
    \betav^{\T} \BB \betav 
	+
    \sum_{j=1}^{\dimp} \bigl\{ a_{j} - \log(1 + a_{j}) \bigr\}
	\\
	& \leq &
    \betav^{\T} \BB \betav
	+ \sum_{j=1}^{\dimp} a_{j}^{2} 
	\\
	& \leq &
	\betav^{\T} \BB \betav + \tr (\BB - \Id_{\dimp})^{2} 
	\leq 
	\betav^{\T} \BB \betav + \errSi^{2}.
\end{EQA}
This implies by Pinsker's inequality 
\begin{EQA}
	\sup_{A} | \P_{0}(A) - \P_{1}(A) |
	& \leq &
	\sqrt{\frac{1}{2} \kullb(\P_{0},\P_{1})} 
	\leq 
	\frac{1}{2} \sqrt{\errSi^{2} + \betav^{\T} \BB \betav} 
\label{TVdPiBvM}
\end{EQA}
as required.
\end{proof}

Notice that the operator norm bound
\begin{EQA}
	\| \Varxi^{-1/2} \Varxib \Varxi^{-1/2} - \Id_{\dimp} \|_{\oper}
	& \leq &
	\rd
\label{SpNormbound}
\end{EQA}
implies for \( \BB = \Varxi^{-1/2} \Varxib \Varxi^{-1/2} \)
\begin{EQA}
	\tr\bigl( \BB - \Id_{\dimp} \bigr)^{2}
	& \leq &
	\dimp \rd^{2},
	\qquad
	\betav^{\T} \BB \betav \leq (1 + \rd) \| \betav \|^{2} .
\label{SpNormbounds}
\end{EQA}

\begin{corollary}
\label{KullbTVoper}
Let \( \P_{0} = \ND(\bvs,\Varxi) \) and \( \P_{1} = \ND(\bvb,\Varxib) \) for 
some non-degenerated matrices \( \Varxi \) and \( \Varxib \) satisfying \eqref{SpNormbound}.
Then 
\begin{EQA}
	\sup_{A} | \P_{0}(A) - \P_{1}(A) |
	& \leq &
	\frac{1}{2} \sqrt{\dimp \rd^{2} + (1 + \rd) \| \betav \|^{2}} 
\label{TVdPiBvMoper}
\end{EQA}
\end{corollary}

For the special case with \( \betav \equiv 0 \), 
we bound for any Borel set \( A \subset \R^{\Meta} \)
\begin{EQA}
	\bigl| \P\bigl( T(\xiv) \in A\bigr) - \P\bigl(T(\bxiv) \in A\bigr) \bigr|
	& \leq &
	\errSi / 2 \, .
\label{PTxiAbxiA}
\end{EQA}
We state a separate corollary for the distribution of the maximum.

\begin{corollary}
\label{Gausscomp}
Let two \( \dimp \)-dimensional zero mean Gaussian vectors 
\( \xiv \sim \ND(0,\Varxi) \) and \( \bxiv \sim \ND(0,\Varxib) \) be given, and 
\eqref{dSiSic2} holds. % for \( \rd \leq 1/2 \). 
Then for any mapping \( T \colon \R^{\dimp} \to \R^{\Meta} \) and any set of values 
\( (\qq_{\etav}) \), the random vectors 
\( \Xv = T(\xiv) \) and \( \Yv = T(\bxiv) \) fulfill
\begin{EQA}
	\bigl| 
		\P\bigl( \max_{\etav} X_{\etav} - \qq_{\etav} > 0 \bigr) 
		- \P\bigl( \max_{\etav} Y_{\etav} - \qq_{\etav} > 0 \bigr) 
	\bigr|
	& \leq &
	\errSi / 2 .
%	\rd \sqrt{\dimp/2} .
\label{PmaxXYTKL}
\end{EQA}
\end{corollary}

\begin{proof}
We simply apply the result of the lemma to the set 
\( A = \{ \xv \in \R^{\dimp} \colon T(\xv) \leq \zv \} \).
\end{proof}

\ifbook{}{
Interestingly, this method can be used for obtaining an anti-concentration bound in the case 
of a homogeneous mapping \( T : \R^{\dimp} \to \R^{\Meta} \).

\begin{theorem}
\label{Tanticoncxi}
Let \( \xiv \sim \ND(0,\Varxi) \) be a Gaussian vector in \( \R^{\dimp} \).
For any homogeneous mapping \( T : \R^{\dimp} \to \R^{\Meta} \), 
and for any \( \zq > 0 \) and \( \Delta \) satisfying \( 0 \leq \Delta /\zq \leq 1 \), it holds
\begin{EQA}
%	\bigl| 
		\P\bigl( \max_{\etav} T_{\etav}(\xiv) \geq \zq \bigr) 
		- \P\bigl( \max_{\etav} T_{\etav}(\xiv) \geq \zq + \Delta \bigr)
%	\bigr|
	& \leq &
	\Delta \, \zq^{-1} \sqrt{\dimp/2} .
	\qquad
\label{Deqqm12p}
\end{EQA}
Moreover, if \( \bxiv \sim \ND(0,\Varxib) \) is another Gaussian vector and 
\eqref{dSiSic2} holds with \( \rd \leq 1/2 \) and some \( \errSi \geq 0 \), then 
\begin{EQA}
	\bigl| 
		\P\bigl( \max_{\etav} T_{\etav}(\xiv) \geq \zq \bigr) 
		- \P\bigl( \max_{\etav} T_{\etav}(\bxiv) \geq \zq + \Delta \bigr)
	\bigr|
	& \leq &
	\errSi/2 + \Delta \, \zq^{-1} \sqrt{\dimp/2}  .
	\qquad
\label{Deqqm12p2}
\end{EQA}
\end{theorem}

\begin{proof}
Given \( \zq \) and \( \Delta \), define \( \bxiv = {\zq/(\zq + \Delta)} \, \xiv \).
It holds by homogeneity of \( T \)
\begin{EQA}
	\P\bigl( \max_{\etav} T_{\etav}(\xiv) \geq \zq + \Delta \bigr)
	& = &
	\P\bigl( \max_{\etav} T_{\etav}(\bxiv) \geq \zq \bigr) .
\label{PmaxetaTetabxi}
\end{EQA}
It is obvious that \( \Var(\bxiv) = (1 + \Delta/\zq)^{-2} \Varxi \).
%and \eqref{SpNormbound} holds with \( \rd = 1 - (1 + \Delta/\zq)^{-2} \).
Now it holds for the KL-divergence between \( \xiv \) and \( \bxiv \)
\begin{EQA}
	\kullb\bigl( \P_{\xiv},\P_{\bxiv} \bigr)
	&=&
	\frac{\dimp}{2} \bigl\{ 2 \Delta/\zq + (\Delta/\zq)^{2} - 2 \log(1 + \Delta/\zq) \bigr\}
	\leq 
	\dimp (\Delta/\zq)^{2} .
\label{2kullbP0P1}
\end{EQA}
Here we used that \( \log(1+\rho) \leq \rho - \rho^{2}/2 \) for \( \rho \leq 1 \).
Now Pinsker's bound \eqref{TVdPiBvM} implies
\begin{EQA}
	&& 
	\nquad
	\bigl| 
		\P\bigl( \max_{\etav} T_{\etav}(\xiv) \geq \zq \bigr) 
		- \P\bigl( \max_{\etav} T_{\etav}(\bxiv) \geq \zq + \Delta \bigr)
	\bigr|
	\\
	& \leq &
	\P\bigl( \max_{\etav} T_{\etav}(\xiv) \geq \zq \bigr) 
		- \P\bigl( \max_{\etav} T_{\etav}(\xiv) \geq \zq + \Delta \bigr)
	\\
	&&
	+ \, \bigl| \P\bigl( \max_{\etav} T_{\etav}(\xiv) \geq \zq + \Delta \bigr) 
		- \P\bigl( \max_{\etav} T_{\etav}(\bxiv) \geq \zq + \Delta \bigr)
 	\bigr|
	\\
	& \leq &
	\errSi/2 + \Delta \, \zq^{-1} \sqrt{\dimp/2}
\label{errGauscompf}
\end{EQA}
and \eqref{Deqqm12p2} follows.
\end{proof}

We also present a simple corollary of the above result which concerns the change in the expectation
\( \E \fs(\xiv) \) for a bounded function \( \fs \).

\begin{lemma}
\label{Lxiqdexi}
Let \( \xiv \sim \ND(0,\Varxi) \) and \( \bxiv \sim \ND(0,\Varxib) \), where 
\( \Varxi, \Varxib \) satisfy \eqref{dSiSic2}.
For any function \( \fs \) on \( \R^{\dimp} \) with \( |\fs(\xv)| \leq 1 \), 
and any \( \delta > 0 \)
it holds
\begin{EQA}
	\bigl| \E \fs(\xiv) - \E \fs\bigl( \bxiv \bigr) \bigr|
	& \leq &
	\errSi .
\label{fxivfaxiv}
\end{EQA}
Also, for any \( \delta \geq 0 \)
\begin{EQA}
	\bigl| \E \fs(\xiv) - \E \fs\bigl( (1+\delta) \xiv \bigr) \bigr|
	& \leq &
	\delta \sqrt{2\dimp} .
\label{fxivfaxiv}
\end{EQA}
\end{lemma}
\begin{proof}
%Represent \( \xiv = \Varxi^{1/2} \gammav \) where \( \gammav \) is standard normal.
%Define \( g(\xv) = \fs((1 + \delta) \Varxi^{1/2} \xv) \).
%Then \( \fs(\xiv) = \fs(\Varxi^{1/2} \gammav) = g\bigl( (1 + \delta)^{-1}\gammav \bigr) \),
%and similarly \( \fs\bigl( (1+\delta) \xiv \bigr) = g(\gammav) \).
in view of \( |\fs(\xv)| \leq 1 \), it holds
\begin{EQA}
	\bigl| \E \fs(\xiv) - \E \fs( \bxiv) \bigr|
	& \leq &
	\int |\fs(\xv)| \cdot 
		\bigl| \phi_{\xiv}(\xv) - \phi_{\bxiv}(\xv) \bigr| d\xv 
	\leq
	\int \bigl| \phi_{\xiv}(\xv) - \phi_{\bxiv}(\xv) \bigr| d\xv \, .
\label{EfsEfsdexi}
\end{EQA}
%where \( \phi_{\delta}(\xv) = (1 + \delta)^{\dimp/2} \phi\bigl( (1+\delta) \xv \bigr) \)
%is the density of \( (1 + \delta)^{-1} \gammav \).
One more use of Pinsker's inequality yields
\begin{EQA}
	\int \bigl| \phi_{\xiv}(\xv) - \phi_{\bxiv}(\xv) \bigr| d\xv
	& = &
	2 \| \P_{\xiv} - \P_{\bxiv} \|_{TV}
	\leq 
	\sqrt{2\kullb(\P_{\xiv},\P_{\bxiv})} ,
%	\leq 
%	\errSi 
%	{\delta} \sqrt{\dimp} 
\label{phi1del}
\end{EQA}
and the assertion \eqref{fxivfaxiv} follows by \( 2 \kullb(\P_{\xiv},\P_{\bxiv}) \leq \errSi^{2} \).
It remains to note that for \( \Varxib = (1 + \delta)^{2} \Varxi \), it holds 
\( \kullb(\P_{\xiv},\P_{\bxiv}) \leq \delta^{2} \dimp \); see \eqref{2kullbP0P1}.
\end{proof}
}

%%%%%%%%%%%%%%%%%%%%%%%%%%%%%%%%%%%%%%%%%%%%%%%%%%%%%%%%%%%%%%%%
}{ % annals
}

\ifims{ %full
}{ % Annals
\begin{supplement}[id=suppA]
%  \sname{Supplement A}
  \stitle{Some result on matrix valued quadratic forms}
  \slink[doi]{COMPLETED BY THE TYPESETTER}
  \sdatatype{.pdf}
  \sdescription{The supplement collects some useful technical facts 
  on Gaussian quadratic forms and Gaussian matrices}
\end{supplement}
}

\bibliography{exp_ts,listpubm-with-url,phdbib}

\begin{thebibliography}{}

\bibitem[Arlot, 2009]{arlot2009}
Arlot, S. (2009).
\newblock Model selection by resampling penalization.
\newblock {\em Electron. J. Statist.}, 3:557--624.

\bibitem[Barron et~al., 1999]{barronmassartbirge1999}
Barron, A., Birg{\'e}, L., and Massart, P. (1999).
\newblock Risk bounds for model selection via penalization.
\newblock {\em Probab. Theory Related Fields}, 113(3):301--413.

\bibitem[Beran, 1986]{beran1986}
Beran, R. (1986).
\newblock Discussion: Jackknife, bootstrap and other resampling methods in
  regression analysis.
\newblock {\em Ann. Statist.}, 14(4):1295--1298.

\bibitem[Birg{\'e}, 2001]{birgelepski}
Birg{\'e}, L. (2001).
\newblock {\em An alternative point of view on Lepski's method}, volume Volume
  36 of {\em Lecture Notes--Monograph Series}, pages 113--133.
\newblock Institute of Mathematical Statistics, Beachwood, OH.

\bibitem[Birg\'e and Massart, 2007]{massartminimalpenalties}
Birg\'e, L. and Massart, P. (2007).
\newblock Minimal penalties for gaussian model selection.
\newblock {\em Probability Theory and Related Fields}, 138(1-2):33--73.

\bibitem[Cai and Low, 2003]{cailow2003}
Cai, T.~T. and Low, M.~G. (2003).
\newblock A note on nonparametric estimation of linear functionals.
\newblock {\em Ann. Statist.}, 31(4):1140--1153.

\bibitem[Cai and Low, 2005]{cailow2005}
Cai, T.~T. and Low, M.~G. (2005).
\newblock On adaptive estimation of linear functionals.
\newblock {\em Ann. Statist.}, 33(5):2311--2343.

\bibitem[Cavalier et~al., 2002]{cavalier2002}
Cavalier, L., Golubev, G.~K., Picard, D., and Tsybakov, A.~B. (2002).
\newblock Oracle inequalities for inverse problems.
\newblock {\em Ann. Statist.}, 30(3):843--874.

\bibitem[Cavalier and Golubev, 2006]{cavaliergolubev2006}
Cavalier, L. and Golubev, Y. (2006).
\newblock Risk hull method and regularization by projections of ill-posed
  inverse problems.
\newblock {\em Ann. Statist.}, 34(4):1653--1677.

\bibitem[Chernozhukov et~al., 2014]{chernozhukov2014}
Chernozhukov, V., Chetverikov, D., and Kato, K. (2014).
\newblock Anti-concentration and honest, adaptive confidence bands.
\newblock {\em Ann. Statist.}, 42(5):1787--1818.

\bibitem[Dalalyan and Salmon, 2012]{Dala2012}
Dalalyan, A.~S. and Salmon, J. (2012).
\newblock Sharp oracle inequalities for aggregation of affine estimators.
\newblock {\em Ann. Statist.}, 40(4):2327--2355.

\bibitem[Gach et~al., 2013]{GaNiSp2012}
Gach, F., Nickl, R., and Spokoiny, V. (2013).
\newblock {Spatially Adaptive Density Estimation by Localised Haar
  Projections}.
\newblock {\em {Annales de l'Institut Henri Poincare - Probability and
  Statistics}}, 49(3):900--914.
\newblock DOI: 10.1214/12-AIHP485; arXiv:1111.2807.

\bibitem[Gine and Nickl, 2010]{GiNi2010}
Gine, E. and Nickl, R. (2010).
\newblock Confidence bands in density estimation.
\newblock {\em Ann. Statist.}, 38(2):1122--1170.

\bibitem[Goldenshluger, 2009]{goldenshluger2009}
Goldenshluger, A. (2009).
\newblock A universal procedure for aggregating estimators.
\newblock {\em Ann. Statist.}, 37(1):542--568.

\bibitem[Kneip, 1994]{Kneip1994}
Kneip, A. (1994).
\newblock Ordered linear smoothers.
\newblock {\em Ann. Statist.}, 22(2):835--866.

\bibitem[Laurent and Massart, 2000]{laurentmassart2000}
Laurent, B. and Massart, P. (2000).
\newblock Adaptive estimation of a quadratic functional by model selection.
\newblock {\em Ann. Statist.}, 28(5):1302--1338.

\bibitem[Lepski, 1990]{lepski1990}
Lepski, O.~V. (1990).
\newblock A problem of adaptive estimation in {G}aussian white noise.
\newblock {\em Teor. Veroyatnost. i Primenen.}, 35(3):459--470.

\bibitem[Lepski, 1991]{lepski1991}
Lepski, O.~V. (1991).
\newblock Asymptotically minimax adaptive estimation. {I}. {U}pper bounds.
  {O}ptimally adaptive estimates.
\newblock {\em Teor. Veroyatnost. i Primenen.}, 36(4):645--659.

\bibitem[Lepski, 1992]{lepski1992}
Lepski, O.~V. (1992).
\newblock Asymptotically minimax adaptive estimation. {II}. {S}chemes without
  optimal adaptation. {A}daptive estimates.
\newblock {\em Teor. Veroyatnost. i Primenen.}, 37(3):468--481.

\bibitem[Lepski et~al., 1997]{lepski1997}
Lepski, O.~V., Mammen, E., and Spokoiny, V.~G. (1997).
\newblock Optimal spatial adaptation to inhomogeneous smoothness: an approach
  based on kernel estimates with variable bandwidth selectors.
\newblock {\em The Annals of Statistics}, 25(3):929--947.

\bibitem[Lepski and Spokoiny, 1997]{lepskispokoiny1997}
Lepski, O.~V. and Spokoiny, V.~G. (1997).
\newblock Optimal pointwise adaptive methods in nonparametric estimation.
\newblock {\em Ann. Statist.}, 25(6):2512--2546.

\bibitem[{Mammen}, 1993]{mammen1993}
{Mammen}, E. (1993).
\newblock {Bootstrap and wild bootstrap for high dimensional linear models.}
\newblock {\em {Ann. Stat.}}, 21(1):255--285.

\bibitem[Massart, 2007]{massartconcentration2007}
Massart, P. (2007).
\newblock {\em Concentration inequalities and model selection}.
\newblock Number 1896 in Ecole d'Et{\'e} de Probabilit{\'e}s de Saint-Flour.
  Springer.

\bibitem[Spokoiny, 2012]{SP2011}
Spokoiny, V. (2012).
\newblock {Parametric estimation. Finite sample theory}.
\newblock {\em Ann. Statist.}, 40(6):2877--2909.
\newblock arXiv:1111.3029.

\bibitem[Spokoiny and Vial, 2009]{spokoinyvial}
Spokoiny, V. and Vial, C. (2009).
\newblock Parameter tuning in pointwise adaptation using a propagation
  approach.
\newblock {\em Ann. Statist.}, 37:2783--2807.

\bibitem[Spokoiny et~al., 2013]{SpWe2012}
Spokoiny, V., Wang, W., and H{\"a}rdle, W. (2013).
\newblock Local quantile regression (with rejoinder).
\newblock {\em J. of Statistical Planing and Inference}, 143(7):1109--1129.
\newblock ArXiv:1208.5384.

\bibitem[Tropp, 2015]{Tropp2014}
Tropp, J.~A. (2015).
\newblock Found. Trends Mach. Learning.
\newblock to appear.

\bibitem[Tsybakov, 2000]{Tsybakov2000835}
Tsybakov, A. (2000).
\newblock On the best rate of adaptive estimation in some inverse problems.
\newblock {\em Comptes Rendus de l'Academie des Sciences - Series I -
  Mathematics}, 330(9):835 -- 840.

\bibitem[Wu, 1986]{wu1986}
Wu, C. F.~J. (1986).
\newblock Jackknife, bootstrap and other resampling methods in regression
  analysis.
\newblock {\em Ann. Statist.}, 14(4):1261--1295.

\end{thebibliography}
\end{document}